\documentclass[12pt]{article}%
\usepackage{amsmath}
\usepackage{amscd}
\usepackage{amsfonts}
\usepackage{amssymb}%
\newtheorem{theorem}{Theorem}[section]
\newtheorem{lemma}[theorem]{Lemma}
\newtheorem{proposition}[theorem]{Proposition}

\newtheorem{definition}[theorem]{Definition}
\newtheorem{corollary}[theorem]{Corollary}
\newtheorem{remark}[theorem]{Remark}
\newtheorem{remarks}[theorem]{Remarks}
\newtheorem{notation}[theorem]{Notation}
\newtheorem{exremark}[theorem]{Extended Remark and Notation}
\newtheorem{conclremarks}[theorem]{Concluding Remarks}
\newenvironment{proof}{\bf Proof. \rm}{$\Box$}
\newcommand{\be}{\begin{equation}}
\newcommand{\ee}{\end{equation}}

\begin{document}

\title{Canonical Models for Representations of Hardy Algebras}
\author{Paul S. Muhly\thanks{Supported in part by grants from the National Science
Foundation and from the U.S.-Israel Binational Science Foundation.}\\Department of Mathematics\\University of Iowa\\Iowa City, IA 52242\\e-mail: muhly@math.uiowa.edu
\and Baruch Solel\thanks{Supported in part by the U.S.-Israel Binational Science
Foundation and by the Fund for the Promotion of Research at the Technion.}\\Department of Mathematics\\Technion\\32000 Haifa, Israel\\e-mail: mabaruch@techunix.technion.ac.il}
\date{}
\maketitle

\section{Introduction}

Our objective in this paper is to describe a model theory for representations
of the Hardy algebras, which we defined and studied in \cite{MS04}, that
generalizes the model theory of Sz.-Nagy and Foia\c{s} \cite{szNF70} for
contraction operators. Our inspiration for this project comes from three
sources. The first is the well-known fact that model theory allows one to
think of a contraction on Hilbert space as a \textquotedblleft
quotient\textquotedblright\ of a \textquotedblleft
projective\textquotedblright\ module over $H^{\infty}(\mathbb{T})$. More
accurately but still incompletely, one views $H^{\infty}(\mathbb{T})$ as an
operator theoretic generalization of the polynomial algebra in one variable
$\mathbb{C}[X]$ and one thinks of the Hilbert space of the contraction as a
module over the algebra it generates, viewing it as a compression of a\ module
over $H^{\infty}(\mathbb{T})$ that is, essentially, a multiplication
representation of $H^{\infty}(\mathbb{T})$ on a vector-valued $H^{2}$-space.
Indeed, the $H^{\infty}(\mathbb{T})$ - $\mathbb{C}[X]$ analogy coupled with
model theory has inspired much of operator theory during the last 40 years -
and more. We find the \textquotedblleft module-over-$H^{\infty}(\mathbb{T}%
)$\textquotedblright\ perspective particularly stimulating and we have been
especially inspired by the work of Douglas and his collaborators (see, e.g.,
\cite{DP89})\ and by the work of Arveson \cite{wAr69, wAr72, wAr98}.

The second source of inspiration for us is the marvelous paper of Pimsner
\cite{mP97} that shows how to build a $C^{\ast}$-algebra, now called a
Cuntz-Pimsner algebra, from a \textquotedblleft coefficient\textquotedblright%
\ $C^{\ast}$-algebra $A$, say, and a certain type of bimodule $E$ over $A$,
known as a $C^{\ast}$-correspondence. These are denoted $\mathcal{O}(E)$. When
$A=\mathbb{C}$ and $E=\mathbb{C}^{n}$, $\mathcal{O}(E)$ is the famous Cuntz
algebra $\mathcal{O}_{n}$. Sitting inside $\mathcal{O}(E)$ is the norm-closed
subalgebra $\mathcal{T}_{+}(E)$ generated by $A$ and $E$ that we call the
\emph{tensor algebra }of $E$ \cite{MS98}. Indeed, $\mathcal{T}_{+}(E)$ is a
completion of the \emph{algebraic} tensor algebra determined by $A$ and $E$.
For the study of representations of tensor algebras and for other purposes, we
were led to consider certain \textquotedblleft weak closures\textquotedblright%
\ of our correspondences $E$ and to form a \textquotedblleft weak
completion\textquotedblright\ of $\mathcal{T}_{+}(E)$, which we called a
\emph{Hardy algebra} and which we denoted $H^{\infty}(E)$ \cite{MS04}. When
$A=\mathbb{C}=E$, the constructs we are discussing are these: The algebraic
tensor algebra is the polynomial algebra $\mathbb{C}[X]$; the tensor algebra
$\mathcal{T}_{+}(E)$ is the disc algebra $A(\mathbb{D})$ viewed as the algebra
of continuous functions on the circle that extend to be analytic on the open
unit disc; and the Hardy algebra, $H^{\infty}(E)$, is $H^{\infty}(\mathbb{T}%
)$. When $A=\mathbb{C}$ and $E=\mathbb{C}^{n}$, the algebraic tensor algebra
is the free algebra in $n$ variables, $\mathbb{C}\langle X_{1},X_{2}%
,\cdots,X_{n}\rangle$; $\mathcal{T}_{+}(E)$ is Popescu's noncommutative disc
algebra \cite{gP91, gP96}; and $H^{\infty}(E)$ is the free semigroup algebra
that he defined in \cite{gP91} and that has been the object of intense study
by Davidson and Pitts, and others \cite{DP98a, kD01}.

And the third source of inspiration comes from the 1947 paper by Hochschild
\cite{gH47}, which shows, among other things, that \emph{every} finite
dimensional algebra over an algebraically closed field may be expressed as a
quotient of a tensor algebra. In fact, in a fashion that is spelled out in
\cite{pM97}, if one is interested in studying the representation theory of
finite dimensional complex algebras, one may assume that the coefficient
algebra is a commutative $C^{\ast}$-algebra. That is, every finite dimensional
algebra is Morita equivalent to a quotient of a graph algebra. By this we mean
the following: Let $G=(G^{0},G^{1},r,s)$ be a countable graph with vertex
space $G^{0}$, edge space $G^{1}$ and range and source maps $r$ and $s$. Then
for the $C^{\ast}$-algebra $A$ we take $c_{0}(G^{0})$ and for $E$ we take (a
completion of) the space of finitely supported functions $\xi$ on $G^{1}$,
which may be view as a bimodule over $A$ via the formula: $a\xi b(\alpha
):=a(r(\alpha))\xi(\alpha)b(s(\alpha))$, $a,b\in A$ and $\alpha\in G^{1}$. If
the graph is finite, then the algebraic tensor algebra is the type of algebra
to which we just referred. \emph{Every} finite dimensional algebra over
$\mathbb{C}$ is naturally Morita equivalent to a quotient of such a tensor
algebra. This perspective has dominated much of finite dimensional algebra
since Gabriel's penetrating study \cite{pG72} of algebras of finite
representation type. (For a recent survey, see \cite{pG92}.) In general, the
Cuntz-Pimsner algebra $\mathcal{O}(E)$ in this setting goes under various
names, depending on the structure of the graph, but for the sake of this
discussion, $\mathcal{O}(E)$ is simply a Cuntz-Krieger algebra first studied
in \cite{CK80}.\ The tensor algebra $\mathcal{T}_{+}(E)$ has been studied by
us in \cite{pM97, MS98, MS99, MS00}. The general theory of Hardy algebras that
we developed in \cite{MS04} was initiated in part to study $H^{\infty}(E)$ in
this setting, and special representations of $H^{\infty}(E)$, when $E$ comes
from a graph, have been studied by Kribs and Power and their co-workers under
the name \textquotedblleft free semi-groupoid algebras\textquotedblright. (See
\cite{KPp03}.)

The three sources of inspiration combined have become the driving force behind
much of our recent work: We want to study tensor algebras and Hardy algebras
in a fashion analogous to the theory of contraction operators on Hilbert space
with an eye to exploiting the insights from finite dimensional algebra in much
the same way that finite dimensional matrix theory and linear algebra inform
operator theory. Although our initial focus was on the interactions between
operator theory and finite dimensional algebra, we soon realized that the
perspective provided significant insights into such things as the theory of
(irreversible) dynamical systems \cite{MS98a, MS00}, the theory of completely
positive maps, quantum Markov processes and other aspects of quantum
probability\cite{MS02, MS03}. Of course, we are not alone in the appreciation
of the impact of Pimsner's insights on these subjects. However, the
perspective from non-self-adjoint operator theory and algebras that has been
the leitmotif of our work led to useful insights that seem not to be easily
accessible from the self adjoint theory.

The theory we present here will be seen to be a direct descendant of the
Sz.-Nagy-Foia\c{s} theory spelled out in \cite{szNF70}. However, there is a
subtle, yet important, distinction. We present a model theory for \emph{some
}representations of our Hardy algebras, not all. We run into the same
difficulties that Popescu encountered in \cite{gP89} and we must limit
ourselves to what he called completely non-coisometric representations. We
adopt his terminology here. Indeed, our analysis owes a great deal to his work.

In the next section we present background information from \cite{MS04} and
elsewhere that we shall use. In particular, we develop the perspective that
the elements in one of our Hardy algebras $H^{\infty}(E)$ can profitably be
studied as functions on the unit ball of the so-called dual of $E$ calculated
with respect to a faithful representation of the underlying $W^{\ast}%
$-algebra. In Section 3, we develop the notion of characteristic operators and
functions for completely non-coisometric representations of $H^{\infty}(E)$
and we show that such representations have canonical models that are (almost)
the exact analogue of the models that Sz.-Nagy and Foia\c{s} built for single
operators. In Section 4, we prove a model-theoretic analogue of Sarason's
original commutant lifting theorem \cite{dS67} and in Section 5 we identify
the relation between invariant subspaces for representations and
factorizations of the characteristic functions. Finally, in Section 6, we
present an example that shows how our theory functions in a special case
related to the classical Sz.-Nagy-Foia\c{s} theory and that helps to clarify
the limitations of the \textquotedblleft completely non
coisometric\textquotedblright\ hypothesis.

\section{Preliminaries\label{preliminaries}}

\subsection{$W^{\ast}$-Correspondences and Hardy Algebras}

We begin by recalling the notion of a $W^{\ast}$-\emph{correspondence}. For
the general theory of Hilbert $C^{\ast}$-modules which we use, we will follow
\cite{cL95}. In particular, a Hilbert $C^{\ast}$-module will be a \emph{right}
Hilbert $C^{\ast}$-module.

\begin{definition}
\label{correspondence}Let $M$ and $N$ be $W^{\ast}$-algebras and let $E$ be a
(right) Hilbert $C^{\ast}$-module over $N$. Then $E$ is called a
\emph{(Hilbert) $W^{\ast}$-module} over $N$ in case it is self dual (i.e.
every continuous $N$-module map from $E$ to $N$ is implemented by an element
of $E$). It is called a \emph{$W^{\ast}$-correspondence} from $M$ to $N$ if it
is also endowed with a structure of a left $M$-module via a normal $\ast
$-homomorphism $\varphi:M\rightarrow\mathcal{L}(E)$.(Here $\mathcal{L}(E)$ is
the algebra of all bounded, adjointable, module maps on $E$ - which is a
$W^{\ast}$-algebra when $E$ is a $W^{\ast}$-module \cite{wP73}). A
\emph{\ $W^{\ast}$-correspondence over $M$} is simply a $W^{\ast}%
$-correspondence from $M$ to $M$.

An isomorphism of $W^{\ast}$-correspondences $E_{1},E_{2}$ from $M$ to $N$ is
an $M,N$-linear, surjective, bimodule map that preserves the inner product. We
shall write $E_{1}\cong E_{2}$ if such an isomorphism exists.
\end{definition}

If $E$ is a $W^{\ast}$-correspondence from $M$ to $N$ and if $F$ is a
$W^{\ast}$-correspondence from $N$ to $Q$, then the balanced tensor product,
$E\otimes_{N}F$ is a $W^{\ast}$-correspondence from $M$ to $Q$. It is defined
as the self-dual extension \cite{wP73} of the Hausdorff completion of the
algebraic balanced tensor product with the internal inner product given by
\[
\langle\xi_{1}\otimes\eta_{1},\xi_{2}\otimes\eta_{2}\rangle=\langle\eta
_{1},\varphi(\langle\xi_{1},\xi_{2}\rangle_{E})\eta_{2}\rangle_{F}%
\]
for all $\xi_{1}$ , $\xi_{2}$ in $E$ and $\eta_{1}$ , $\eta_{2}$ in $F$. The
left and right actions of $M$ and $Q$ are defined by
\[
\varphi_{E\otimes_{N}F}(a)(\xi\otimes\eta)b=\varphi_{E}(a)\xi\otimes\eta b
\]
for $a$ in $M$, $b$ in $Q$, $\xi$ in $E$ and $\eta$ in $F$.

If $\sigma$ is a normal representation of $N$ on a Hilbert space $H$ and $E$
is a $W^{\ast}$-correspondence from $M$ to $N$, then $H$ can be viewed as a
$W^{\ast}$-correspondence from $N$ to $\mathbb{C}$ and $E\otimes_{N}H$ is a
Hilbert space (with a normal representation of $M$ on it). \ Of course,
$E\otimes_{N}H$, also denoted $E\otimes_{\sigma}H$, is nothing but the Hilbert
space of the representation of $M$ that is \emph{induced }by $\sigma$,
$E$-$Ind_{N}^{M}\sigma$, in the sense of Rieffel's pioneering studies
\cite{mR74a, mR74b}. (See \cite[p. 36 ff.]{RW98} for the general theory.) It
is defined by the equation
\[
E\text{-}Ind_{N}^{M}\sigma(a)(\xi\otimes h)=a\xi\otimes h,\;\;\ \xi\otimes
h\in E\otimes_{\sigma}H,\ a\in M\text{.}%
\]
To lighten the formulas that appear in this paper, we adopt the following
notation throughout.

\begin{notation}
\label{Indnotation}If $E$ is a Hilbert $W^{\ast}$-module over a von Neumann
algebra $N$, if $\sigma$ a normal representation of $N$ on the Hilbert space
$H$ and if $\mathcal{A}$ is any subalgebra of $\mathcal{L}(E)$, then we shall
write $\sigma^{E}$ for the restriction of $E$-$Ind_{N}^{\mathcal{L}(E)}\sigma$
to $\mathcal{A}$, and for $a\in\mathcal{A}$, we shall often abbreviate
$\sigma^{E}(a)$ as $a\otimes I_{H}$.
\end{notation}

Note also that, given an operator $R\in\sigma(M)^{\prime}$, the map that maps
$\xi\otimes h$ in $E\otimes_{\sigma}H$ to $\xi\otimes Rh$ is a bounded linear
operator and we write $I_{E}\otimes R$ for it. \ In fact, Theorem 6.23 of
\cite{mR74a} shows that the commutant of $\sigma^{E}(\mathcal{L}(E))$ is
$\{I_{E}\otimes R\mid R\in\sigma(M)^{\prime}\}$.

If $\{E_{\alpha}\}$ is a family of $W^{\ast}$-correspondences from $M$ to $N$
then one defines the direct sum $\oplus E_{\alpha}$ as in \cite{wP73}. It is a
$W^{\ast}$-module over $N$ and one defines a left module structure (making it
a $W^{\ast}$-correspondence) in a natural way. Combining this observation
about direct sums with the notion of tensor products leads us to the Fock
space construction: Given a $W^{\ast}$-correspondence $E$ over $M$, the
\emph{full Fock space} over $E$, $\mathcal{F}(E)$, is defined to be $M\oplus
E\oplus\ E^{\otimes2}\oplus\cdots$. It is also a $W^{\ast}$-correspondence
over $M$ with the left action $\varphi_{\infty}$ (or $\varphi_{E,\infty}$)
given by the formula
\[
\varphi_{\infty}(a)=diag(a,\varphi(a),\varphi^{(2)}(a),\cdots),
\]
where $\varphi^{(n)}(a)(\xi_{1}\otimes\xi_{2}\otimes\cdots\xi_{n}%
)=(\varphi(a)\xi_{1})\otimes\xi_{2}\otimes\cdots\xi_{n}$ . For $\;\xi\in E\;$
we write $T_{\xi}$ for the creation operator on $\mathcal{F}(E)$ : $T_{\xi
}\eta=\xi\otimes\eta,\;\eta\in\mathcal{F}(E)$.$\;$Then $T_{\xi}$ is a
continuous, adjointable operator in $\mathcal{L}(\mathcal{F}(E))$. The
\emph{norm closed} subalgebra of $\mathcal{L}(\mathcal{F}(E))$ generated by
all the $T_{\xi}$'s and $\varphi_{\infty}(A)$ is called \emph{the tensor
algebra} of $E$ and is denoted $\mathcal{T}_{+}(E)$ (\cite{MS98}). Since
$\mathcal{F}(E)$ is a Hilbert $W^{\ast}$-module, $\mathcal{L}(\mathcal{F}(E))$
is a $W^{\ast}$-algebra \cite{wP73}. Hence the following definition from
\cite{MS04} makes sense.

\begin{definition}
If $E$ is a $W^{\ast}$-correspondence over a $W^{\ast}$-algebra then closure
of $\mathcal{T}_{+}(E)$ in the $w^{\ast}$-topology on $\mathcal{L}%
(\mathcal{F}(E))$ is called the \emph{Hardy algebra of }$E$, and is denoted
$H^{\infty}(E)$.
\end{definition}

The $w^{\ast}$-continuous, completely contractive representations of this
algebra are our principal objects of study.

\subsection{Representations}

Recall that a $W^{\ast}$-correspondence $E$ over a $W^{\ast}$-algebra $M$
carries a natural weak topology, called the $\sigma$-topology (see
\cite{BDH88}). This the topology defined by the functionals$\;f(\cdot
)=\sum_{n=1}^{\infty}\omega_{n}(\langle\eta_{n},\cdot\rangle)$, where the
$\eta_{n}$ lie in $E$ , the $\omega_{n}$ lie in the pre-dual of $M$, $M_{\ast
}$, and where $\sum\Vert\omega_{n}\Vert\Vert\eta_{n}\Vert<\infty$.

\begin{definition}
Let $E$ be a $W^{\ast}$-correspondence over a $W^{\ast}$-algebra $N$ and let
$H$ be a Hilbert space.

\begin{enumerate}
\item[(1)] A \emph{completely contractive covariant representation} of $E$
(or, simply, a \emph{representation} of $E$) in $B(H)$ is a pair, $(T,\sigma
)$, such that

\begin{enumerate}
\item[(a)] $\sigma$ is a normal representation of $N$ in $B(H)$.

\item[(b)] $T$ is a linear, completely contractive map from $E$ to $B(H)$ that
is continuous with respect to the $\sigma$-topology of \cite{BDH88} on $E$ and
the $\sigma$-weak topology on $B(H)$, and

\item[(c)] $T$ is a bimodule map in the sense that $T(\varphi(a)\xi
b)=\sigma(a)T(\xi)\sigma(b),$ $\xi\in E$, and $a,b\;\in\;N$.
\end{enumerate}

\item[(2)] A completely contractive covariant representation $(T,\sigma)$ of
$E$ in $B(H)$ is called \emph{isometric} in case
\[
T(\xi)^{\ast}T(\eta)=\sigma(\langle\xi,\eta\rangle),
\]
for all $\xi,\eta$ in $E$.
\end{enumerate}
\end{definition}

The theory developed in \cite{MS98} applies here to prove that if a
representation $(T,\sigma)$ of $E$ is given, then it determines a contraction
$\tilde{T}:E\otimes_{\sigma}H\rightarrow H$ defined by the formula
\[
\tilde{T}(\xi\otimes h)=T(\xi)h.
\]
Moreover, for every $a$ in $N$ we have
\begin{equation}
\tilde{T}(\varphi(a)\otimes I)=\tilde{T}\sigma^{E}(\varphi(a))=\sigma
(a)\tilde{T}, \label{Ttilde}%
\end{equation}
i.e., $\tilde{T}$ intertwines $\sigma$ and $\sigma^{E}\circ\varphi$. In fact,
it is shown in \cite{MS98} that there is a bijection between representations
$(T,\sigma)$ of $E$ and intertwining operators $\tilde{T}$ of $\sigma$ and
$\sigma^{E}\circ\varphi$.

It is also shown in \cite{MS98} that $(T,\sigma)$ is isometric if and only if
$\tilde{T}$ is an isometry.

\begin{remark}
\label{GenPowers}In addition to $\tilde{T}$ we also require the
\textquotedblleft generalized higher powers\textquotedblright\emph{\ }of
$\tilde{T}$. These are\emph{\ }maps$\;\tilde{T}_{n}:E^{\otimes n}\otimes
H\rightarrow H\;$defined by the equation$\;\tilde{T}_{n}(\xi_{1}\otimes
\ldots\otimes\xi_{n}\otimes h)=T(\xi_{1})\cdots T(\xi_{n})h$, $\xi_{1}%
\otimes\ldots\otimes\xi_{n}\otimes h\in E^{\otimes n}\otimes H$. We call
$\tilde{T}_{n}$ the $n^{th}$\emph{-power} or the $n^{th}$\emph{-generalized
power} of $\tilde{T}.$ An important role in our analysis is played by the
following formula which is valid for all positive integers $m$ and
$n$:$\;\tilde{T}_{n+m}=\tilde{T}_{n}(I_{n}\otimes\tilde{T}_{m})=\tilde{T}%
_{m}(I_{m}\otimes\tilde{T}_{n})$, where $I_{n}$ is the identity map on
$E^{\otimes n}$ \cite{MS99}. It will also be convenient to write $T_{n}%
(\xi)=T(\xi_{1})\cdots T(\xi_{n})$ for $\xi=\xi_{1}\otimes\cdots\otimes\xi
_{n}\in E^{\otimes n}$, so that $\tilde{T}_{n}(\xi\otimes h)=T_{n}(\xi
)h=T(\xi_{1})\cdots T(\xi_{n})h$ for $h\in H$.
\end{remark}

The theory developed in \cite{MS98} shows that there is a bijective
correspondence between covariant representations of $E$ and completely
contractive representations $\rho$ of $\mathcal{T}_{+}(E)$ with the property
that $\rho\circ\varphi_{\infty}$ is a normal representation of $N$. (Given
$\rho$, let $T(\xi):=\rho(T_{\xi})$ and let $\sigma(\cdot)=\rho(\varphi
_{\infty}(\cdot))$ then $(T,\sigma)$ is a representation of $E$ and we write
$\rho:=T\times\sigma$.) However, \emph{only certain of these extend from
}$\mathcal{T}_{+}(E)$ \emph{to }$H^{\infty}(E)$. The full story has yet to be
understood, but an initial analysis may be found in \cite{MS04}. Aspects of
the analysis in \cite{MS04} will play a role in this paper.

The representations of $H^{\infty}(E)$ that are \textquotedblleft
induced\textquotedblright\emph{\ }by representations of $M$ play a central
role in our theory, where they serve as analogues of \emph{pure} isometries.
This is made clear in \cite{MS99} and \cite{MS04} and will be developed
further here.

\begin{definition}
\label{inducedRep}Let $E$ be a correspondence over a $W^{\ast}$-algebra $M$
and let $\sigma_{0}$ be a (normal) representation of $M$ on a Hilbert space
$H$. The representation of $H^{\infty}(E)$ on $\mathcal{F}(E)\otimes
_{\sigma_{0}}H$ \emph{induced }by $\sigma_{0}$ is defined to be the
restriction to $H^{\infty}(E)$ of $\sigma_{0}^{\mathcal{F}(E)}$.
\end{definition}

Observe that the covariant representation $(T,\sigma)$ determined by
$\sigma_{0}^{\mathcal{F}(E)}$ is given by the formulae%
\begin{equation}
\sigma=\sigma_{0}^{\mathcal{F}(E)}\circ\varphi_{\infty}=\varphi_{\infty
}\otimes I_{H} \label{ind1}%
\end{equation}
and%
\begin{equation}
T(\xi)=\sigma_{0}^{\mathcal{F}(E)}(T_{\xi})=T_{\xi}\otimes I_{H}\text{,}
\label{ind2}%
\end{equation}
$\xi\in E$. We also say that $(T,\sigma)$ is \emph{induced by }$\sigma_{0} $.

\begin{remark}
\label{faithful}It follows from Theorem 6.23 of \cite{mR74a} that
$\sigma^{\mathcal{F}(E)}$ is a faithful representation of $H^{\infty}(E)$ if
$\sigma$ is a faithful representation of $M$. Most of the time, we will be
dealing with faithful representations of $M$, and when non-faithful
representations may arise we will go to great lengths to supplement them to
yield faithful representations. (See Definition \ref{Supplemental} and the
discussion related to it.)
\end{remark}

\subsection{Duals and Commutants}

In order to identify the commutant of an induced representation, we introduced
concept of \textquotedblleft duality\textquotedblright\ for correspondences in
\cite{MS04}. Since it plays an important role in the present investigation, we
outline its salient features. Given a $W^{\ast}$-correspondence $E$ over the
$W^{\ast}$-algebra $M$ and given a \emph{faithful} normal representation
$\sigma$ of $M$ on a Hilbert space $H$, we set
\[
E^{\sigma}=\{\eta\in B(H,E\otimes_{\sigma}H)\mid\eta\sigma(a)=(\varphi
(a)\otimes I_{H})\eta,\;a\in M\}.
\]
Then $E^{\sigma}$ is a bimodule over $\sigma(M)^{\prime}$ where the right
action is defined by $\eta S=\eta\circ S$ and the left action by $S\cdot
\eta=(I_{E}\otimes S)\circ\eta$, for $\eta\in E^{\sigma}$ and $S\in
\sigma(M)^{\prime}$. In fact, $E^{\sigma}$ is a $W^{\ast}$-correspondence over
$\sigma(M)^{\prime}$, where the inner product is defined by the formula
$\langle\eta_{1},\eta_{2}\rangle=\eta_{1}^{\ast}\eta_{2}$. This correspondence
is called the \emph{$\sigma$-dual (correspondence) of $E$}. Write $\iota$ for
the identity representation of $\sigma(M)^{\prime}$ on $H$. Then we may form
the $W^{\ast}$-correspondence $(E^{\sigma})^{\iota}$ over $\sigma
(M)^{\prime\prime}=\sigma(M)$. Since $\sigma$ is faithful we can view this as
a correspondence over $M$. As we shall outline, $(E^{\sigma})^{\iota}$ is
naturally isomorphic to $E$ in a way that sets up an isomorphism between the
commutant of the representation of $H^{\infty}(E)$ induced by $\sigma$ and the
image of $H^{\infty}(E^{\sigma})$ under the representation induced by $\iota$.
The latter acts on $\mathcal{F}(E^{\sigma})\otimes_{\iota}H$.

For a given $\xi\in E$ we define the operator $L_{\xi}:H\rightarrow
E\otimes_{\sigma}H$ by the equation $L_{\xi}h=\xi\otimes h$. It is evident
that $L_{\xi}$ is a bounded operator and that its adjoint is given by the
formula $L_{\xi}^{\ast}(\zeta\otimes h)=\sigma(\langle\xi,\zeta\rangle)h$ for
$\zeta\in E$ and $h\in H$.

\begin{proposition}
\label{dual}\hfill

\begin{enumerate}
\item[(i)] \cite[Theorem 3.6]{MS04}For every $\xi\in E$ let $\hat{\xi
}:H\rightarrow E^{\sigma}\otimes_{\iota}H$ be defined by adjoint equation,
\[
\hat{\xi}^{\ast}(\eta\otimes h)=L_{\xi}^{\ast}(\eta(h))\in H,
\]
$\eta\otimes h\in E^{\sigma}\otimes H$. Then $\hat{\xi}\in(E^{\sigma})^{\iota
}$ and the map $\xi\mapsto\hat{\xi}$ is an isomorphism of $W^{\ast}%
$-correspondences (that is, it is a bimodule map and an isometry).

\item[(ii)] \cite[Lemma 3.7]{MS04}For two $W^{\ast}$-correspondences $E_{1}$
and $E_{2}$ over $M$,
\[
(E_{1}\oplus E_{2})^{\sigma}\cong E_{1}^{\sigma}\oplus E_{2}^{\sigma}%
\]
and
\[
(E_{1}\otimes_{M}E_{2})^{\sigma}\cong E_{2}^{\sigma}\otimes_{\sigma
(M)^{\prime}}E_{1}^{\sigma}.
\]
The second isomorphism is given by the map that sends $\eta_{2}\otimes\eta
_{1}\in E_{2}^{\sigma}\otimes_{\sigma(M)^{\prime}}E_{1}^{\sigma}$ to
$(I_{E_{1}}\otimes\eta_{2})\eta_{1}\in(E_{1}\otimes_{M}E_{2})^{\sigma}$.
\end{enumerate}
\end{proposition}

Concerning part (i) of Proposition \ref{dual}, it should be noted that since
$\eta\in E^{\sigma}$, $\eta$ is an operator from $H$ to $E\otimes_{\sigma}H$.
Thus $\eta(h)\in E\otimes_{\sigma}H$ for all $h\in H\ $and $L_{\xi}^{\ast
}(\eta(h))$ makes good sense as an element of $H$.

With the notation we have established, we also have

\begin{proposition}
\label{U}In the notation of Proposition \ref{dual}, the formula
\[
U_{k}(\xi\otimes h)=\hat{\xi}(h)\text{,}%
\]
$\xi\in E^{\otimes k}$, $h\in H$, defines a Hilbert space isomorphism $U_{k} $
from $E^{\otimes k}\otimes_{\sigma}H$ onto $(E^{\sigma})^{\otimes k}%
\otimes_{\iota}H$. The inverse is given by the formula $U^{\ast}(\eta\otimes
h)=\eta(h)$, $\eta\otimes h\in(E^{\sigma})^{\otimes k}\otimes_{\iota}H$. The
direct sum of the $U_{k}$, $U:=\sum_{k\geq0}^{\oplus}U_{k} $, is a Hilbert
space isomorphism from $\mathcal{F}(E)\otimes_{\sigma}H$ onto $\mathcal{F}%
(E^{\sigma})\otimes_{\iota}H$.
\end{proposition}

The following result, \cite[Theorem 3.9]{MS04}, identifies the commutant of an
induced representation in the fashion promised. The theorem is an analogue of
the assertion that the commutant of the unilateral shift is the weakly closed
algebra generated by the unilateral shift. In Theorem~\ref{commutant} it will
be generalized to the \textquotedblleft model-theoretic\textquotedblright%
\ version of the commutant lifting theorem proved by Sarason \cite{dS67}.

\begin{theorem}
\label{comminduced}Let $E$ be a correspondence over the $W^{\ast}$-algebra $M$
and let $\sigma:M\rightarrow B(H)$ be a faithful normal representation of $M$
on the Hilbert space $H$. Write $\sigma^{\mathcal{F}(E)}$ for the
representation of $H^{\infty}(E)$ on $\mathcal{F}(E)\otimes_{\sigma}H$ induced
by $\sigma$, write $\iota^{\mathcal{F}(E^{\sigma})}$ for the representation of
$H^{\infty}(E^{\sigma})$ on $\mathcal{F}(E^{\sigma})\otimes_{\iota}H$ induced
by the identity representation $\iota$ of $\sigma(M)^{\prime}$ on $H$ and
write $U:\mathcal{F}(E)\otimes_{\sigma}H\rightarrow\mathcal{F}(E^{\sigma
})\otimes_{\iota}H$ for the Hilbert space isomorphism described in Proposition
\ref{U}. Then the commutant of $\sigma^{\mathcal{F}(E)}(H^{\infty}(E))$ is
$U^{\ast}\iota^{\mathcal{F}(E^{\sigma})}(H^{\infty}(E^{\sigma}))U$.
\end{theorem}

\begin{exremark}
\label{FunctionalRep}One of the principal achievements of \cite{MS04} was the
representation of elements of $H^{\infty}(E)$ as functions on the open unit
ball of $E$. This representation plays a role here, but \emph{with a twist}.
To understand what we need in more detail, assume that $\sigma$ is a faithful
representation of $M$ in $B(H)$ and let $\eta$ be an operator in the open unit
ball of $E^{\sigma}$, then $\eta^{\ast}:E\otimes_{\sigma}H\rightarrow H$
intertwines $\varphi(a)\otimes I_{H}$ and $\sigma(a)$ for every $a\in M$.
Thus, there is a representation $(T,\sigma)$ of $E$ such that $\eta^{\ast
}=\tilde{T}$ \cite[Lemma 2.16]{MS98}. Since $\Vert\tilde{T}\Vert<1$ the
representation $T\times\sigma$ of $\mathcal{T}_{+}(E)$ on $H $ can be extended
to a $\sigma$-weakly continuous representation, also written $T\times\sigma$,
of $H^{\infty}(E)$ (see \cite[Corollary 2.14]{MS04}). So, given $X\in
H^{\infty}(E)$, we define
\[
X(\eta)=(T\times\sigma)(X)\in B(H).
\]
That is, each $X\in H^{\infty}(E)$ gives a $B(H)$-function defined on the open
unit ball of $E^{\sigma}$. The properties of this functional representation of
$H^{\infty}(E)$ are explored in \cite{MS04} . We point out, however, that in
general the functional representation of $H^{\infty}(E)$ is not faithful. That
is, $X(\eta)$ can vanish for all $\eta$ in the open unit ball of $E^{\sigma}$
without $X=0$ \cite{MS04}. Nevertheless, the function theoretic point of view
proves very effective for studying and unifying a wide variety of problems in
operator theory. In particular, in \cite{MS04}, we proved a general version of
the Nevanlinna-Pick interpolation theorem, which contains an enormous number
of operator theoretic variants of the classical result as a special cases.

In this paper, we shall use the identification of $E$ with $(E^{\sigma
})^{\iota}$ through the map $\xi\mapsto\hat{\xi}$ in part (i) of Proposition
\ref{dual} to view elements of $H^{\infty}(E^{\sigma})$ as functions on the
open unit ball of $E$. More importantly, we shall use the spatial
identification of the commutant of $\sigma^{\mathcal{F}(E)}(H^{\infty}(E))$
with $\iota^{\mathcal{F}(E^{\sigma})}(H^{\infty}(E^{\sigma}))$, given in terms
of $U$ and described in Theorem \ref{comminduced}, to view elements in
$(\sigma^{\mathcal{F}(E)}(H^{\infty}(E)))^{\prime}$ as functions on the open
unit ball of $E$.

Thus, we adopt the following notation: If $\Psi\in(\sigma^{\mathcal{F}%
(E)}(H^{\infty}(E)))^{\prime}$, then $\hat{\Psi}$ will denote the element in
$H^{\infty}(E^{\sigma})$ defined by the formula%
\begin{equation}
\hat{\Psi}:=(\iota^{\mathcal{F}(E^{\sigma})})^{-1}(U\Psi U^{\ast})\text{,}
\label{hat}%
\end{equation}
where $U:\mathcal{F}(E)\otimes_{\sigma}H\rightarrow\mathcal{F}(E^{\sigma
})\otimes_{\iota}H$ is the Hilbert space isomorphism defined in Proposition
\ref{U}. Note that $\iota^{\mathcal{F}(E^{\sigma})}$ is faithful since $\iota$
is (Remark \ref{faithful}). We shall also write equation (\ref{hat}) as%
\begin{equation}
\hat{\Psi}\otimes I_{H}=U\Psi U^{\ast}\text{.} \label{hatprime}%
\end{equation}
We shall then want to evaluate $\hat{\Psi}$ on the open unit ball of $E$. On
the other hand, given an element $\Xi\in H^{\infty}(E^{\sigma})$, we shall
write $\check{\Xi}$ for the operator in the commutant of $\sigma
^{\mathcal{F}(E)}(H^{\infty}(E))$ given by the formula%
\begin{equation}
\check{\Xi}:=U^{\ast}\iota^{\mathcal{F}(E^{\sigma})}(\Xi)U=U^{\ast}(\Xi\otimes
I_{H})U\text{.} \label{check}%
\end{equation}
Thus, evidently, we have $(\hat{\Psi}\check{)}=\Psi$ and $(\check{\Xi}\hat
{)}=\Xi$.

This notation is, of course, suggestive of the idea that the Hilbert space
isomorphism $U$ in Proposition \ref{U} should be viewed as some sort of
generalized Fourier transform. The analogy turns out to be more than one built
from notation. Accordingly, we shall call $U:\mathcal{F}(E)\otimes_{\sigma
}H\rightarrow\mathcal{F}(E^{\sigma})\otimes_{\iota}H$ the \emph{Fourier
transform determined by }$\sigma$. Also, given $\Psi\in(\sigma^{\mathcal{F}%
(E)}(H^{\infty}(E)))^{\prime}$, we shall $\hat{\Psi}$ the \emph{Fourier
transform} of $\Psi$, if $\Xi\in H^{\infty}(E^{\sigma})$, then $\check{\Xi}$
will be called the \emph{inverse Fourier transform }of $\Xi$.
\end{exremark}

\section{Characteristic Operators and Characteristic Functions of
Representations\label{COCFR}}

In the model theory for a single contraction operator on Hilbert space, the
role of the characteristic operator function is to \textquotedblleft
locate\textquotedblright\ the Hilbert space of the operator in the Hilbert
space of its minimal isometric dilation. In \cite{MS98} we successfully
constructed isometric dilations of representations of $H^{\infty}(E)$.
(Actually, in \cite{MS98} we worked with $C^{\ast}$-correspondences over
$C^{\ast}$-algebras. Adjustments necessary to handle representations of
$H^{\infty}(E)$, when $E$ is a $W^{\ast}$-correspondence, were made in
\cite{MS04}.) We therefore begin by briefly recapping aspects of the theory we
shall use.

\subsection{Isometric Dilations}

Let $E$ be a $W^{\ast}$-correspondence over a $W^{\ast}$-algebra $M$ and let
$(T,\sigma)$ be a completely contractive covariant representation of $E$ on a
Hilbert space $H$. Then $(T,\sigma)$ has a \textquotedblleft minimal isometric
dilation\textquotedblright, $(V,\rho)$, defined as follows. Recall that the
map $\tilde{T}:E\otimes_{\sigma}H\rightarrow H$ defined by the equation
$\tilde{T}(\xi\otimes h)=T(\xi)h$ is a contraction that satisfies the equation
$\tilde{T}(\varphi(a)\otimes I_{H})=\sigma(a)\tilde{T}$. We set $\Delta
:=(I-\tilde{T}^{\ast}\tilde{T})^{1/2}$ (in $B(E\otimes_{\sigma}H)$),
$\Delta_{\ast}:=(I-\tilde{T}\tilde{T}^{\ast})^{1/2}$ (in $B(H)$),
$\mathcal{D}:=\overline{\Delta(E\otimes_{\sigma}H)}$ and $\mathcal{D}_{\ast
}:=\overline{\Delta_{\ast}(H)}$. Observe that on account of the intertwining
equation $\tilde{T}(\varphi(a)\otimes I_{H})=\sigma(a)\tilde{T}$,
$\mathcal{D}_{\ast}$ reduces $\sigma$, while $\mathcal{D}$ reduces
$\varphi(\cdot)\otimes I_{H}=\sigma^{E}\circ\varphi(\cdot)$. Also we write
$D(\xi):=\Delta\circ L_{\xi}:H\rightarrow E\otimes_{\sigma}H$, $\xi\in E$,
where, recall, $L_{\xi}:H\rightarrow E\otimes_{\sigma}H$ is the map $L_{\xi
}h=\xi\otimes h$, $h\in H$, $\xi\in E$. Note, too, that $T(\xi)=\tilde{T}\circ
L_{\xi}$.

The representation space $K$ of $(V,\rho)$ is
\begin{align*}
K  &  =H\oplus\mathcal{D}\oplus(E\otimes_{\sigma_{1}}\mathcal{D}%
)\oplus(E^{\otimes2}\otimes_{\sigma_{1}}\mathcal{D})\oplus...\\
&  =H\oplus\mathcal{F}(E)\otimes_{\sigma_{1}}\mathcal{D}%
\end{align*}
where $\sigma_{1}$ is defined to be the restriction to $\mathcal{D}$ of
$\varphi(\cdot)\otimes I_{H}$. The representation $\rho$, in the isometric
dilation $(V,\rho)$ for $(T,\sigma)$, is defined to be $\rho=\sigma
\oplus\sigma_{1}^{\mathcal{F}(E)}\circ\varphi_{\infty}$. That is,
$\rho=diag(\sigma,\sigma_{1},\sigma_{2},\ldots)$ where $\sigma_{k+1}%
(\cdot)=\sigma_{1}^{E^{\otimes k}}\circ\varphi_{k}(\cdot)=\varphi_{k}%
(\cdot)\otimes I_{\mathcal{D}}$ acting on $E^{\otimes k}\otimes_{\sigma_{1}%
}\mathcal{D}$. The map $V:E\rightarrow B(K)$ is defined in terms of the
matrix
\begin{equation}
\label{vxi}V(\xi)=\left(
\begin{array}
[c]{clll}%
T(\xi) & 0 & 0 & \cdots\\
D(\xi) & 0 & 0 & \cdots\\
0 & L_{\xi} & 0 & \\
0 & 0 & L_{\xi} & \\
&  &  & \ddots\\
&  &  &
\end{array}
\right)  \text{,}%
\end{equation}
where we abuse notation slightly and write $L_{\xi}$ also for the map from
$E^{\otimes m}\otimes_{\sigma_{1}}\mathcal{D}$ to $E^{\otimes(m+1)}%
\otimes_{\sigma_{1}}\mathcal{D}$ defined by the equation $L_{\xi}(\eta\otimes
h)=(\xi\otimes\eta)\otimes h$, $\eta\otimes h\in E^{\otimes m}\otimes
_{\sigma_{1}}\mathcal{D}$.

\begin{definition}
\label{MinIsoDilation}Let $E$ be a $W^{\ast}$-correspondence over the
$W^{\ast}$-algebra $M$ and let $(T,\sigma)$ be a completely contractive
covariant representation of $E$ on the Hilbert space $H$. Then the isometric
covariant representation $(V,\rho)$ just constructed is called the
\emph{minimal isometric dilation of }$(T,\sigma)$.
\end{definition}

The representation $(V,\rho)$ \emph{is} minimal in the sense that the smallest
subspace of $K$ that contains $H$ and reduces the set of operators
$\{V(\xi)\mid\xi\in E\}\cup\rho(M)$ is all of $K$. Thus the terminology is
justified. We note also that $(V,\rho)$ is unique up to unitary equivalence
\cite[Proposition 3.2]{MS98}.

If we let $\tilde{V}:E\otimes_{\rho}K\rightarrow K$ be the map that sends
$\xi\otimes k$ to $V(\xi)k$, then $\tilde{V}$ be written as the infinite
matrix
\begin{equation}
\tilde{V}=\left(
\begin{array}
[c]{llll}%
\tilde{T} & 0 & 0 & \cdots\\
\Delta & 0 & 0 & \\
0 & I & 0 & \\
0 & 0 & I & \\
&  &  & \ddots\\
&  &  &
\end{array}
\right)  \text{,} \label{Vtilde}%
\end{equation}
where the identity operators are interpreted as the maps that identify
$E\otimes_{\sigma_{n+1}}(E^{\otimes n}\otimes_{\sigma_{1}}\mathcal{D})$ with
$E^{\otimes(n+1)}\otimes_{\sigma_{1}}\mathcal{D}$. It is then an easy
calculation to see $\tilde{V}^{\ast}\tilde{V}=I$ on $K$, so that $\tilde{V}$
is an isometry (which confirms our assertion that $(V,\rho)$ is an isometric
dilation of $(T,\sigma)$), and that
\begin{equation}
\tilde{V}\tilde{V}^{\ast}=\left(
\begin{array}
[c]{llll}%
\tilde{T}\tilde{T}^{\ast} & \tilde{T}\Delta^{\ast} & 0 & \cdots\\
\Delta\tilde{T}^{\ast} & \Delta^{2} & 0 & \\
0 & 0 & I & \\
&  &  & \ddots\\
&  &  &
\end{array}
\right)  \text{,} \label{Vvstar}%
\end{equation}
a calculation that we shall use in a moment. Let $\tilde{T}_{n}:E^{\otimes
n}\otimes H\rightarrow H$ be the $n^{th}$-generalized power of $\tilde{T}$
(Remark \ref{GenPowers}) and similarly let $\tilde{V}_{n}$, mapping
$E^{\otimes n}\otimes K$ to $K$ be the $n^{th}$-generalized power of
$\tilde{V}$. Then, of course, each $\tilde{T}_{n}$ is a contraction, while
each $\tilde{V}_{n}$ is an isometry. Also, as we mentioned in Remark
\ref{GenPowers}, $\tilde{V}_{n+m}=\tilde{V}_{n}(I_{n}\otimes\tilde{V}%
_{m})=\tilde{V}_{m}(I_{m}\otimes\tilde{V}_{n})$, where $I_{n}$ is the identity
map on $E^{\otimes n}$. The importance of the $\tilde{V}_{n}$ for our purposes
is that they implement endomorphisms of the \emph{commutant} of $\rho(M)$.
Indeed, if we set
\[
L(x)=\tilde{V}(I_{E}\otimes x)\tilde{V}^{\ast}\text{,}%
\]
$x\in\rho(M)^{\prime}$, then $L$ is an endomorphism of $\rho(M)^{\prime} $
and
\[
L^{n}(x)=\tilde{V}_{n}(I_{E}\otimes x)\tilde{V}_{n}^{\ast}\text{,}%
\]
for all $n\geq0$ and $x\in\rho(M)^{\prime}$\cite[Lemma 2.3]{MS99}. It follows
easily that for a subspace $\mathcal{M}$ of $K$ that is invariant under
$\rho(M)$, the range of $L^{n}(P_{\mathcal{M}})$ is the span
\[
\overline{span}\{V(\xi_{1})\cdots V(\xi_{n})h:h\in\mathcal{M}\text{,}\;\xi
_{1},...,\xi_{n}\in E\}.
\]

\begin{definition}
\label{wandering} A subspace $\mathcal{M}$ of $K$ that is invariant for
$\rho(M)$ is called a \emph{wandering} \emph{subspace}, and the projection
$P_{\mathcal{M}}$ of $K$ onto $\mathcal{M}$ is called a \emph{wandering
projection},\emph{\ }if for every $n\neq m$, $L^{n}(P_{\mathcal{M}})$ and
$L^{m}(P_{\mathcal{M}})$ are orthogonal projections. For such a subspace we
shall write $L_{\infty}(\mathcal{M})\ $for the range of $\sum_{n\geq0}%
^{\oplus}L^{n}(P_{\mathcal{M}})$.
\end{definition}

Note that, whenever $\mathcal{M}\subseteq K$ is a wandering subspace, the map
$W_{\mathcal{M}}:\mathcal{F}(E)\otimes_{\rho}\mathcal{M}\rightarrow L_{\infty
}(\mathcal{M})$ defined by sending $\xi_{1}\otimes\cdots\otimes\xi_{n}\otimes
k\in E^{\otimes n}\otimes_{\rho}\mathcal{M}$ to $V(\xi_{1})\cdots V(\xi
_{n})k\in L_{\infty}(\mathcal{M})$ is a Hilbert space isometry. Note, too,
that for $a\in M$ and $\xi\in E$, we have
\begin{equation}
W_{\mathcal{M}}(\varphi_{\infty}(a)\otimes I_{\mathcal{M}})=\rho
(a)W_{\mathcal{M}} \label{W1}%
\end{equation}
and
\begin{equation}
W_{\mathcal{M}}(T_{\xi}\otimes I_{\mathcal{M}})=V(\xi)W_{\mathcal{M}}.
\label{W2}%
\end{equation}

We also write $P_{n}$ for $\tilde{V}_{n}\tilde{V}_{n}^{\ast}$, so that
$P_{n}=L^{n}(I)$. Of course $P_{1}$ is given by the matrix (\ref{Vvstar}).
Then $\{P_{n}\}_{n=1}^{\infty}$ is a decreasing sequence of projections in
$\rho(M)^{\prime}$. We set $Q_{n}=P_{n}-P_{n+1}$ and $Q_{0}=I-P_{1}$, so that
$\sum_{k=0}^{\infty}Q_{k}=I-P_{\infty}$, where $P_{\infty}=\wedge P_{n}$. By
\cite[Corrolary 2.4]{MS99}, $Q_{0}$ is a wandering projection, $Q_{k}%
=L^{k}(Q_{0})$ and $Q_{\infty}:=\sum_{k=0}^{\infty}L^{k}(Q_{0})=\sum
_{k=0}^{\infty}Q_{k}=I-P_{\infty}$.

\begin{lemma}
\label{Qn}With the notation just established, we have for every $\xi\in E$ and
$m\geq0$,
\[
V(\xi)Q_{m}=Q_{m+1}V(\xi)
\]
and
\[
V(\xi)Q_{\infty}=Q_{\infty}V(\xi).
\]

\end{lemma}

\begin{proof}
For $k\in K$ we have $V(\xi)Q_{m}k=\tilde{V}(\xi\otimes Q_{m}k)=\tilde
{V}(I\otimes Q_{m})(\xi\otimes k)=\tilde{V}(I\otimes Q_{m})\tilde{V}^{*}%
\tilde{V}(\xi\otimes k)=Q_{m+1}V(\xi)k.$
\end{proof}

If we let $\rho_{0}$ be the restriction of $\rho$ to the range of $Q_{0}$,
then it follows from \cite[Theorem 2.9]{MS99} that $(V,\rho)$ may be written
as the direct sum
\[
(V,\rho)=(V_{ind},\rho_{ind})\oplus(V_{\infty},\rho_{\infty})
\]
where $(V_{ind},\rho_{ind})$ is (unitarily equivalent to) the representation
of $E$ that is induced by $\rho_{0}$, while $(V_{\infty},\rho_{\infty})$ is
the restriction to $P_{\infty}(K)$ and is fully coisometric in the sense of
\cite{MS98, MS99, MS04}, meaning that $\tilde{V}_{\infty}$ is a coisometry.
Thus, $\tilde{V}_{\infty}$ is a unitary operator on $P_{\infty}(K)$.

\subsection{C.N.C. and C.$_{0}$ Representations}

Our goal is to describe how $H$ sits in the dilation space $K$. The analysis
we present follows Sz.-Nagy and Foia\c{s}, as one might imagine. \ However,
there are some important refinements that are due to Popescu \cite{gP89} and
we need to extend these to our situation. As a first step, we have the
following observation, which may be \textquotedblleft dug out
of\textquotedblright\ \cite{MS04} (see Lemma 7.8, in particular.) However,
since we need a bit more than is explicit there, we present a proof.

\begin{lemma}
\label{Q0}Write $K_{0}$ for the range, $Q_{0}(K)$, of the projection $Q_{0}$. Then

\begin{itemize}
\item[(i)] $K_{0}=\overline{Q_{0}(H)}=\overline{\{\Delta_{\ast}^{2}%
h\oplus(-\Delta\tilde{T}^{\ast}h):h\in H\}}\subseteq H\oplus\mathcal{D}$.

\item[(ii)] The map $u$ that sends $\Delta_{\ast}^{2}h\oplus(-\Delta\tilde
{T}^{\ast}h)$ to $\Delta_{\ast}h$ is an isometry from $K_{0}$ onto
$\mathcal{D}_{\ast}$.

\item[(iii)] The equation $\rho(a)u=\sigma(a)u=u\rho(a)$ holds for all $a\in
M$.
\end{itemize}
\end{lemma}

\begin{proof}
>From the minimality of $K$ it follows that $I_{K}=\vee_{n=0}^{\infty}%
L^{n}(P_{H})=P_{H}\vee P_{1}$. Since $Q_{0}$ and $P_{1}$ are orthogonal, by
definition, we have $Q_{0}(K)=\overline{Q_{0}(H)}$. The other equality follows
when we write $Q_{0}$ matricially as
\[
Q_{0}=I-\tilde{V}\tilde{V}^{\ast}=\left(
\begin{array}
[c]{cccc}%
I_{H}-\tilde{T}\tilde{T}^{\ast} & -\tilde{T}\Delta & 0 & \ldots\\
-\Delta\tilde{T}^{\ast} & I_{{}}-\Delta^{2} & 0 & \\
0 & 0 & 0 & \\
\vdots &  &  & \ddots\\
&  &  &
\end{array}
\right)  \text{,}%
\]
as we may, by equation (\ref{Vvstar}). This proves (i). For (ii) we compute:
\begin{align*}
\langle\Delta_{\ast}^{2}h\oplus(-\Delta\tilde{T}^{\ast}h),\Delta_{\ast}%
^{2}h\oplus(-\Delta\tilde{T}^{\ast}h)\rangle &  =\langle\Delta_{\ast}%
^{4}h,h\rangle+\langle\tilde{T}\Delta^{2}\tilde{T}^{\ast}h,h\rangle\\
&  =\langle\Delta_{\ast}^{2}(\Delta_{\ast}^{2}+\tilde{T}\tilde{T}^{\ast
})h,h\rangle=\langle\Delta_{\ast}^{2}h,h\rangle\text{,}%
\end{align*}
which proves the assertion. The proof of part (iii) is immediate from the
following computation, which is valid for all $a\in M$ and $h\in H$:
\begin{align*}
\rho(a)(\Delta_{\ast}^{2}h\oplus(-\Delta\tilde{T}^{\ast}h))  &  =\sigma
(a)\Delta_{\ast}^{2}h\oplus(\varphi(a)\otimes I_{H})(-\Delta\tilde{T}\ast h)\\
&  =\Delta_{\ast}^{2}\sigma(a)h\oplus(-\Delta(\varphi(a)\otimes I_{H}%
)\tilde{T}^{\ast}h)=\Delta_{\ast}^{2}\sigma(a)h\oplus(-\Delta\tilde{T}^{\ast
}\sigma(a)h).
\end{align*}

\end{proof}

The following terminology is adopted from \cite{gP89, gP89a}, which, in turn,
derives from the work of Sz.-Nagy and Foia\c{s} (see \cite{szNF70}).

\begin{definition}
\label{c0}\hfill

\begin{enumerate}
\item[(i)] A covariant representation $(T,\sigma)$ will be called a
$\emph{C.}_{\emph{0}}$\emph{-representation} if $P_{\infty}=0$ (equivalently,
if $K=L_{\infty}(K_{0})$).

\item[(ii)] A covariant representation $(T,\sigma)$ will be called
\emph{completely non coisometric} (abbreviated c.n.c. ) in case $K=L_{\infty
}(K_{0})\vee L_{\infty}(\mathcal{D})$.
\end{enumerate}
\end{definition}

\begin{remark}
\label{CNCDecomp}It is\ shown in Remark 7.2 of \cite{MS04} that given a
covariant representation $(T,\sigma)$ of $E$ on a Hilbert space $H$, then $H $
may be written as $H=H_{1}\oplus H_{2}$ so that if $T$ and $\sigma$ are
written as matrices relative to this decomposition, then
\[
\sigma=\left(
\begin{array}
[c]{cc}%
\sigma_{1} & 0\\
0 & \sigma_{2}%
\end{array}
\right)  \text{,}%
\]
i.e., $\sigma$ is reduced by $H_{1}$ and $H_{2}$, and
\[
T(\cdot)=\left(
\begin{array}
[c]{cc}%
T_{1}(\cdot) & 0\\
X(\cdot) & T_{2}(\cdot)
\end{array}
\right)  \text{,}%
\]
where $(T_{1},\sigma_{1})$ is a covariant representation that is c.n.c. and
where $(T_{2},\sigma_{2})$ is a covariant representation with the property
that all the generalized powers of $\tilde{T}_{2}$ are coisometries. Further,
$H_{2}$ may be described as $\{h\in H\mid\left\Vert \tilde{T}_{n}^{\ast
}h\right\Vert =\left\Vert h\right\Vert $ \emph{for all} $n\}$, i.e., $H_{2}$
is the largest space on which all the generalized powers $\tilde{T}_{n}^{\ast
}$ act isometrically. Thus $(T,\sigma)$ is c.n.c. if and only if there is no
non-zero vector $h$ such that $\left\Vert \tilde{T}_{n}^{\ast}h\right\Vert
=\left\Vert h\right\Vert $ for all $n$.
\end{remark}

For our purpose here, the significance of the concept \textquotedblleft
c.n.c.\textquotedblright\ is the condition in the second of the following two
lemmas. The first is Proposition 7.15 of \cite{MS04}, while the second is
Lemma 7.10 of \cite{MS04}.

\begin{lemma}
\label{C0}Let $(T,\sigma)$ be a covariant representation of a $W^{\ast}%
$-correspondence on a Hilbert space $H$ and let $(V,\rho)$ be its minimal
isometric dilation acting on $K$. Then the following conditions are equivalent.

\begin{enumerate}
\item[(i)] $(T,\sigma)$ is of class $C._{0}$, i.e. $P_{\infty}=0$.

\item[(ii)] $\wedge\tilde{V}_{k}\tilde{V}_{k}^{\ast}=0$, which happens if and
only if $\Vert\tilde{V}_{k}^{\ast}k\Vert\rightarrow0$ for all $k\in K $.

\item[(iii)] $\tilde{T}_{k}\tilde{T}_{k}^{\ast}\rightarrow0$ in the weak
operator topology on $B(H)$, which happens if and only if $\Vert\tilde{T}%
_{k}^{\ast}h\Vert\rightarrow0$ for all $h\in H$.

\item[(iv)] $(V,\rho)$ is an induced representation.
\end{enumerate}

\noindent So, in particular, if $\Vert\tilde{T}\Vert<1$ then $(T,\sigma)$ is a
$C._{0}$-representation.
\end{lemma}

\begin{lemma}
\label{cnc}\hfill

\begin{enumerate}
\item[(i)] Every $C._{0}$-representation is c.n.c.

\item[(ii)] A representation is c.n.c if and only if $P_{\infty}%
(K)=\overline{P_{\infty}(L_{\infty}(\mathcal{D}))}$, which happens if and only
if $P_{\infty}(H)\subseteq\overline{P_{\infty}(L_{\infty}(\mathcal{D}))}$.
\end{enumerate}
\end{lemma}

We record here for the sake of reference the following statement, which is
part of Theorem 7.3 of \cite{MS04}.

\begin{theorem}
\label{CNCExtension}If $(T,\sigma)$ is a completely contractive covariant
representation of a $W^{\ast}$-correspondence on a Hilbert space $H$, and if
$(T,\sigma)$ is completely non-coisometric, then $T\times\sigma$ extends to an
ultraweakly continuous, completely contractive representation of the Hardy
algebra, $H^{\infty}(E)$, on $H$.
\end{theorem}

\subsection{Characteristic Operators}

We now turn to the construction of the characteristic operator and the
characteristic function associated to a covariant representation. At the
outset, we do not require that the representation is \emph{c.n.c.} We fix a
completely contractive covariant representation $(T,\sigma)$ acting on the
Hilbert space $H$. We maintain the notation just developed. However, we shall
write $W_{\infty}$ for the Hilbert space isomorphism that we would have
written $W_{K_{0}}$ earlier in order to lighten the notation. So $W_{\infty}$
is a Hilbert space isomorphism from $\mathcal{F}(E)\otimes_{\rho}K_{0}$ onto
$L_{\infty}(K_{0})$ that satisfies (\ref{W1}) and (\ref{W2}) (with $K_{0}$ in
place of $\mathcal{M}$). We also write $u$ for the isometry from $K_{0}$ onto
$\mathcal{D}_{\ast}$ described in Lemma~\ref{Q0}. It induces an isometry,
written $I_{\mathcal{F}(E)}\otimes u$ from $\mathcal{F}(E)\otimes K_{0}$ onto
$\mathcal{F}(E)\otimes\mathcal{D}_{\ast}$.

\begin{definition}
\label{charop} Let $(T,\sigma)$ be a completely contractive covariant
representation of the $W^{\ast}$-correspondence $E$ over the $W^{\ast}%
$-algebra $M$ and let $(V,\rho)$ be the minimal isometric dilation of
$(T,\sigma)$. Also, in the notation just established, let $\tau_{1}$ be the
restriction of $\rho$ to $\mathcal{D}$ and let $\tau_{2}$ be the restriction
of $\rho$ (or $\sigma$) to $\mathcal{D}_{\ast}$. Then the operator $\Theta
_{T}$ defined from $\mathcal{F}(E)\otimes_{\rho}\mathcal{D}$ to $\mathcal{F}%
(E)\otimes_{\rho}\mathcal{D}_{\ast}$ by the equation%
\begin{equation}
\Theta_{T}:=(I_{\mathcal{F}(E)}\otimes u)\circ W_{\infty}^{\ast}(I-P_{\infty
})W_{\mathcal{D}} \label{Theta_T}%
\end{equation}
is called the \emph{characteristic operator }of the representation\emph{\ }%
$(T,\sigma)$.
\end{definition}

\begin{remarks}
\label{CharOpBackground}\hfill

\begin{enumerate}
\item[(i)] Evidently, $\Theta_{T}$ is a contraction. Indeed, since
$I_{\mathcal{F}(E)}\otimes u$, $W_{\infty}$ and $W_{\mathcal{D}}$ are all
isometries, the \textquotedblleft only\textquotedblright\ things that keep
$\Theta_{T}$ from being an isometry are the relations among the range of
$W_{\infty}$, the range of $I-P_{\infty}$ and $W_{\mathcal{D}}$. Further,
given the calculations involving $W_{\infty}$, $I-P_{\infty}$ and
$W_{\mathcal{D}}$ that we have made so far, it is clear that $\Theta_{T}$
carries \emph{some }information about the location of $H$ in the space of the
minimal isometric dilation of $(T,\sigma)$. Our goal is to show that under the
assumption that our representation is c.n.c., it carries all the information
and is a complete unitary invariant for the representation $(T,\sigma)$.

\item[(ii)] We frequently will want to refer to the entire system,
$(\Theta_{T},\mathcal{D},\mathcal{D}_{\ast},\tau_{1},\tau_{2})$, as the
characteristic operator for the covariant representation $(T,\sigma)$.

\item[(iii)] By definition, $\tau_{2}$ is the restriction of $\sigma$ to
$\mathcal{D}_{\ast}$. By definition of the minimal isometric dilation of
$(T,\sigma)$, $(V,\rho)$, $\tau_{1}$ really is the restriction of $\sigma
\circ\varphi$ to $\mathcal{D}$ regarded as the zero$^{th}$ component in the
natural decomposition of $\mathcal{F}(E)\otimes_{\sigma_{1}}\mathcal{D} $. See
Definition \ref{MinIsoDilation}.

\item[(iv)] Although $\Theta_{T}$ is defined to be a map between the two
Hilbert spaces, $\mathcal{F}(E)\otimes\mathcal{D}$ and $\mathcal{F}%
(E)\otimes\mathcal{D}_{\ast}$, which are different, in general, we shall
occasionally identify $\Theta_{T}$ with the $2\times2$ operator matrix%
\[
\left(
\begin{array}
[c]{cc}%
0 & 0\\
\Theta_{T} & 0
\end{array}
\right)
\]
in $B(\mathcal{F}(E)\otimes(\mathcal{D}\oplus\mathcal{D}_{\ast}))$.
\end{enumerate}
\end{remarks}

Several basic properties of $\Theta_{T}$ are established in the following lemma.

\begin{lemma}
\label{Y}The characteristic operator $\Theta_{T}$ is a contraction that
satisfies the equations
\begin{equation}
(\varphi_{\infty}(a)\otimes I_{\mathcal{D}_{\ast}})\Theta_{T}=\Theta
_{T}(\varphi_{\infty}(a)\otimes I_{\mathcal{D}})\text{,}\;\;\;a\in M
\label{Ymod}%
\end{equation}
and
\begin{equation}
\Theta_{T}(T_{\xi}\otimes I_{\mathcal{D}})=(T_{\xi}\otimes I_{\mathcal{D}%
_{\ast}})\Theta_{T}\text{,}\;\;\xi\in E\text{.} \label{commute}%
\end{equation}
That is, $\Theta_{T}$ intertwines the representations of $H^{\infty}(E)$
induced by $\tau_{1}$ and $\tau_{2}$. Further, if $(T,\sigma)$ is a $C._{0}%
$-representation, then $Q_{\infty}=I$, i.e., $P_{\infty}=0$, and $\Theta_{T}$
is an isometry from $\mathcal{F}(E)\otimes\mathcal{D}$ into $\mathcal{F}%
(E)\otimes\mathcal{D}_{\ast}$.
\end{lemma}

\begin{proof}
We already have remarked that $\Theta_{T}$ is a contraction. The other parts
of the lemma are immediate consequences of equation (\ref{W1}), Lemma~\ref{Q0}%
, Lemma~\ref{Qn} and the equations $W_{\mathcal{D}}(T_{\xi}\otimes
I_{\mathcal{D}})=V(\xi)W_{\mathcal{D}}$ and $W_{\infty}^{\ast}V(\xi)=(T_{\xi
}\otimes I_{K_{0}})W_{\infty}^{\ast}$, which are easy to check.
\end{proof}

As we shall show in Theorem \ref{InnerCriterion}, there is a conditioned
converse to the last assertion in Lemma \ref{Y}.

The representations $\tau_{1}$ and $\tau_{2}$, defined above, need not be
faithful. Indeed, they need not even be jointly faithful. This will have to be
accommodated in our analysis. Accordingly, we let $e$ be the central
projection in $M$ such that $Ker(\tau_{1}\oplus\tau_{2})=eM$. The following
lemma reveals its significance.

\begin{lemma}
\label{tau}The projection $e$ is the largest central projection $q$ in $M$
such that the operator $\sigma(q)\tilde{T}$ is a partial isometry with initial
space $\varphi(q)E\otimes H$ and final space $\sigma(q)H$.
\end{lemma}

\begin{proof}
The projection $e$ is the largest central projection $q$ with $\tau
_{1}(q)=\tau_{2}(q)=0$. But this holds if and only if both the restriction of
$\sigma(q)$ to $\Delta_{\ast}H$ and the restriction of $\varphi(q)\otimes
I_{H}$ to $\Delta(E\otimes H)$ are equal to zero. This is equivalent to the
requirements that $\sigma(q)(I_{H}-\tilde{T}\tilde{T}^{\ast})=0$ and
$(\varphi(q)\otimes I)(I_{E\otimes H}-\tilde{T}^{\ast}\tilde{T})=0$. Since
$\sigma(q)\tilde{T}=\tilde{T}(\varphi(q)\otimes I)$, the proof is complete.
\end{proof}

\begin{corollary}
If either $\Vert\tilde{T}\Vert<1$ or $M$ is a factor, then $\tau_{1}\oplus
\tau_{2}$ is faithful and $e=0$.
\end{corollary}

\subsection{Characteristic Functions}

The technology involving the theory of duality that was developed in
\cite{MS04}, and is summarized in Section \ref{preliminaries}, requires
faithful representations of the $W^{\ast}$-algebras in question. Since
$\tau_{1}\oplus\tau_{2}$ need not be faithful, we will \textquotedblleft
supplement\textquotedblright\ it to build a faithful representation of $M$.
For this purpose, we introduce the following terminology.

\begin{definition}
\label{Supplemental}For $i=1,2$, let $\tau_{i}:M\rightarrow B(\mathcal{E}%
_{i})$ be a normal representation of $M$ on $\mathcal{E}_{i}$ and let $e$ be
the central projection such that $\ker(\tau_{1}\oplus\tau_{2})=eM$. Chose a
faithful representation $\pi_{0}$ of $M$ on a Hilbert space $H_{0}$ and let
$\tau_{0}$ be the representation of $M$ on $\pi_{0}(e)H_{0}$ obtained by
restricting $\pi_{0}$ to $eM$. Form the Hilbert space $\mathcal{E}:=\pi
_{0}(H_{0})\oplus\mathcal{E}_{1}\oplus\mathcal{E}_{2}$ and let $\tau:=\tau
_{0}\oplus\tau_{1}\oplus\tau_{2}$ be the (necessarily faithful) representation
of $M$ on $\mathcal{E}$. Then we call $\mathcal{E}$ a \emph{supplemental
space} for the pair of representations $\tau_{1}$ and $\tau_{2}$, we shall
call the representation $\tau$ of $M$ on $\mathcal{E}$ a \emph{supplemental
representation} and we shall simply call the pair $(\mathcal{E},\tau)$ a
\emph{supplement }for\emph{\ }$\tau_{1}$ and $\tau_{2}$.
\end{definition}

Evidently, if $\tau_{1}$ and $\tau_{2}$ are jointly faithful, then
$(\mathcal{E}_{1}\oplus\mathcal{E}_{2},\tau_{1}\oplus\tau_{2})$ is the only
possible supplement for\emph{\ }$\tau_{1}$ and $\tau_{2}$. We shall see
shortly that the use of supplemental spaces and representations is a matter of
convenience only and that the constructs we consider do not depend in any
material way on the choice of $\pi_{0}$ used to define them.

Suppose, now, that $(\Theta_{T},\mathcal{D},\mathcal{D}_{\ast},\tau_{1}%
,\tau_{2})$ is the characteristic operator determined by a covariant
representation $(T,\sigma)$ of $E$. We fix once and for all a supplement
$(\mathcal{G},\tau)$ for $\tau_{1}$ and $\tau_{2}$ and we consider
$\mathcal{F}(E)\otimes_{\tau}\mathcal{G}$ as written as the direct sum%
\begin{equation}
\mathcal{F}(E)\otimes_{\tau}\mathcal{G}=(\mathcal{F}(E)\otimes_{\pi_{0}}%
H_{0})\oplus(\mathcal{F}(E)\otimes_{\tau_{1}}\mathcal{D})\oplus(\mathcal{F}%
(E)\otimes_{\tau_{2}}\mathcal{D}_{\ast})\text{.} \label{DecompFEG}%
\end{equation}
Corresponding to this direct sum decomposition of $\mathcal{F}(E)\otimes
_{\sigma}\mathcal{G}$, we shall identify $\Theta_{T}$ with the block matrix%
\begin{equation}
\left(
\begin{array}
[c]{ccc}%
0 & 0 & 0\\
0 & 0 & 0\\
0 & \Theta_{T} & 0
\end{array}
\right)  \text{.} \label{MatrixFEG}%
\end{equation}
Since $\Theta_{T}$ satisfies equations (\ref{commute}) and (\ref{Ymod}), it
follows that this block matrix actually lies in the commutant of
$\tau^{\mathcal{F}(E)}(H^{\infty}(E))$. Hence we may take its Fourier
transform relative to $\tau$ as in Remark \ref{FunctionalRep}, obtaining an
element $\hat{\Theta}_{T}\in H^{\infty}(E^{\tau})$ such that
\begin{equation}
\hat{\Theta}_{T}\otimes I_{\mathcal{G}}=U\Theta_{T}U^{\ast}\text{,}
\label{ThetaHatT}%
\end{equation}
where $U$ is the Fourier transform from $\mathcal{F}(E)\otimes_{\tau
}\mathcal{G}$ onto $\mathcal{F}(E^{\tau})\otimes_{\iota}\mathcal{G}$ defined
in Proposition \ref{U}. Since elements of $H^{\infty}(E^{\tau})$ may be viewed
as functions on the unit ball of $E$ (see Remark \ref{FunctionalRep}), we will
think of $\hat{\Theta}_{T}$ as being so represented when we wish. The
following lemma records some of the properties of this transform and shows
that it does not really depend on the choice of $\tau$ and $\mathcal{G}$.

\begin{lemma}
\label{char}Let $\hat{\Theta}_{T}$ be the element of $H^{\infty}(E^{\tau})$
defined in equation (\ref{ThetaHatT}) using the Fourier transform $U$ from
$\mathcal{F}(E)\otimes_{\tau}\mathcal{G}$ onto $\mathcal{F}(E^{\tau}%
)\otimes_{\iota}\mathcal{G}$. Also let $q_{1}$ be the projection from
$\mathcal{G}$ onto $\mathcal{D}$ and $q_{2}$ be the projection onto
$\mathcal{D}_{\ast}$. Then both $q_{1}$ and $q_{2}$ lie in $\tau(M)^{\prime}$, and

\begin{enumerate}
\item[(i)] $U^{\ast}(q_{i}\otimes I_{\mathcal{G}})U=I_{\mathcal{F}(E)}\otimes
q_{i}$, $i=1,2$.

\item[(ii)] $\hat{\Theta}_{T}=q_{2}\hat{\Theta}_{T}q_{1}$ and, if $(T,\sigma)$
is a $C_{\cdot0}$-representation, then $\hat{\Theta}_{T}^{\ast}\hat{\Theta
}_{T}=q_{1}$.

\item[(iii)] For every $\xi\in E$ with $\Vert\xi\Vert<1$, $q_{2}\hat{\Theta
}_{T}(\xi)q_{1}=\hat{\Theta}_{T}(\xi)$.
\end{enumerate}
\end{lemma}

\begin{proof}
To prove (i), recall first that, for $\eta_{1},\ldots,\eta_{k}$ in $E^{\tau}$
and $h\in\mathcal{G}$,
\[
U^{\ast}(\eta_{1}\otimes\ldots\otimes\eta_{k}\otimes h)=(I_{E^{\otimes(k-1)}%
}\otimes\eta_{1})\cdots(I_{E}\otimes\eta_{k-1})\eta_{k}(h).
\]
For $q\in\tau(M)^{\prime}$ and $\eta\in E^{\tau}$, we have $q\cdot\eta
=(I_{E}\otimes q)\eta$. (This is the left action of $\tau(M)^{\prime}$ on
$E^{\tau}$.) Thus, for such $q$,
\[
U^{\ast}(q\otimes I_{\mathcal{G}})(\eta_{1}\otimes\ldots\otimes\eta_{k}\otimes
h)=U^{\ast}(q\eta_{1}\otimes\ldots\otimes\eta_{k}\otimes h)=
\]%
\[
=(I_{E^{\otimes k}}\otimes q)(I_{E^{\otimes(k-1)}}\otimes\eta_{1})\cdots
(I_{E}\otimes\eta_{k-1})\eta_{k}(h)=(I_{E^{\otimes k}}\otimes q)U^{\ast}%
(\eta_{1}\otimes\ldots\otimes\eta_{k}\otimes h).
\]
This proves (i). From the construction of the operator $\Theta_{T}$ above it
follows that $\Theta_{T}=(I_{\mathcal{F}(E)}\otimes q_{2})\Theta
_{T}(I_{\mathcal{F}(E)}\otimes q_{1})$. Thus, using (i), $U\Theta_{T}U^{\ast
}=U(I_{\mathcal{F}(E)}\otimes q_{2})U^{\ast}U\Theta_{T}U^{\ast}%
U(I_{\mathcal{F}(E)}\otimes q_{1})U^{\ast}=(q_{2}\otimes I_{\mathcal{G}%
})U\Theta_{T}U^{\ast}(q_{1}\otimes I_{\mathcal{G}}).$ Since $\hat{\Theta}%
_{T}\otimes I_{\mathcal{G}}=U\Theta_{T}U^{\ast}$, we proved (ii). For $X\in
H^{\infty}(E^{\tau})$, $X(\xi)$ is the image, under a certain representation
of $H^{\infty}(E^{\tau})$ defined by $\xi$, of $X$. Thus the map $X\mapsto
X(\xi)$ is multiplicative and it carries elements of $\tau(M)^{\prime}$ to
themselves. Part (iii) thus follows from part (ii).
\end{proof}

The lemma shows that $q_{2}\hat{\Theta}_{T}(\xi)q_{1}=\hat{\Theta}_{T}(\xi) $
for all $\xi$ in the open unit ball of $E$ and so we may view $\hat{\Theta
}_{T}$ as a function from the open unit ball of $E$ to $B(\mathcal{D}%
,\mathcal{D}_{\ast})$.

The properties of $\hat{\Theta}_{T}$ will be formalized in the following definition.

\begin{definition}
\label{chardata}Given a $W^{\ast}$-algebra $M$ and a $W^{\ast}$-correspondence
$E$ over $M$, a \emph{characteristic function} is a system $(\Theta
,\mathcal{E}_{1},\mathcal{E}_{2},\tau_{1},\tau_{2})$ with the following properties:

\begin{enumerate}
\item[(i)] For $i=1,2$, $\mathcal{E}_{i}$ is a Hilbert space and $\tau_{i}$ is
a representation of $M$ on $\mathcal{E}_{i}$.

\item[(ii)] If $(\mathcal{E},\tau)$ is a supplement for $\tau_{1}$ and
$\tau_{2}$, and if $q_{i}$ is the projection of $\mathcal{E}$ onto
$\mathcal{E}_{i}$, $i=1,2$, then $\Theta$ is a contraction in $H^{\infty
}(E^{\tau})$ satisfying $\Theta=q_{2}\Theta q_{1}$.
\end{enumerate}

If, in addition, $\Theta$ satisfies $\Theta^{\ast}\Theta=q_{1}$ then $\Theta$
will be called an \emph{inner characteristic function}.
\end{definition}

Very often we shall write $\Theta$ for the tuple $(\Theta,\mathcal{E}%
_{1},\mathcal{E}_{2},\tau_{1},\tau_{2})$. Also, given a characteristic
function, we shall freely use the notation set in Definition~\ref{chardata}
(i.e. $\mathcal{E}_{i}$, $\tau_{i}$ and $q_{i}$).

\begin{definition}
\label{characteristic}If $(T,\sigma)$ is a covariant representation of the
$W^{\ast}$-correspondence $E$, then the system $(\hat{\Theta}_{T}%
,\mathcal{D},\mathcal{D}_{\ast},\tau_{1},\tau_{2})$ defined by equation
(\ref{ThetaHatT}), or simply $\hat{\Theta}_{T}$, will be called \emph{the
characteristic function of the representation $(T,\sigma)$}.
\end{definition}

The following result is familiar from the theory of single operators.\ It is
the \textquotedblleft converse\textquotedblright\ of Lemma \ref{Y}.

\begin{theorem}
\label{InnerCriterion}Let $E$ be a $W^{\ast}$-correspondence over a $W^{\ast}%
$-algebra $M$ and let $(T,\sigma)$ be a c.n.c. representation of $E $ on the
Hilbert space $H$. Then the characteristic function $\hat{\Theta}_{T}$ of the
covariant representation $(T,\sigma)$ is inner if and only if $(T,\sigma)$ is
a $C_{\cdot0}$-representation.
\end{theorem}

\begin{proof}
Lemma \ref{Y} shows that if $(T,\sigma)$ is a $C_{\cdot0}$ representation,
then $\Theta_{T}$ is an isometry. Consequently, $\hat{\Theta}_{T}$ is inner.
To prove the converse, observe that from the definition of $\Theta_{T}$,
equation (\ref{Theta_T}), $\Theta_{T}$ is an isometry if and only if
$L_{\infty}(\mathcal{D})\subseteq L_{\infty}(K_{0})$. However, by our
assumption that $(T,\sigma)$ is c.n.c., we know by definition (Definition
\ref{c0}) that $L_{\infty}(\mathcal{D})\vee L_{\infty}(K_{0})=K$. Hence, if
$\Theta_{T}$ is an isometry, so that $L_{\infty}(\mathcal{D})\subseteq
L_{\infty}(K_{0})$, we conclude that $L_{\infty}(K_{0})=K$. Hence by
definition (Definition \ref{c0}), $(T,\sigma)$ is a $C_{\cdot0}$
representation. Since $\Theta_{T}$ is an isometry if and only if $\hat{\Theta
}_{T}$ is inner, the proof is complete.
\end{proof}

\subsection{Pointwise Evaluations}

Of course several natural questions arise at this point: Is every
characteristic function the characteristic function of some representation? If
so, how is the representation constructed? What is the level of uniqueness
among the constructs? Before tackling these, we first compute the values
$\hat{\Theta}_{T}(\xi)$ for the characteristic function of a covariant
representation $(T,\sigma)$. The calculations will play roles in the sequel.
The initial step of our analysis is the following computation.

\begin{lemma}
\label{Sxi} Let $P_{\mathcal{D}}$ (resp. $P_{\mathcal{D}_{\ast}}$) denote the
projection of $\mathcal{F}(E)\otimes_{\tau}\mathcal{D}=\mathcal{D}%
\oplus(E\otimes_{\tau}\mathcal{D})\oplus\cdots$ onto the zeroth summand,
$\mathcal{D}$ (resp. the projection of $\mathcal{F}(E)\otimes_{\tau
}\mathcal{D}_{\ast}$ onto the zeroth summand $\mathcal{D}_{\ast}$). Also, for
$\xi\in E$, $\left\Vert \xi\right\Vert \leq1$, write $L_{\xi^{\otimes k}}$ for
the operator from $\mathcal{F}(E)\otimes\mathcal{D}_{\ast}$ to $\mathcal{F}%
(E)\otimes\mathcal{D}_{\ast}$ defined by formula $L_{\xi^{\otimes k}}%
\eta\otimes h=\xi^{\otimes k}\otimes\eta\otimes h$, when $k\geq1$, and let
$L_{\xi^{\otimes0}}$ be the identity operator. Then for every $\xi$ in the
open unit ball of $E$ and every $g\in\mathcal{D}$
\[
\hat{\Theta}_{T}(\xi)g=\sum_{k=0}^{\infty}P_{\mathcal{D}}L_{\xi^{\otimes k}%
}^{\ast}(I_{\mathcal{F}(E)}\otimes u)W_{\infty}^{\ast}Q_{\infty}g\text{,}%
\]
where, recall, $W_{\infty}:\mathcal{F}(E)\otimes K_{0}\rightarrow K$ and
$u:K_{0}\rightarrow\mathcal{D}_{\ast}$ are the isometries defined above.
\end{lemma}

\begin{proof}
Note first that, since $\Vert\xi\Vert<1$, the sum converges. To establish the
formula we shall fix such a $\xi$ and show that for every $R\in H^{\infty
}(E^{\tau})$ and every $g\in\mathcal{D}$,
\begin{equation}
R(\xi)g=\sum_{k=0}^{\infty}P_{\mathcal{G}}L_{\xi^{\otimes k}}^{\ast}U^{\ast
}(R\otimes I_{\mathcal{G}})Ug \label{evaluation}%
\end{equation}
where, recall, $\mathcal{G}$ is $\pi_{0}(e)H_{0}\oplus\mathcal{D}%
\oplus\mathcal{D}_{\ast}$ and $U$ is the Fourier transform from $\mathcal{F}%
(E)\otimes_{\tau}\mathcal{G}$ to $\mathcal{F}(E^{\tau})\otimes_{\iota
}\mathcal{G}$, while $P_{\mathcal{G}}$ is the projection of $\mathcal{F}%
(E)\otimes_{\tau}\mathcal{G}$ onto the zeroth summand. When $R=\hat{\Theta
}_{T}$ we will obtain the desired result since $U^{\ast}(\hat{\Theta}%
_{T}\otimes I_{\mathcal{G}})U=\Theta_{T}$. Suppose first that $R=\varphi
_{\infty}(b)\in H^{\infty}(E^{\tau})$ (with $b\in M$). Then $R(\xi)=b$ by
definition. Computing the right hand side of (\ref{evaluation}) we get first
$U^{\ast}(\varphi_{\infty}(b)\otimes I_{\mathcal{G}})Ug=U^{\ast}bg=bg$ and,
thus, the only non zero term in the sum is the one corresponding to $k=0$. In
this event the sum is then equal to $bg$, proving the equation for constant
functions. Now fix $m\geq1$ , let $\eta=\eta_{1}\otimes\eta_{2}\cdots
\otimes\eta_{m}\in(E^{\tau})^{\otimes m}$ and let $R=T_{\eta}\in H^{\infty
}(E^{\tau})$. Then, from the definition of $R(\xi)$,
\[
R(\xi)=(T_{\eta_{1}})(\xi)\cdots(T_{\eta_{m}})(\xi)=(L_{\xi}^{\ast}\eta
_{1})\cdots(L_{\xi}^{\ast}\eta_{m})
\]
where $\eta_{i}$ is viewed as a map from $\mathcal{G}$ into $E\otimes_{\tau
}\mathcal{G}$ and, thus, $L_{\xi}^{\ast}\eta_{i}\in B(\mathcal{G})$.

To compute the right hand side of (\ref{evaluation}) in this case we first
compute $U^{\ast}(T_{\eta}\otimes I_{\mathcal{G}})Ug=U^{\ast}(\eta\otimes
g)=(I_{(E^{\tau})^{\otimes(m-1)}}\otimes\eta_{1})\cdots(I_{E^{\tau}}%
\otimes\eta_{m-1})\eta_{m}(g)$. It then follows that the only non zero term in
the sum is the one that corresponds to $k=m$. A simple computation shows that
\[
L_{\xi^{\otimes m}}^{\ast}(I_{(E^{\tau})^{\otimes(m-1)}}\otimes\eta_{1}%
)\cdots(I_{E^{\tau}}\otimes\eta_{m-1})\eta_{m}(g)=(L_{\xi}^{\ast}\eta
_{1})\cdots(L_{\xi}^{\ast}\eta_{m})g.
\]
This, by linearity, proves (\ref{evaluation}) for a $\sigma$-weakly dense
subset of $H^{\infty}(E^{\tau})$. Since both sides of the equation are
$\sigma$-weakly continuous (as a function of $R$), equation (\ref{evaluation}) follows.
\end{proof}

To use lemma \ref{Sxi} to calculate the values of $\hat{\Theta}_{T}(\xi)$, we
compute the series appearing in the lemma term by term. For $k=0$ we have
$P_{\mathcal{D}_{\ast}}\Theta_{T}g=P_{\mathcal{D}_{\ast}}uW_{\infty}^{\ast
}Q_{\infty}g=uQ_{0}g$ for\ all $g\in\mathcal{D}$. Suppose $g=\Delta
(\theta\otimes h)$, $\theta\otimes h\in E\otimes_{\tau}H$. Then
\[
uQ_{0}g=uQ_{0}\Delta(\theta\otimes h)=u(-\tilde{T}\Delta^{2}(\theta\otimes
h)+(I_{\mathcal{D}}-\Delta^{2})\Delta(\theta\otimes h))=
\]%
\[
=u(-\Delta_{\ast}^{2}\tilde{T}(\theta\otimes h)+\Delta\tilde{T}^{\ast}%
\tilde{T}(\theta\otimes h))=-\Delta_{\ast}\tilde{T}(\theta\otimes
h)=-\tilde{T}\Delta(\theta\otimes h)=-\tilde{T}g.
\]
Since vectors $g$ of the form $\Delta(\theta\otimes h)$ generate $\mathcal{D}%
$, we see that
\begin{equation}
uQ_{0}|\mathcal{D}=-\tilde{T}|\mathcal{D} \label{Y0}%
\end{equation}
and, thus, the zeroth term in the expression of $\Theta(\xi)$ is $-\tilde{T}$.
To compute the other terms recall first, from equation (\ref{Vtilde}), that we
can write $\tilde{V}^{\ast}$ matricially as
\begin{equation}
\tilde{V}^{\ast}=\left(
\begin{array}
[c]{cccc}%
\tilde{T}^{\ast} & \Delta & 0 & \ldots\\
0 & 0 & I & \\
0 & 0 & 0 & \\
&  &  & \ddots
\end{array}
\right)  :H\oplus\mathcal{D}\oplus\cdots\rightarrow E\otimes H\oplus
E\otimes\mathcal{D}\oplus\cdots
\end{equation}
Thus, for $g\in\mathcal{D}$, $\tilde{V}^{\ast}g=\Delta g$ and $\tilde{V}%
_{2}^{\ast}g=(I_{E}\otimes\tilde{V}^{\ast})\tilde{V}^{\ast}g=(I_{E}%
\otimes\tilde{T}^{\ast})\Delta g$. In fact, for every $k\geq2$,
\[
\tilde{V}_{k}^{\ast}g=(I_{E^{\otimes(k-1)}}\otimes\tilde{T}^{\ast}%
)\cdots(I_{E}\otimes\tilde{T}^{\ast})\Delta g
\]
for $g\in\mathcal{D}$.

The next term ($k=1$) applied to $g=\Delta(\theta\otimes h)$ is
\[
L_{\xi}^{\ast}(I_{E}\otimes u)W_{\infty}^{\ast}Q_{\infty}\Delta(\theta\otimes
h)=L_{\xi}^{\ast}(I_{E}\otimes u)W_{\infty}^{\ast}\tilde{V}(I_{E}\otimes
Q_{0})\tilde{V}^{\ast}\Delta(\theta\otimes h)=
\]%
\[
=L_{\xi}^{\ast}(I_{E}\otimes uQ_{0})\tilde{V}^{\ast}\Delta(\theta\otimes h).
\]
Using the comments above, $\tilde{V}^{\ast}\Delta(\theta\otimes h)=\Delta
^{2}(\theta\otimes h)$. Also, for $h\in H$, we have
\[
uQ_{0}h=u(\Delta_{\ast}^{2}h\oplus(-\Delta\tilde{T}^{\ast}h))=\Delta_{\ast
}h\text{,}%
\]
by lemma~\ref{Q0}. Hence $L_{\xi}^{\ast}(I_{E}\otimes uQ_{0})V_{\infty}^{\ast
}\Delta(\theta\otimes h)=L_{\xi}^{\ast}(I_{E}\otimes\Delta_{\ast})\Delta
^{2}(\theta\otimes h)=\Delta_{\ast}L_{\xi}^{\ast}\Delta^{2}(\theta\otimes h)$.
It follows that the term that corresponds to $k=1$ in the expression of
$\hat{\Theta}_{T}(\xi)$ is $\Delta_{\ast}L_{\xi}^{\ast}\Delta$. Continuing in
this fashion, we see that for $k\geq2$ and $g=\Delta(\theta\otimes h)$, we
have
\[
L_{\xi^{\otimes k}}^{\ast}(I_{E^{\otimes k}}\otimes u)W_{\infty}^{\ast}%
\tilde{V}_{k}(I_{E^{\otimes k}}\otimes Q_{0})\tilde{V}_{k}^{\ast}%
g=L_{\xi^{\otimes k}}^{\ast}(I_{E^{\otimes k}}\otimes uQ_{0})\tilde{V}%
_{k}^{\ast}\Delta(\theta\otimes h)=
\]%
\[
=L_{\xi^{\otimes k}}^{\ast}(I_{E^{\otimes k}}\otimes\Delta_{\ast
})(I_{E^{\otimes(k-1)}}\otimes\tilde{T}^{\ast})\cdots(I_{E}\otimes\tilde
{T}^{\ast})\Delta g=\Delta_{\ast}(L_{\xi}^{\ast}\tilde{T}^{\ast})^{k-1}L_{\xi
}^{\ast}\Delta g.
\]
Thus the $k$th term in the expression of $\hat{\Theta}_{T}(\xi)$ is
$\Delta_{\ast}(L_{\xi}^{\ast}\tilde{T}^{\ast})^{k-1}L_{\xi}^{\ast}\Delta$. We
now summarize the discussion above.

\begin{theorem}
\label{expression}The values of the characteristic function $\hat{\Theta}_{T}
$ on the open unit ball of $E$ can be written as
\[
\hat{\Theta}_{T}(\xi)=-\tilde{T}|\mathcal{D}+\sum_{k=1}^{\infty}\Delta_{\ast
}(L_{\xi}^{\ast}\tilde{T}^{\ast})^{k-1}L_{\xi}^{\ast}\Delta|\mathcal{D}%
=-\tilde{T}|\mathcal{D}+\Delta_{\ast}(I-L_{\xi}^{\ast}\tilde{T})^{-1}L_{\xi
}^{\ast}\Delta|\mathcal{D}\text{.}%
\]

\end{theorem}

\begin{remark}
Theorem \ref{expression} may be viewed as a realization formula associated
with the unitary operator matrix
\[
\left(
\begin{array}
[c]{cc}%
-\tilde{T}|\mathcal{D} & \Delta_{\ast}\\
\Delta|\mathcal{D} & \tilde{T}^{\ast}%
\end{array}
\right)  :\mathcal{D}\oplus H\rightarrow\mathcal{D}_{\ast}\oplus
(E\otimes_{\sigma}H).
\]
(See e.g. \cite{AM02} .) Evidently, it is an exact analogue of the formula for
the characteristic operator function for a single contraction operator
\cite{szNF70}.
\end{remark}

\subsection{Models from Characteristic Functions}

Suppose we are given a characteristic function $(\Theta,\mathcal{E}%
_{1},\mathcal{E}_{2},\tau_{1},\tau_{2})$ and form $\check{\Theta}:=U^{\ast
}(\Theta\otimes I_{\mathcal{E}})U$ where, recall, $\mathcal{E}:=\mathcal{\pi
}_{0}(e)H\oplus\mathcal{E}_{1}\oplus\mathcal{E}_{2}$ is the Hilbert space
described in Definition \ref{chardata} and $U:\mathcal{F}(E)\otimes_{\tau
}\mathcal{E}\rightarrow\mathcal{F}(E^{\tau})\otimes_{\iota}\mathcal{E}$ is the
Fourier transform described in Proposition \ref{U} and Remark
\ref{FunctionalRep}. Then $\check{\Theta}$ commutes with the operators
$T_{\xi}\otimes I_{\mathcal{E}}$ and $\varphi_{\infty}(a)\otimes
I_{\mathcal{E}}$ for $\xi\in E$ and $a\in M$. Since $\Theta=q_{2}\Theta q_{1}%
$, we can use the argument of the proof of Lemma~\ref{char} (i) to show that
$U^{\ast}(q_{i}\otimes I_{\mathcal{E}})U=I_{\mathcal{F}(E)}\otimes q_{i}$,
$i=1,2$, and, thus, $\check{\Theta}(\mathcal{F}(E)\otimes\mathcal{E}%
_{1})=U^{\ast}(\Theta\mathcal{F}(E^{\tau})\otimes\mathcal{E})=U^{\ast}%
(q_{2}\otimes I)(\Theta\mathcal{F}(E^{\tau})\otimes\mathcal{E})\subseteq
\mathcal{F}(E)\otimes\mathcal{E}_{2}$. It follows that, for $\xi\in E$ and
$a\in M$,
\begin{equation}
\check{\Theta}(T_{\xi}\otimes I_{\mathcal{E}_{1}})=(T_{\xi}\otimes
I_{\mathcal{E}_{2}})\check{\Theta} \label{Yxi}%
\end{equation}
and
\begin{equation}
\check{\Theta}(\varphi_{\infty}(a)\otimes I_{\mathcal{E}_{1}})=(\varphi
_{\infty}(a)\otimes I_{\mathcal{E}_{2}})\check{\Theta} \label{Ya}%
\end{equation}

Our objective is to show that there is a covariant representation $(T,\sigma)$
of $E$ such that $\Theta_{T}=\check{\Theta}$. To this end, we write
$\Delta_{\check{\Theta}}=(I_{\mathcal{F}(E)\otimes\mathcal{E}_{1}}%
-\check{\Theta}^{\ast}\check{\Theta})^{1/2}\in B(\mathcal{F}(E)\otimes
_{\tau_{1}}\mathcal{E}_{1})$ and set
\begin{equation}
K(\Theta):=(\mathcal{F}(E)\otimes_{\tau_{2}}\mathcal{E}_{2})\oplus
\overline{\Delta_{\check{\Theta}}(\mathcal{F}(E)\otimes_{\tau_{1}}%
\mathcal{E}_{1})}\subseteq\mathcal{F}(E)\otimes_{\tau}\mathcal{E} \label{kth}%
\end{equation}
and
\begin{equation}
H(\Theta):=((\mathcal{F}(E)\otimes_{\tau_{2}}\mathcal{E}_{2})\oplus
\overline{\Delta_{\check{\Theta}}(\mathcal{F}(E)\otimes_{\tau_{1}}%
\mathcal{E}_{1})})\ominus\{\check{\Theta}\xi\oplus\Delta_{\check{\Theta}}%
\xi\mid\xi\in\mathcal{F}(E)\otimes_{\tau_{1}}\mathcal{E}_{1}\}. \label{hth}%
\end{equation}
Note that if $\Theta$ is inner, then $\check{\Theta}^{\ast}\check{\Theta
}=U^{\ast}(q_{1}\otimes I_{\mathcal{E}})U=I_{\mathcal{F}(E)}\otimes q_{1}$ and
so $\Delta_{\check{\Theta}}=0$. Thus, in this case $K(\Theta)=\mathcal{F}%
(E)\otimes_{\tau_{2}}\mathcal{E}_{2}$ and $H(\Theta)=(\mathcal{F}%
(E)\otimes_{\tau_{2}}\mathcal{E}_{2})\ominus\check{\Theta}(\mathcal{F}%
(E)\otimes_{\tau_{1}}\mathcal{E}_{1})$.

We shall also write $P_{\Theta}$ for the projection from $K(\Theta)$ onto
$H(\Theta)$.

\begin{lemma}
\label{charep}Let $\Theta$ be a characteristic function and let $\check
{\Theta}$, $K(\Theta)$ and $H(\Theta)$ be the operator and spaces just
defined. For every $a\in M$ and $\xi\in E$ we define the operators $S_{\Theta
}(\xi)$ and $\psi_{\Theta}(a)$ on $\Delta_{\check{\Theta}}(\mathcal{F}%
(E)\otimes\mathcal{E}_{1})$ by the formulae%
\begin{equation}
S_{\Theta}(\xi)\Delta_{\check{\Theta}}g=\Delta_{\check{\Theta}}(T_{\xi}\otimes
I_{\mathcal{E}_{1}})g,\;\;g\in\mathcal{F}(E)\otimes\mathcal{E}_{1} \label{sth}%
\end{equation}
and
\begin{equation}
\psi_{\Theta}(a)\Delta_{\check{\Theta}}g=\Delta_{\check{\Theta}}%
(\varphi_{\infty}(a)\otimes I_{\mathcal{E}_{1}})g,\;\;g\in\mathcal{F}%
(E)\otimes\mathcal{E}_{1}. \label{psith}%
\end{equation}
Also, we define the following operators on $K(\Theta)$:%
\begin{equation}
V_{\Theta}(\xi)=(T_{\xi}\otimes I_{\mathcal{E}_{2}})\oplus S_{\Theta}(\xi)
\label{vth}%
\end{equation}
and
\begin{equation}
\rho_{\Theta}(a)=(\varphi_{\infty}(a)\otimes I_{\mathcal{E}_{2}})\oplus
\psi_{\Theta}(a). \label{rhoth}%
\end{equation}
Then

\begin{enumerate}
\item[(i)] $(S_{\Theta},\psi_{\Theta})$ and $(V_{\Theta},\rho_{\Theta}) $ are
isometric covariant representations of $E$ on $\overline{\Delta_{\check
{\Theta}}(\mathcal{F}(E)\otimes\mathcal{E}_{1})}$ and $K(\Theta)$ respectively.

\item[(ii)] The space $K(\Theta)\ominus H(\Theta)$ is invariant for
$(V_{\Theta},\rho_{\Theta})$ and, thus, the compression of $(V_{\Theta}%
,\rho_{\Theta})$ to $H(\Theta)$, which we denote by $(T_{\Theta}%
,\sigma_{\Theta})$, is a completely contractive covariant representation of
$E$. Explicitly,
\begin{equation}
T_{\Theta}(\xi)=P_{\Theta}V_{\Theta}(\xi)|H(\Theta),\;\;\xi\in E \label{tth}%
\end{equation}
and
\begin{equation}
\sigma_{\Theta}(a)=P_{\Theta}\rho_{\Theta}(a)|H(\Theta),\;\;a\in M.
\label{sigth}%
\end{equation}

\end{enumerate}
\end{lemma}

\begin{proof}
In (i) it is enough to prove the statement about $(S_{\Theta},\psi_{\Theta})$.
We shall write $\Delta$ for $\Delta_{\check{\Theta}}$. Then, for $\xi\in E$ ,
$a,b\in M$ and $g\in\mathcal{F}(E)\otimes\mathcal{E}_{1}$, $S_{\Theta}(a\xi
b)\Delta g=\Delta(T_{a\xi b}\otimes I_{\mathcal{E}_{1}})g=\Delta
(\varphi_{\infty}(a)T_{\xi}\varphi_{\infty}(b)\otimes I_{\mathcal{E}_{1}%
})g=\psi_{\Theta}(a)S_{\Theta}(\xi)\psi_{\Theta}(b)\Delta g.$ This proves the
covariance property. Since $(\varphi_{\infty}(a)\otimes I_{\mathcal{E}_{2}%
})\check{\Theta}=\check{\Theta}(\varphi_{\infty}(a)\otimes I_{\mathcal{E}_{1}%
})$, $\varphi_{\infty}(a)\otimes I_{\mathcal{E}_{1}}$ commutes with $\Delta$
and $\psi_{\Theta}$ is a $\ast$-representation of $M$. To show that the
representation is isometric we compute for $\eta_{i}\otimes h_{i}%
\in\mathcal{F}(E)\otimes\mathcal{E}_{1}$, $i=1$ and $2$,
\begin{multline*}
\langle S_{\Theta}(\xi_{1})\Delta(\eta_{1}\otimes h_{1}),S_{\Theta}(\xi
_{2})\Delta(\eta_{2}\otimes h_{2})\rangle\\
=\langle\Delta(\xi_{1}\otimes\eta_{1}\otimes h_{1}),\Delta(\xi_{2}\otimes
\eta_{2}\otimes h_{2}\rangle\\
=\langle\xi_{1}\otimes\eta_{1}\otimes h_{1},\xi_{2}\otimes\eta_{2}\otimes
h_{2}\rangle-\langle\check{\Theta}(\xi_{1}\otimes\eta_{1}\otimes h_{1}%
),\check{\Theta}(\xi_{2}\otimes\eta_{2}\otimes h_{2})\rangle\\
=\langle\eta_{1}\otimes h_{1},\varphi_{\infty}(\langle\xi_{1},\xi_{2}%
\rangle)\eta_{2}\otimes h_{2}\rangle-\langle\xi_{1}\otimes\check{\Theta}%
(\eta_{1}\otimes h_{1}),\xi_{2}\otimes\check{\Theta}(\eta_{2}\otimes
h_{2})\rangle\\
=\langle\eta_{1}\otimes h_{1},\varphi_{\infty}(\langle\xi_{1},\xi_{2}%
\rangle)\eta_{2}\otimes h_{2}\rangle-\langle\check{\Theta}(\eta_{1}\otimes
h_{1}),(\varphi_{\infty}(\langle\xi_{1},\xi_{2}\rangle)\otimes I_{\mathcal{E}%
_{2}})\check{\Theta}(\eta_{2}\otimes h_{2})\rangle\\
=\langle\eta_{1}\otimes h_{1},\varphi_{\infty}(\langle\xi_{1},\xi_{2}%
\rangle)\eta_{2}\otimes h_{2}\rangle-\langle\check{\Theta}(\eta_{1}\otimes
h_{1}),\check{\Theta}(\varphi_{\infty}(\langle\xi_{1},\xi_{2}\rangle)\otimes
I_{\mathcal{E}_{1}})(\eta_{2}\otimes h_{2})\rangle\\
=\langle\Delta^{2}(\eta_{1}\otimes h_{1}),(\varphi_{\infty}(\langle\xi_{1}%
,\xi_{2}\rangle)\otimes I_{\mathcal{E}_{1}})(\eta_{2}\otimes h_{2})\rangle\\
=\langle\Delta(\eta_{1}\otimes h_{1}),\psi_{\Theta}(\langle\xi_{1},\xi
_{2}\rangle)\Delta(\eta_{2}\otimes h_{2})\rangle.
\end{multline*}
This shows that the representation is isometric. To prove (ii) all we have to
show is the invariance of $K(\Theta)\ominus H(\Theta)=\{\check{\Theta}%
g\oplus\Delta g:g\in\mathcal{F}(E)\otimes\mathcal{E}_{1}\}$ under the
representation $(V_{\Theta},\rho_{\Theta})$. However, this is an immediate
application of equations (\ref{Yxi}) and (\ref{Ya}).
\end{proof}

\begin{definition}
\label{CanonicalModel}Let $(\Theta,\mathcal{E}_{1},\mathcal{E}_{2},\tau
_{1},\tau_{2})$ be a characteristic function. Then the covariant
representation $(T_{\Theta},\sigma_{\Theta})$ on $H(\Theta)$ defined from
$\Theta$ in Lemma \ref{charep} is called the \emph{canonical model constructed
from }$\Theta$. If $(\Theta,\mathcal{E}_{1},\mathcal{E}_{2},\tau_{1},\tau
_{2})$ is the characteristic function of a covariant representation
$(T,\sigma)$, i.e., if $\Theta=\hat{\Theta}_{T}$, then $(T_{\Theta}%
,\sigma_{\Theta})$ will be called the \emph{canonical model for }$(T,\sigma)$.
\end{definition}

We begin to justify this terminology in the following Theorem.

\begin{theorem}
\label{TSisT}Let $(T,\sigma)$ be a c.n.c. covariant representation of $E$,
with characteristic operator $\Theta_{T}$. Let $\Theta:=\hat{\Theta}_{T}$ be
the associated characteristic function and $(T_{\Theta},\sigma_{\Theta})$ be
the canonical model for $(T,\sigma)$. Then $(T,\sigma)$ and $(T_{\Theta
},\sigma_{\Theta})$ are unitarily equivalent.
\end{theorem}

\begin{proof}
Let $H$ be the representation space of $(T,\sigma)$ and recall the definition
of $\Theta_{T}$ in Definition \ref{Theta_T}. Note that in the notation of
Lemma \ref{charep}, $\Theta_{T}=\check{\Theta}$. Write
\[
\Phi_{1}=W_{\infty}(I_{\mathcal{F}(E)}\otimes u^{\ast}):\mathcal{F}%
(E)\otimes\mathcal{D}_{\ast}\rightarrow K
\]
where $K$ and $\mathcal{D}_{\ast}$ are the spaces associated with $(T,\sigma)$
and its minimal isometric dilation, and where $W_{\infty}$ and $u$ are the
operators defined in the discussion preceding Definition \ref{Theta_T}. Then
$\Phi_{1}$ is an isometry whose range is $L_{\infty}(K_{0}) $. We also define
$\Phi_{2}:\Delta_{\check{\Theta}}(\mathcal{F}(E)\otimes\mathcal{D})\rightarrow
P_{\infty}(K)$ by the equation
\[
\Phi_{2}(\Delta_{\check{\Theta}}x)=P_{\infty}(W_{\mathcal{D}}x)\;,\;\;x\in
\mathcal{F}(E)\otimes\mathcal{D}.
\]
Since the representation is c.n.c., $P_{\infty}(L_{\infty}(\mathcal{D}%
))=P_{\infty}(K)$ by part (ii) of Lemma \ref{cnc} and so $\Phi_{2}$ is
surjective. We show that it is an isometry. For this we compute
\[
\Vert\Delta_{\check{\Theta}}\xi\Vert^{2}=\langle(I-\check{\Theta}^{\ast}%
\check{\Theta})\xi,\xi\rangle=\Vert\xi\Vert^{2}-\Vert(I_{\mathcal{F}%
(E)}\otimes u_{\ast})\check{\Theta}\xi\Vert^{2}.
\]
By definition of $\check{\Theta}=\Theta_{T}$ (equation (\ref{Theta_T})), the
last expression is equal to
\[
\Vert\xi\Vert^{2}-\Vert W_{\infty}^{\ast}Q_{\infty}W_{\mathcal{D}}\xi\Vert
^{2}=\Vert W_{\mathcal{D}}\xi\Vert^{2}-\Vert Q_{\infty}W_{\mathcal{D}}\xi
\Vert^{2}=\Vert P_{\infty}W_{\mathcal{D}}\xi\Vert^{2}.
\]
Thus $\Phi_{2}$ is a unitary operator onto $P_{\infty}(K)$. Setting $\Phi
=\Phi_{1}\oplus\Phi_{2}$ we obtain a unitary operator from $K(\Theta)$ onto
$K$.

Next we show that $\Phi$ maps $H(\Theta)$ onto $H$. Fix $x\in\mathcal{F}%
(E)\otimes\mathcal{D}$. Then by definition,
\[
\Phi(\check{\Theta}x\oplus\Delta_{\check{\Theta}}x)=W_{\infty}(I_{\mathcal{F}%
(E)}\otimes u^{\ast})\check{\Theta}x+P_{\infty}(W_{\mathcal{D}}x).
\]
So, if $x\in\mathcal{D}$, with $\mathcal{D}$ regarded as the zero$^{th}$
summand of $\mathcal{F}(E)\otimes\mathcal{D}$, we find from the definition of
$\check{\Theta}=\Theta_{T}$ (equation (\ref{Theta_T})) that $\Phi
(\check{\Theta}x\oplus\Delta_{\check{\Theta}}x)=W_{\infty}(I_{\mathcal{F}%
(E)}\otimes u^{\ast})\check{\Theta}x+P_{\infty}(W_{\mathcal{D}}x)=Q_{\infty
}x+P_{\infty}x=x$. Since $\mathcal{D}$ is orthogonal to $H$, we see that
$\Phi(\check{\Theta}x\oplus\Delta_{\check{\Theta}}x)\in H^{\perp}$. If
$n\geq1$, then for $x=\xi\otimes d\in E^{\otimes n}\otimes\mathcal{D}$, we
also have
\begin{align*}
W_{\infty}(I_{\mathcal{F}(E)}\otimes u^{\ast})\check{\Theta}x+P_{\infty
}(W_{\mathcal{D}}x)  &  =W_{\infty}(I_{\mathcal{F}(E)}\otimes u^{\ast}%
)(\xi\otimes\check{\Theta}d)+P_{\infty}(V_{n}(\xi)d)\\
&  =V_{n}(\xi)Q_{\infty}d+V_{n}(\xi)P_{\infty}d=V_{n}(\xi)d\in H^{\perp}.
\end{align*}
Thus, we find that $\Phi(K(\Theta)\ominus H(\Theta))=\sum^{\oplus}%
V_{n}(E^{\otimes n})\mathcal{D}=K\ominus H$, and it follows that $\Phi$ maps
$H(\Theta)$ onto $H$.

Notice also that for $\xi\in E$
\begin{equation}
\Phi_{1}(T_{\xi}\otimes I)=W_{\infty}(T_{\xi}\otimes u^{\ast})=V(\xi
)W_{\infty}(I\otimes u^{\ast})=V(\xi)\Phi_{1}\text{,} \label{W}%
\end{equation}
while
\begin{multline*}
\Phi_{2}(S(\xi)\Delta_{\check{\Theta}}x)=\Phi_{2}(\Delta_{\check{\Theta}%
}(T_{\xi}\otimes I)x)=P_{\infty}(W_{\mathcal{D}}(T_{\xi}\otimes I)x)\\
=P_{\infty}(V(\xi)W_{\mathcal{D}}x)=V(\xi)P_{\infty}(W_{\mathcal{D}}%
x)=V(\xi)\Phi_{2}(\Delta_{\check{\Theta}}x).
\end{multline*}
Thus $\Phi$ intertwines $V$ and $V_{\Theta}$. To show that $\Phi$ also
intertwines $\rho$ and $\rho_{\Theta}$, we let $a\in M$ and compute:%
\begin{multline*}
\Phi_{1}(\varphi_{\infty}(a)\otimes I_{\mathcal{D}_{\ast}})=W_{\infty
}(I_{\mathcal{F}(E)}\otimes u^{\ast})(\varphi_{\infty}(a)\otimes I)=W_{\infty
}(\varphi_{\infty}(a)\otimes I)(I_{\mathcal{F}(E)}\otimes u^{\ast})\\
=\rho(a)W_{\infty}(I_{\mathcal{F}(E)}\otimes u^{\ast})\text{,}%
\end{multline*}
and, for $x\in\mathcal{F}(E)\otimes\mathcal{D}$,
\begin{multline*}
\Phi_{2}(\psi_{\Theta}(a)(\Delta_{\check{\Theta}}x))=\Phi_{2}(\Delta
(\varphi_{\infty}(a)\otimes I_{\mathcal{D}})x)=P_{\infty}(W_{\mathcal{D}%
}(\varphi_{\infty}(a)\otimes I)x)\\
=P_{\infty}(\rho(a)W_{\mathcal{D}}x)=\rho(a)P_{\infty}W_{\mathcal{D}}%
x=\rho(a)\Phi_{2}(\Delta_{\check{\Theta}}x).
\end{multline*}
It follows that the restriction of $\Phi$ to $H(\Theta)$ gives the desired equivalence.
\end{proof}

\begin{definition}
\label{Canonical Isomorph}Let $(T,\sigma)$ be a c.n.c. representation of the
$W^{\ast}$-correspondence on the Hilbert space $H$. Let $\Theta:=\hat{\Theta
}_{T}$ be the characteristic function for $(T,\sigma)$ and let $(T_{\Theta
},\sigma_{\Theta})$ on $H(\Theta)$ be the canonical model built from $\Theta$.
Then the Hilbert space isomorphism $\Phi$ from the Hilbert space $K$ of the
minimal isometric dilation of $(T,\sigma)$ to $K(\Theta)$ constructed in the
proof of Theorem \ref{TSisT} will be called the \emph{canonical (Hilbert
space) isomorphism (implementing a unitary equivalence between }$(T,\sigma)$
\emph{and }$(T_{\Theta},\sigma_{\Theta}))$ or simply the \emph{canonical
equivalence} for short.
\end{definition}

\begin{remark}
\label{Forward}Given a general characteristic function $(\Theta,\mathcal{E}%
_{1},\mathcal{E}_{2},\tau_{1},\tau_{2})$, the isometric representation
$(V_{\Theta},\rho_{\Theta})$ on $K(\Theta)$ defined by equations (\ref{vth})
and (\ref{rhoth}) is an isometric dilation of $(T_{\Theta},\sigma_{\Theta})$
by definition. In general, it need not be minimal. However, it will be under
hypotheses that we discuss shortly. See Lemma \ref{isomdil}.
\end{remark}

\subsection{Isomorphic Characteristic Functions}

\begin{definition}
\label{isomchar}Let $(\Theta,\mathcal{E}_{1},\mathcal{E}_{2},\tau_{1},\tau
_{2})$ and $(\Theta^{\prime},\mathcal{E}_{1}^{\prime},\mathcal{E}_{2}^{\prime
},\tau_{1}^{\prime},\tau_{2}^{\prime})$ be two characteristic functions. We
say that they are \emph{isomorphic} if there are Hilbert space isomorphisms
$W_{i}:\mathcal{E}_{i}\rightarrow\mathcal{E}_{i}^{\prime}$ that intertwine
$\tau_{i}$ and $\tau_{i}^{\prime}$, $i=1$ and $2$, and satisfy the equation
\begin{equation}
\check{\Theta}^{\prime}=(I_{\mathcal{F}(E)}\otimes W_{2})\check{\Theta
}(I_{\mathcal{F}(E)}\otimes W_{1}^{\ast}). \label{Ytheta}%
\end{equation}

\end{definition}

It follows easily from the way in which a characteristic function is
associated to a representation that if two c.n.c. representations are
(unitarily) equivalent then the associated characteristic functions are
isomorphic in the sense of Definition \ref{isomchar}. Conversely, a moment's
reflection on Lemma \ref{charep} and Proposition \ref{TSisT} reveals
immediately that given two isomorphic characteristic functions, the associated
representations are unitarily equivalent. We may therefore summarize our
analysis to this point in the following theorem that asserts that the
isomorphism class of a characteristic function of a c.n.c. representation is a
complete unitary invariant for the representation.

\begin{theorem}
\label{equiv}Two c.n.c. representations are unitarily equivalent if and only
if the associated characteristic functions are isomorphic.
\end{theorem}

\begin{remark}
\label{Thetaeq}The notion of isomorphism between two characteristic functions
$\Theta$ and $\Theta^{\prime}$ was defined using the operators $\check{\Theta
}$ and $\check{\Theta}^{\prime}$. One can also write an isomorphism directly
in terms of $\Theta$ and $\Theta^{\prime}$. For this, note first that if
Hilbert space isomorphisms $W_{i}:\mathcal{E}_{i}\rightarrow\mathcal{E}%
_{i}^{\prime}$ intertwining $\tau_{i}$ and $\tau_{i}^{\prime}$, $i=1,2$,
exist, then $\tau_{1}\oplus\tau_{2}$ and $\tau_{1}^{\prime}\oplus\tau
_{2}^{\prime}$ have the same kernels. So, if we choose a common representation
$\pi_{0}$ to define the supplements $\mathcal{E}$ and $\mathcal{E}^{\prime}$
for these representations, then the $W_{i}$'s may be extended to a Hilbert
space isomorphism $W:\mathcal{E}\rightarrow\mathcal{E}^{\prime}$ that
intertwines $\tau$ and $\tau^{\prime}$. On the other hand, if such a $W$
exists, then it restricts to give $W_{i}$'s that intertwine $\tau_{i}$ and
$\tau_{i}^{\prime}$. Also, equation (\ref{Ytheta}) is equivalent to the
equation
\[
\Theta^{\prime}\otimes I_{\mathcal{E}^{\prime}}=C(\Theta\otimes I_{\mathcal{E}%
})C^{\ast}%
\]
where $C$ is the unitary operator $C=U^{\prime}(I_{\mathcal{F}(E)}\otimes
W)U^{\ast}:F(E^{\tau})\otimes_{\iota}\mathcal{E}\rightarrow F(E^{\tau^{\prime
}})\otimes_{\iota^{\prime}}\mathcal{E}^{\prime}$ and $U$ and $U^{\prime}$ are
the evident Fourier transforms. In fact, one can show that for $\eta
\in(E^{\otimes k})^{\tau}=(E^{\tau})^{\otimes k}\subseteq F(E^{\tau})$ and
$h\in\mathcal{E}$,
\[
C(\eta\otimes h)=(I\otimes W)\eta W^{\ast}\otimes Wh\text{,}%
\]
where $(I\otimes W)\eta W^{\ast}$ is a map from $\mathcal{E}^{\prime}$ to
$E^{\otimes k}\otimes\mathcal{E}^{\prime}$ that lies in the $\tau$-dual of
$E^{\otimes k}$, which may be identified with $(E^{\tau})^{\otimes k}$ by
Proposition \ref{dual}. Consequently, the map $X\mapsto X^{\prime}$, where
$X^{\prime}\otimes I_{\mathcal{E}^{\prime}}=C(X\otimes I_{\mathcal{E}}%
)C^{\ast}$, is an isomorphism of $H^{\infty}(E^{\tau})$ onto $H^{\infty
}(E^{\tau^{\prime}})$. Once we use this map to identify the two algebras, we
see that the two characteristic functions are isomorphic in the sense of
Definition \ref{isomchar} if they are identified via this map. Since we do not
use this remark in the rest of the paper, we shall omit further details.
\end{remark}

\subsection{Models and Characteristic Functions: Completing the Circle}

\begin{lemma}
\label{cond}Let $(T,\sigma)$ be a c.n.c. representation of the $W^{\ast}%
$-correspondence on a Hilbert space, let $\mathcal{D}$ and $\mathcal{D}_{\ast
}$ be the defect spaces, let $\Theta=\Theta_{T}$ be its characteristic
operator and let $\Delta:=\Delta_{\Theta_{T}}=(I-\Theta^{\ast}\Theta)^{1/2}$ Then:

\begin{enumerate}
\item[(i)] There is no non zero vector $x\in\mathcal{D}$ such that
$x=P_{\mathcal{D}}\Theta^{\ast}P_{\mathcal{D}_{\ast}}\Theta x$.

\item[(ii)] $\overline{\Delta(\mathcal{F}(E)\otimes_{\sigma_{1}}\mathcal{D}%
)}=\overline{\Delta((\mathcal{F}(E)\otimes_{\sigma_{1}}\mathcal{D}%
)\ominus\mathcal{D})}$, where $\sigma_{1}=\sigma\circ\varphi$.
\end{enumerate}
\end{lemma}

\begin{proof}
It follows from the proof of Theorem~\ref{expression} (see equation
(\ref{Y0})) that $P_{\mathcal{D}_{\ast}}\Theta|\mathcal{D}=-\tilde{T}$. So (i)
amounts to the fact that the kernel of the positive operator $D=(I-T^{\ast
}T)^{1/2}$ restricted to the range of $D$ (i.e. to $\mathcal{D}$) is trivial.
Since this is obvious, (i) is proved. To prove (ii) note first that
$P_{\infty}(K)=\overline{span}\{V(\xi)P_{\infty}(k):\;\xi\in E,k\in
K\}=\overline{span}\{V(\xi)P_{\infty}(k):\xi\in E,k\in L_{\infty}%
(\mathcal{D})\}=\overline{span}\{P_{\infty}(V(\xi)k):\xi\in E,k\in L_{\infty
}(\mathcal{D})\}=\overline{P_{\infty}(L_{\infty}(\mathcal{D})\ominus
\mathcal{D})}$. So if $x\in\mathcal{F}(E)\otimes_{\sigma_{1}}\mathcal{D}$ and
if $\Phi_{2}$ is the isometry defined in Proposition~\ref{TSisT}, then
$\Phi_{2}(\Delta x)$ lies in $P_{\infty}(K)$. Hence $\Phi_{2}(\Delta x)=\lim
P_{\infty}y_{n}$ for some $y_{n}\in L_{\infty}(\mathcal{D})\ominus\mathcal{D}$
and so $\Delta x=\lim\Phi_{2}^{\ast}P_{\infty}y_{n}=\lim\Delta_{Y}%
(W_{\mathcal{D}}^{\ast}y_{n})$. It follows that $\Delta x\in\overline
{\Delta((\mathcal{F}(E)\otimes_{\sigma_{1}}\mathcal{D})\ominus\mathcal{D})}$.
\end{proof}

\begin{definition}
\label{PurePredictable}Let $\Theta=(\Theta,\mathcal{E}_{1},\mathcal{E}%
_{2},\tau_{1},\tau_{2})$ be a characteristic function and let $\Delta
:=(I-\check{\Theta}^{\ast}\check{\Theta})^{1/2}$.

\begin{enumerate}
\item[(i)] We say that $\Theta$ is \emph{pure} if there is no non-zero vector
$x$ in $\mathcal{E}_{1}$ so that $x=P_{\mathcal{E}_{1}}\check{\Theta}^{\ast
}P_{\mathcal{E}_{2}}\check{\Theta}x$.

\item[(ii)] We say that $\Theta$ is \emph{predictable }in case
\[
\overline{\Delta(\mathcal{F}(E)\otimes_{\tau_{1}}\mathcal{E}_{1})}%
=\overline{\Delta((\mathcal{F}(E)\otimes_{\tau_{1}}\mathcal{E}_{1}%
)\ominus\mathcal{E}_{1})}\text{.}%
\]

\end{enumerate}
\end{definition}

\begin{remark}
\label{ExplPredictable}The reason for the term \textquotedblleft
predictable\textquotedblright\ derives from the role of Hardy spaces in the
setting of prediction theory. Recall that if $M=\mathbb{C}=E$, then the Fock
space $\mathcal{F}(E)$ may be identified with the Hardy space $H^{2}%
(\mathbb{T})$. So, if $\mathcal{E}_{1}=\mathcal{E}_{2}=\mathbb{C}$ also, then
$\mathcal{F}(E)\otimes_{\tau_{1}}\mathcal{E}_{1}=H^{2}(\mathbb{T})$ as well,
and a characteristic function is simply a function $\theta\in H^{\infty
}(\mathbb{T})$ such that $\left\Vert \theta\right\Vert \leq1$, i.e., $\theta$
is a Schur function. (The function $\theta$ is pure if and only if $\theta$ is
not constant, by the maximum modulus principle.) The function $\delta
:=(1-|\theta|^{2})^{1/2}$ lies in $L^{\infty}(\mathbb{T)}$. To say that
$\theta$ is predictable is the same thing as saying that $\overline{\delta
H^{2}(\mathbb{T})}=\overline{\delta H_{0}^{2}(\mathbb{T})}$, where $H_{0}%
^{2}(\mathbb{T})$ is the space of those functions in $H^{2}(\mathbb{T})$ that
vanish at the origin. The connection with prediction theory is this: Suppose
$\{\xi_{n}\}_{n\in\mathbb{Z}}$ is a stationary Gaussian process with
covariance matrix $\{\Hat{\delta}(n-m)\}_{n,m\in\mathbb{Z}}$. Then the future,
$\bigvee_{n>0}\xi_{n}$, is contained in the past, $\bigvee_{n\leq0}\xi_{n}$,
i.e., the process $\{\xi_{n}\}_{n\in\mathbb{Z}}$ is predictable, if and only
if $\overline{\delta H^{2}(\mathbb{T})}=\overline{\delta H_{0}^{2}%
(\mathbb{T})}$. We note in passing that $\theta$ is predictable if and only if
$\overline{\delta H^{2}(\mathbb{T})}=L^{2}(\mathbb{T)}$ and that this is also
equivalent to the assertion that $\ln(\delta)\notin L^{1}(\mathbb{T)}$ by
Szeg\"{o}'s theorem.
\end{remark}

\begin{remark}
\label{rem}Let $\Theta$ be a characteristic function. Note that, for all
$\xi,\zeta$ in $E^{\otimes n}$, $\check{\Theta}$ commutes with both $T_{\xi
}\otimes I_{\mathcal{E}_{1}}$ and $T_{\zeta}^{\ast}T_{\xi}\otimes
I_{\mathcal{E}_{1}}$, since $T_{\zeta}^{\ast}T_{\xi}\in\varphi_{\infty}(M) $.
Thus $(T_{\zeta}^{\ast}\otimes I)\check{\Theta}(T_{\xi}\otimes I)=(T_{\zeta
}T_{\xi}^{\ast}\otimes I)\check{\Theta}=\check{\Theta}(T_{\zeta}^{\ast}\otimes
I)(T_{\xi}\otimes I)$. It follows that $(T_{\zeta}^{\ast}\otimes
I)\check{\Theta}$ and $\check{\Theta}(T_{\zeta}^{\ast}\otimes I)$ are equal
when restricted to $E^{\otimes m}\otimes\mathcal{E}_{1}$ for $m\geq n$.
\end{remark}

\begin{lemma}
\label{isomdil}Let $\Theta=(\Theta,\mathcal{E}_{1},\mathcal{E}_{2},\tau
_{1},\tau_{2})$ be a characteristic function that is pure and predictable.
Form its canonical model $(T,\sigma):=(T_{\Theta},\sigma_{\Theta})$ on the
Hilbert space $H(\Theta)$ and the isometric representation $(V,\rho
):=(V_{\Theta},\rho_{\Theta})$ on the Hilbert space $K(\Theta)$ as described
in Lemma \ref{charep}. Then $(V,\rho)$ is \emph{minimal} as an isometric
dilation of $(T,\sigma)$.
\end{lemma}

\begin{proof}
We already know that $(V,\rho)$ is an isometric dilation of $(T,\sigma)$ by
definition. So we need only prove minimality. For this, write $\mathcal{K} $
for the subspace
\[
\mathcal{K}=\overline{span}\{V(\xi)H(\Theta):\xi\in E\}.
\]
We shall show that $\mathcal{K}=K(\Theta)$. Fix a vector $x\in K(\Theta
)\ominus\mathcal{K}$. Since $x$ is orthogonal to $H(\Theta)$, we can write
$x=\check{\Theta}w_{0}+\Delta w_{0}$ for some $w_{0}\in\mathcal{F}%
(E)\otimes_{\tau_{1}}\mathcal{E}_{1}$, where as usual $\Delta:=(I-\check
{\Theta}^{\ast}\check{\Theta})^{1/2}$. For every $n\geq1$ and every $\xi\in
E^{\otimes n}$, $V(\xi)^{\ast}x\in H(\Theta)^{\perp}$ and we can find
$w(\xi)\in\mathcal{F}(E)\otimes\mathcal{E}_{1}$ such that
\[
V(\xi)^{\ast}(\check{\Theta}w_{0}+\Delta w_{0})=\check{\Theta}w(\xi)+\Delta
w(\xi).
\]
We now write $S$ for the operator $S_{\Theta}$ in Lemma \ref{charep} and
conclude from the previous equation that $(T_{\xi}^{\ast}\otimes
I)\check{\Theta}w_{0}=\check{\Theta}w(\xi)$ and $S(\xi)^{\ast}\Delta
w_{0}=\Delta w(\xi)$. Hence, for every $\xi,\zeta$ in $E^{\otimes n}$ we have
\[
\check{\Theta}^{\ast}(T_{\zeta}T_{\xi}^{\ast}\otimes I)\check{\Theta}%
w_{0}=\check{\Theta}^{\ast}\check{\Theta}(T_{\zeta}\otimes I)w(\xi)
\]
and
\[
\Delta S(\zeta)S(\xi)^{\ast}\Delta w_{0}=\Delta^{2}(T_{\zeta}\otimes
I)w(\xi),
\]
where we used the facts that $\check{\Theta}$ commutes with $T_{\zeta}\otimes
I$ and that, by definition, $S(\zeta)\Delta=\Delta(T_{\zeta}\otimes I)$).
Adding these two equations gives
\begin{equation}
\check{\Theta}^{\ast}(T_{\zeta}T_{\xi}^{\ast}\otimes I)\check{\Theta}%
w_{0}+\Delta S(\zeta)S(\xi)^{\ast}\Delta w_{0}=(T_{\zeta}\otimes
I)w(\xi)\text{.}%
\end{equation}
We shall write $e_{i}$ (respectively, $f_{i}$) for the projection of
$\mathcal{F}(E)\otimes_{\tau_{1}}\mathcal{E}_{1}$ (respectively,
$\mathcal{F}(E)\otimes_{\tau_{2}}\mathcal{E}_{2}$) onto $E^{\otimes i}%
\otimes_{\tau_{1}}\mathcal{E}_{1}$ (respectively, $E^{\otimes i}\otimes
_{\tau_{2}}\mathcal{E}_{2}$). Note that, for $\zeta\in E^{\otimes n}$ as
above, we have $e_{i}(T_{\zeta}\otimes I)w(\xi)=0$ if $i<n$. Thus, for $i<n$,
\[
e_{i}(\check{\Theta}^{\ast}(T_{\zeta}T_{\xi}^{\ast}\otimes I)\check{\Theta
}w_{0}+\Delta S(\zeta)S(\xi)^{\ast}\Delta w_{0})=0.
\]
It will be convenient to write $(R,\phi)$ for the (isometric) representation
of $E$ on $\mathcal{F}(E)\otimes_{\tau_{2}}\mathcal{E}_{2}$ defined by
$R(\xi)=T_{\xi}\otimes I_{\mathcal{E}_{2}}$ (for $\xi\in E$) and
$\phi(a)=\varphi_{\infty}(a)\otimes I$ for $a\in M$. Then the maps $\tilde
{R}_{n}:E^{\otimes n}\otimes_{\tau_{2}\circ\varphi_{\infty}}\mathcal{F}%
(E)\otimes_{\tau_{2}}\mathcal{E}_{2}\rightarrow\mathcal{F}(E)\otimes_{\tau
_{2}}\mathcal{E}_{2}$ are defined in the usual way. For $\zeta,\xi\in
E^{\otimes n}$ we write $\zeta\otimes\xi^{\ast}$ for the operator
$\zeta\otimes\xi^{\ast}$ on $E^{\otimes n}$ defined by the formula
$(\zeta\otimes\xi^{\ast})\xi^{\prime}=\zeta\langle\xi,\xi^{\prime}\rangle$.
The $C^{\ast}$-algebra generated by these operators is written $K(E^{\otimes
n})$ and it is $\sigma$-weakly dense in the $W^{\ast} $-algebra $\mathcal{L}%
(E^{\otimes n})$. We have $S(\zeta)S(\xi)^{\ast}=\tilde{S}_{n}((\zeta
\otimes\xi^{\ast})\otimes I)\tilde{S}_{n}^{\ast}$ and $T_{\zeta}T_{\xi}^{\ast
}\otimes I=\tilde{R}_{n}((\zeta\otimes\xi^{\ast})\otimes I)\tilde{R}_{n}%
^{\ast}$. Hence, for every $K\in K(E^{\otimes n})$ and every $i<n$,
\[
e_{i}(\check{\Theta}^{\ast}\tilde{R}_{n}(K\otimes I_{\mathcal{F}%
(E)\otimes\mathcal{E}_{2}})\tilde{R}_{n}^{\ast}\check{\Theta}w_{0}%
+\Delta\tilde{S}_{n}(K\otimes I_{\Delta(\mathcal{F}(E)\otimes\mathcal{E}_{1}%
)})\tilde{S}_{n}^{\ast}\Delta w_{0})=0.
\]
Noting that $I_{E^{\otimes n}}$ is in the $\sigma$-weak closure of
$K(E^{\otimes n})$ we conclude that
\[
e_{i}(\check{\Theta}^{\ast}\tilde{R}_{n}\tilde{R}_{n}^{\ast}\check{\Theta
}w_{0}+\Delta\tilde{S}_{n}\tilde{S}_{n}^{\ast}\Delta w_{0})=0
\]
for $i<n$. But $\tilde{R}_{n}\tilde{R}_{n}^{\ast}=\sum_{j=n}^{\infty}f_{j}$,
on the one hand, and $\tilde{S}_{n}\tilde{S}_{n}^{\ast}=I$ by our assumption
that $\Theta$ is predictable. Thus $e_{i}(\check{\Theta}^{\ast}(\sum
_{j=n}^{\infty}f_{j})\check{\Theta}w_{0}+\Delta^{2}w_{0})=0$ and, since
$\Delta^{2}=I-\check{\Theta}^{\ast}(\sum_{j=0}^{\infty}f_{j})\check{\Theta}$,
we have $e_{i}(w_{0}-\check{\Theta}^{\ast}(\sum_{j=0}^{n-1}f_{j})\check
{\Theta}w_{0})=0$. But also $(\sum_{j=0}^{n-1}f_{j})\check{\Theta}w_{0}%
=(\sum_{j=0}^{n-1}f_{j})\check{\Theta}(\sum_{k=0}^{n-1}e_{k})w_{0}$ and we get
the following equation, for every $i<n$,
\begin{equation}
e_{i}w_{0}-e_{i}\check{\Theta}^{\ast}(\sum_{j=0}^{n-1}f_{j})\check{\Theta
}(\sum_{k=0}^{n-1}e_{k})w_{0}=0. \label{w}%
\end{equation}
Setting $n=1$ and $i=0$ we obtain in particular the equation $e_{0}w_{0}%
=e_{0}\check{\Theta}^{\ast}f_{0}\check{\Theta}e_{0}w_{0}$. Since $\Theta$ is
assumed to be pure, $e_{0}w_{0}=0$. Now set $n=2$ and $i=1$ in equation
(\ref{w}) and use the fact that $f_{0}\check{\Theta}w_{0}=f_{0}\check{\Theta
}e_{0}=0$ to conclude that
\begin{equation}
e_{1}w_{0}=e_{1}\check{\Theta}^{\ast}f_{1}\check{\Theta}e_{1}w_{0}. \label{ew}%
\end{equation}
In order to \textquotedblleft bootstrap\textquotedblright\ purity to this
equation we first fix $\zeta\in E$ and, using Remark~\ref{rem}, we compute
\[
(T_{\zeta}^{\ast}\otimes I)e_{1}w_{0}=(T_{\zeta}^{\ast}\otimes I)e_{1}%
\check{\Theta}^{\ast}f_{1}\check{\Theta}e_{1}w_{0}=e_{0}\check{\Theta}^{\ast
}f_{0}(T_{\zeta}^{\ast}\otimes I)\check{\Theta}e_{1}w_{0}=
\]%
\[
=e_{0}\check{\Theta}^{\ast}f_{0}\check{\Theta}(T_{\zeta}^{\ast}\otimes
I)e_{1}w_{0}.
\]
Now we can appeal to the purity of $\Theta$ to conclude that $(T_{\zeta}%
^{\ast}\otimes I)e_{1}w_{0}=0$. Since this holds for all $\zeta\in E$,
$e_{1}w_{0}=0$. Continuing in this way we see that $e_{n}w_{0}=0$ for all
$n\geq0$. Thus $w_{0}=0$ and, consequently, $x=0$.
\end{proof}

\begin{lemma}
\label{tech}Let $\Theta$ be a characteristic function that is pure and
predictable and adopt the notation from Lemma \ref{charep}. For $i\geq1$ set
\[
\mathcal{K}_{i}:=\overline{span}\{V_{\Theta}(\xi)h\mid\xi\in E^{\otimes
i},\;h\in H(\Theta)\;\}
\]
and for $j\geq0$ set%
\[
\mathcal{M}_{j}:=\{\check{\Theta}x+\Delta_{\check{\Theta}}x:x\in E^{\otimes
j}\otimes\mathcal{E}_{1}\}\text{,}%
\]
where, for $j=0$, $E^{\otimes0}\otimes\mathcal{E}_{1}$ is $\mathcal{E}_{1}$.
Then,
\[
\mathcal{M}_{0}=(I_{K(\Theta)}-P_{\Theta})(\mathcal{K}_{1}).
\]

\end{lemma}

\begin{proof}
As usual, write $\Delta$ for $(I-\check{\Theta}^{\ast}\check{\Theta})^{1/2} $.
First we note that the map taking $x\in\mathcal{F}(E)\otimes\mathcal{E}_{1}$
to $\check{\Theta}x+\Delta x\in K(\Theta)$ is an isometry, since
$\check{\Theta}^{\ast}\check{\Theta}+\Delta^{2}=I$, and, consequently, that
for $i\neq j$, $\mathcal{M}_{i}$ is orthogonal to $\mathcal{M}_{j}$. Also, we
note that for $x\in E^{\otimes j}\otimes\mathcal{E}_{1}$ and $\xi\in E$,
$V_{\Theta}(\xi)(\check{\Theta}x+\Delta x)=(T_{\xi}\otimes I)\check{\Theta
}x+S_{\Theta}(\xi)\Delta x=\check{\Theta}(T_{\xi}\otimes I)x+\Delta(T_{\xi
}\otimes I)x$. Hence $V_{\Theta}(E)\mathcal{M}_{j}\subseteq\mathcal{M}_{j+1}$,
where we abbreviate $\overline{span}\{V_{\Theta}(\xi)x\mid\xi\in
E,\ x\in\mathcal{M}_{j}\}$ by $V_{\Theta}(E)\mathcal{M}_{j}$. It is also clear
that $V_{\Theta}(E)\mathcal{K}_{i}\subseteq\mathcal{K}_{i+1}$.

Next we show that for $j\geq1$, $\mathcal{K}_{1}$ is orthogonal to
$\mathcal{M}_{j}$. Indeed, let $j\geq1$, let $\zeta\in E$, let $\theta\in
E^{\otimes(j-1)}$ and let $h\in\mathcal{E}_{1}$. Then, for $\xi\in E$, we have
$V(\xi)^{\ast}(\check{\Theta}(\zeta\otimes\theta\otimes h)+\Delta(\zeta
\otimes\theta\otimes h))=(T_{\xi}^{\ast}\otimes I)\check{\Theta}(\zeta
\otimes\theta\otimes h)+S_{\Theta}(\xi)^{\ast}\Delta(\zeta\otimes\theta\otimes
h)$. Using Remark~\ref{rem} and the fact that $\Delta(\zeta\otimes
\theta\otimes h)=\Delta(T_{\zeta}\otimes I)(\theta\otimes h)=S_{\Theta}%
(\zeta)\Delta(\theta\otimes h)$ we find that $V(\xi)^{\ast}(\check{\Theta
}(\zeta\otimes\theta\otimes h)+\Delta(\zeta\otimes\theta\otimes h))=\check
{\Theta}(T_{\xi}^{\ast}\otimes I)(\zeta\otimes\theta\otimes h)+S_{\Theta}%
(\xi)^{\ast}S_{\Theta}(\zeta)\Delta(\theta\otimes h)=\check{\Theta}(\langle
\xi,\zeta\rangle\theta\otimes h)+\Delta(\langle\xi,\zeta\rangle\theta\otimes
h)\in H(\Theta)^{\perp}$.

It follows that $\mathcal{K}_{1}$ is orthogonal to $\mathcal{M}_{j}$, $j\geq
1$. Since $\mathcal{M}_{j}=(I-P_{\Theta})(\mathcal{M}_{j})$, we conclude that
$(I-P_{\Theta})\mathcal{K}_{1}$ is orthogonal to $\mathcal{M}_{j}$ for all
$j\geq1$. But it is also orthogonal to $H(\Theta)$ and we have $K(\Theta
)=H(\Theta)\oplus\sum_{j=0}^{\infty}\oplus\mathcal{M}_{j}$. Thus
\begin{equation}
(I-P_{\Theta})(\mathcal{K}_{1})\subseteq\mathcal{M}_{0}. \label{dir}%
\end{equation}
>From (\ref{dir}) it follows that $\mathcal{K}_{1}\subseteq\mathcal{M}%
_{0}\oplus H(\Theta)$. Applying $V_{\Theta}(E)$ to this we find that
$\mathcal{K}_{2}\subseteq\mathcal{M}_{1}\oplus\mathcal{K}_{1}$. A second
application of $V_{\Theta}(E)$ yields $\mathcal{K}_{3}\subseteq(\mathcal{M}%
_{2}\oplus\mathcal{M}_{1})+\mathcal{K}_{1}$. Continuing by induction we find
that for every $i\geq2$,
\begin{equation}
\mathcal{K}_{i}\subseteq\mathcal{K}_{1}+\sum_{j=1}^{i-1}\oplus\mathcal{M}_{j}.
\label{km}%
\end{equation}
Now suppose $y\in\mathcal{M}_{0}\ominus(I-P_{\Theta})(\mathcal{K}_{1})$. Then
$y=(I-P_{\Theta})y\in\mathcal{K}_{1}^{\perp}$. Since $y\in\mathcal{M}_{0}$,
$y$ is also orthogonal to $\mathcal{M}_{j}$ for every $j\geq1$. By (\ref{km}),
$y$ is orthogonal to $\mathcal{K}_{i}$ for every $i\geq1$. But $y\in
H(\Theta)^{\perp}$ and, by the minimality of $(V_{\Theta},\rho_{\Theta})$,
$H(\Theta)+\sum\mathcal{K}_{i}$ is dense in $K(\Theta)$. Thus $y=0$ and this,
combined with the inclusion (\ref{dir})) completes the proof.
\end{proof}

\begin{lemma}
\label{ei}Let $(\Theta,\mathcal{E}_{1},\mathcal{E}_{2},\tau_{1},\tau_{2})$ be
a pure and predictable characteristic function, let $(T,\sigma)=(T_{\Theta
},\sigma_{\Theta})$ be its canonical model acting on $H=H(\Theta)$, and let
$\mathcal{D}$ and $\mathcal{D}_{\ast}$ be the defect spaces associated with
$(T,\sigma)$. Then:

\begin{enumerate}
\item[(i)] The spaces $\mathcal{E}_{1}$ and $\mathcal{D}$ are isomorphic as
left $M$-modules;i.e. there is a unitary operator $W_{1}:\mathcal{E}%
_{1}\rightarrow\mathcal{D}$ such that, for every $a\in M$,
\[
W_{1}\tau_{1}(a)=(\varphi(a)\otimes I_{H})W_{1}.
\]

\item[(ii)] The spaces $\mathcal{E}_{2}$ and $\mathcal{D}_{*}$ are isomorphic
as left $M$-modules ; i.e. there is a unitary operator $W_{2} : \mathcal{E}%
_{2} \rightarrow\mathcal{D}_{*} $ such that, for every $a\in M$,
\[
W_{2}\tau_{2}(a)= \sigma(a)W_{2} .
\]

\end{enumerate}
\end{lemma}

\begin{proof}
Write $(V,\rho)$ for the minimal isometric dilation of $(T,\sigma)$ as
constructed in (\ref{vxi}) and the discussion preceding it. The representation
space of $(V,\rho)$ is $K=H\oplus(\mathcal{F}(E)\otimes_{\sigma_{1}%
}\mathcal{D})$. From the uniqueness of the minimal isometric dilation
\cite[Proposition 3.2]{MS98} and Lemma \ref{isomdil}, it follows that there is
a unitary operator $W:K(\Theta)\rightarrow K$ such that $W$ maps $H(\Theta)$
onto $H$ and satisfies the equations $V(\xi)W=WV_{\Theta}(\xi)$, $\xi\in E$,
and $\rho(a)W=W\rho_{\Theta}(a)$, $a\in M$. Write $W_{1}h=W(\check{\Theta
}h+\Delta h)$ for $h\in\mathcal{E}_{1}$, where $\Delta:=(I-\check{\Theta
}^{\ast}\check{\Theta})^{1/2}$. Then, in the notation of Lemma~\ref{tech},
$W_{1}(\mathcal{E}_{1})=W\mathcal{M}_{0}=W(I-P_{H(\Theta)})\mathcal{K}%
_{1}=W(I-P_{H(\Theta)})V_{\Theta}(E)H(\Theta)=(I-P_{H(\Theta)})WV_{\Theta
}(E)W^{\ast}WH(\Theta)=(I-P_{H(\Theta)})V(E)H=\mathcal{D}$, where the last
equality follows from equation (\ref{vxi}). Recall that the map $x\mapsto
\check{\Theta}x+\Delta x$ is an isometry defined on $\mathcal{F}%
(E)\otimes\mathcal{E}_{1}$. Hence $W_{1}$ is indeed a unitary operator from
$\mathcal{E}_{1}$ onto $\mathcal{D}$. Now fix $a\in M$ and $h\in
\mathcal{E}_{1}$ and recall that $\mathcal{D}\subseteq E\otimes H$ and
$\rho(a)|\mathcal{D}=(\varphi(a)\otimes I_{H(S)})|\mathcal{D}$. We have
\begin{multline*}
(\varphi(a)\otimes I_{H})W_{1}h=\rho(a)W(\check{\Theta}h+\Delta h)=W\rho
_{\Theta}(a)(\check{\Theta}h+\Delta h)\\
=W((\varphi_{\infty}(a)\otimes I)\check{\Theta}h+\Delta\tau_{1}(a)h)=W(\check
{\Theta}\tau_{1}(a)h+\Delta\tau_{1}(a)h)=W_{1}\tau_{1}(a)h.
\end{multline*}
This proves (i).

To prove the other assertion, recall first from Lemma~\ref{Q0} that $K_{0}$ is
the range of the projection $I-\tilde{V}\tilde{V}^{\ast}$ (in fact, we can
write $K_{0}=K\ominus V(E)K$) and there is an isometry $u$ from $K_{0}$ onto
$\mathcal{D}_{\ast}$. Note that we may view $\mathcal{E}_{2}$ as the first
summand of $\mathcal{F}(E)\otimes\mathcal{E}_{2}$ and that when we do, we can
write $\mathcal{E}_{2}=(\mathcal{F}(E)\otimes\mathcal{E}_{2})\ominus
\overline{span}\{(T_{\xi}\otimes I)(\mathcal{F}(E)\otimes\mathcal{E}_{2}%
)\mid\xi\in E\}$. Since $S_{\Theta}(E)\Delta(\mathcal{F}(E)\otimes
\mathcal{E}_{1})=\Delta((\mathcal{F}(E)\otimes\mathcal{E}_{1})\ominus
\mathcal{E}_{1})=\Delta(\mathcal{F}(E)\otimes\mathcal{E}_{1})$, we have
$\mathcal{E}_{2}=K(\Theta)\ominus V_{\Theta}(E)K(\Theta)=W^{\ast}K\ominus
W^{\ast}V(E)WW^{\ast}K=W^{\ast}(K\ominus V(E)K)=W^{\ast}K_{0}$. Thus, setting
$W_{2}=uW|\mathcal{E}_{2}$, we obtain a unitary operator from $\mathcal{E}%
_{2}$ onto $\mathcal{D}_{\ast}$. Finally, for $a\in M$ and $h\in
\mathcal{E}_{2}\subseteq K(\Theta)$,
\[
W_{2}\tau_{2}(a)h=uW\rho_{\Theta}(a)h=u\rho(a)Wh=\sigma(a)W_{2}h
\]
where the last equality follows from Lemma~\ref{Q0} (iii).
\end{proof}

\begin{theorem}
\label{isomS}Let $(\Theta,\mathcal{E}_{1},\mathcal{E}_{2},\tau_{1},\tau_{2})$
be a pure and predictable characteristic function and let $(T,\sigma
)=(T_{\Theta},\sigma_{\Theta})$ on $H:=H(\Theta)$ be the associated canonical
model. Then this representation is c.n.c and its characteristic function
$(\hat{\Theta}_{T},\mathcal{D},\mathcal{D}_{\ast},(\varphi\otimes
I_{H})|\mathcal{D},\sigma|\mathcal{D}_{\ast})$ is isomorphic to $(\Theta
,\mathcal{E}_{1},\mathcal{E}_{2},\tau_{1},\tau_{2})$.
\end{theorem}

\begin{proof}
We continue with the notation of the proof of Lemma \ref{ei}. In particular,
$W_{1}$ will denote the Hilbert space isomorphism from $\mathcal{E}_{1}$ to
$\mathcal{D}$ constructed there, while $W_{2}$ will denote the Hilbert space
isomorphism from $\mathcal{E}_{2}$ to $\mathcal{D}_{\ast}$. Also, $W$ will be
the unitary operator from $K(\Theta)$ onto $K$, where $K$ is the space of the
minimal isometric dilation $(V,\rho)$ of $(T,\sigma)$ as in the proof of
Lemma~\ref{ei}. It is shown there that $W$ maps $\mathcal{E}_{2}$ onto $K_{0}$
and it intertwines $V_{\Theta}$ and $V$. Thus it maps $\mathcal{F}%
(E)\otimes\mathcal{E}_{2}$ onto $Q_{\infty}(K)$. Since $W(H(\Theta))=H$,
$H\cap P_{\infty}(K)=W(H(\Theta)\cap\overline{\Delta(\mathcal{F}%
(E)\otimes\mathcal{E}_{1})}$. But if $y\in H(\Theta)\cap\overline
{\Delta(\mathcal{F}(E)\otimes\mathcal{E}_{1})}$ then, for every $x\in
\mathcal{F}(E)\otimes\mathcal{E}_{1}$, $y$ is orthogonal to $\check{\Theta
}x+\Delta x$ and also $y$ is orthogonal to $\check{\Theta}x\in\mathcal{F}%
(E)\otimes\mathcal{E}_{2}$. Thus $y$ is orthogonal to $\Delta x$ for every
such $x$ and it follows that $y=0$. Hence $H\cap P_{\infty}(K)=\{0\}$ and,
consequently, $(T,\sigma)$ is a c.n.c. representation.

Since $W$ maps $\mathcal{F}(E)\otimes\mathcal{E}_{2}$ onto $Q_{\infty}(K)$, it
follows that $Q_{\infty}W\Delta x=0$ for $x\in\mathcal{F}(E)\otimes
\mathcal{E}_{1}$. Also, recall that $W_{2}=uW|\mathcal{E}_{2}$ and, for
$\xi\in\mathcal{F}(E)$ and $h\in\mathcal{E}_{2}$ we have $(I_{\mathcal{F}%
(E)}\otimes W_{2})(\xi\otimes h)=\xi\otimes uWh=(I\otimes u)(\xi\otimes
Wh)=(I\otimes u)W_{\infty}^{\ast}Q_{\infty}V(\xi)Wh=(I\otimes u)W_{\infty
}^{\ast}Q_{\infty}W(\xi\otimes h)$. Thus
\[
I_{\mathcal{F}(E)}\otimes W_{2}=(I_{\mathcal{F}(E)}\otimes u)W_{\infty}^{\ast
}Q_{\infty}W.
\]
So from the definition of $\Theta_{T}$, Definition \ref{charop}, we find that
for every $h\in\mathcal{E}_{1}$,
\begin{align*}
\Theta_{T}W_{1}h  &  =\Theta_{T}W(\check{\Theta}h+\Delta h)=\Theta_{T}%
W\check{\Theta}h\\
&  =(I_{\mathcal{F}(E)}\otimes u)W_{\infty}^{\ast}Q_{\infty}W\check{\Theta
}h=(I_{\mathcal{F}(E)}\otimes W_{2})\check{\Theta}h\text{.}%
\end{align*}
Hence, for $\xi\otimes d\in\mathcal{F}(E)\otimes\mathcal{D}$ and
$h:=W_{1}^{\ast}d\in\mathcal{E}_{1}$, we have $(I_{\mathcal{F}(E)}\otimes
W_{2})\check{\Theta}(I_{\mathcal{F}(E)}\otimes W_{1}^{\ast})(\xi\otimes
d)=(I_{\mathcal{F}(E)}\otimes W_{2})\check{\Theta}(\xi\otimes
h)=(I_{\mathcal{F}(E)}\otimes W_{2})(T_{\xi}\otimes I_{\mathcal{E}_{2}}%
)\check{\Theta}h=(T_{\xi}\otimes I_{\mathcal{D}_{\ast}})(I_{\mathcal{F}%
(E)}\otimes W_{2})\check{\Theta}h=(T_{\xi}\otimes I_{\mathcal{D}_{\ast}%
})\Theta_{T}W_{1}h=\Theta_{T}(T_{\xi}\otimes I_{\mathcal{D}})d=\Theta_{T}%
(\xi\otimes d)$. Therefore
\[
(I_{\mathcal{F}(E)}\otimes W_{2})\check{\Theta}(I_{\mathcal{F}(E)}\otimes
W_{1}^{\ast})=\Theta_{T}\text{,}%
\]
as was to be proved.
\end{proof}

\section{Commutants of Models}

In \cite[Theorem 4.4]{MS98} we proved a commutant lifting theorem for
completely contractive representations of tensor algebras. The analysis there
extends without difficulty to $\sigma$-weakly continuous representations of
Hardy algebras. However, with the analysis in \cite{MS04} available to us and
the results of the preceding section, it is possible to give a refined version
of the commutant lifting theorem, at least in the context of $C_{\cdot0}$
representations. The theorem we shall prove in this section generalizes
Theorem 6.1 of \cite{gP89}.

First recall that if $(T,\sigma)$ is a $C_{\cdot0}$ representation of $E$ on a
Hilbert space $H$, if $\Theta=\hat{\Theta}_{T}$ is the characteristic function
associated to the characteristic operator $(\Theta_{T},\mathcal{D}%
,\mathcal{D}_{\ast},\tau_{1},\tau_{2})$, and if $(T_{\Theta},\sigma_{\Theta})$
is the canonical model built from $\Theta$, then the Hilbert space of the
minimal isometric dilation of $(T_{\Theta},\sigma_{\Theta})$, $K(\Theta)$, is
$\mathcal{F}(E)\otimes_{\tau_{2}}\mathcal{D}_{\ast}$, by virtue of Theorems
\ref{InnerCriterion} and \ref{TSisT}. (A bit more completely, Theorem
\ref{InnerCriterion} guarantees that $\hat{\Theta}_{T}$ is inner if
$(T,\sigma)$ is $C_{\cdot0}$. Also, Lemma \ref{C0} guarantees that the minimal
isometric dilation of $(T,\sigma)$ is an induced representation if (and only
if) $(T,\sigma)$ is $C_{\cdot0}$. And, Theorem \ref{TSisT} identifies the form
of that induced representation.) The model space $H(\Theta)$ is $(\mathcal{F}%
(E)\otimes_{\tau_{2}}\mathcal{D}_{\ast})\ominus\Theta_{T}(\mathcal{F}%
(E)\otimes_{\tau_{1}}\mathcal{D})$ in this case. Recall, too, that
$(\mathcal{G},\tau)$ is a fixed supplement of $\tau_{1}$ and $\tau_{2}$ and
that $\mathcal{F}(E)\otimes_{\tau}\mathcal{G}$ decomposes as $\mathcal{F}%
(E)\otimes_{\tau}\mathcal{G}=(\mathcal{F}(E)\otimes_{\pi_{0}}H_{0}%
)\oplus(\mathcal{F}(E)\otimes_{\tau_{1}}\mathcal{D})\oplus(\mathcal{F}%
(E)\otimes_{\tau_{2}}\mathcal{D}_{\ast})$\ (equation (\ref{DecompFEG})). A
moment's reflection reveals that if $v_{2}$ is the isometric embedding of
$\mathcal{D}_{\ast}$ in $\mathcal{G}$ that sends $d_{\ast}$ in $\mathcal{D}%
_{\ast}$ to $(0,0,d_{\ast})^{tr}$, then $I\otimes v_{2}$ is an isometric
embedding of $\mathcal{F}(E)\otimes_{\tau_{2}}\mathcal{D}_{\ast}$ in
$\mathcal{F}(E)\otimes_{\tau}\mathcal{G}$ that intertwines the two induced
representations of $H^{\infty}(E)$ and that maps $H(\Theta)$ onto the space%
\[
(\mathcal{F}(E)\otimes_{\tau}\mathcal{G})\ominus\Theta_{T}(\mathcal{F}%
(E)\otimes_{\tau}\mathcal{G})\text{,}%
\]
where here $\Theta_{T}$ is treated as the matrix in equation (\ref{MatrixFEG}%
). On the other hand, the canonical equivalence $\Phi$ from $K(\Theta
)=\mathcal{F}(E)\otimes_{\tau_{2}}\mathcal{D}_{\ast}$ to the Hilbert space $K$
of the minimal isometric dilation $(V,\rho)$ of $(T,\sigma)$ is a Hilbert
space isomorphism that intertwines $V\times\rho$ and the induced
representation $\tau_{2}^{\mathcal{F}(E)}$, maps $H(\Theta) $ onto $H$ and
implements a unitary equivalence between $(T,\sigma)$ and $(T_{\Theta}%
,\sigma_{\Theta})$ (see Theorem \ref{TSisT}). Hence, if $U:\mathcal{F}%
(E)\otimes_{\tau}\mathcal{G}\rightarrow\mathcal{F}(E^{\tau})\otimes_{\iota
}\mathcal{G}$ is the Fourier transform from Remark \ref{FunctionalRep}, and if
$U_{0}$ is the composition $U_{0}:=U(I\otimes v_{2})(\Phi^{-1}|H)$, i.e., if
$U_{0}$ is built from the following diagram%
\[
H\subseteq K\overset{\Phi^{-1}}{\longrightarrow}K(\Theta)\overset{I\otimes
v_{2}}{\longrightarrow}\mathcal{F}(E)\otimes_{\tau}\mathcal{G}\overset
{U}{\longrightarrow}\mathcal{F}(E^{\tau})\otimes_{\iota}\mathcal{G}\text{,}%
\]
then $U_{0}$ is an isometry mapping $H$ into $\mathcal{F}(E^{\tau}%
)\otimes_{\iota}\mathcal{G}$ and has the property that for every $\Xi\in
H^{\infty}(E^{\tau})$, $U_{0}^{\ast}(\Xi\otimes I_{\mathcal{G}})U_{0}$
commutes with $T\times\sigma(H^{\infty}(E))$.

\begin{theorem}
\label{commutant}Let $\pi$ be a completely contractive $\sigma$-weakly
continuous representation of $H^{\infty}(E)$ on the Hilbert space $H$ such
that the associated covariant representation of $E$, $(T,\sigma)$, is a
$C_{\cdot0}$-representation. Let $U_{0}:H\rightarrow\mathcal{F}(E^{\tau
})\otimes_{\iota}\mathcal{G}$ be the isometric embedding just described. Then
for every $X\in B(H)$ that commutes with $\pi(H^{\infty}(E))$, there is an
$\Xi\in H^{\infty}(E^{\tau})$ such that

\begin{enumerate}
\item[(i)] $\Vert\Xi\Vert=\Vert X \Vert$, and

\item[(ii)] $X=U_{0}^{\ast}(\Xi\otimes I_{\mathcal{G}})U_{0}$.
\end{enumerate}
\end{theorem}

\begin{proof}
We have already noted that every $X$ of the form in (ii) commutes with
$\pi(H^{\infty}(E))$ and of course $\left\Vert X\right\Vert \leq\left\Vert
\Xi\otimes I_{\mathcal{E}}\right\Vert =\left\Vert \Xi\right\Vert $ since
$\iota^{\mathcal{F}(E^{\tau})}$ is faithful by Remark \ref{faithful}. But the
converse results from \cite[Theorem 4.4]{MS98} as follows. Given $X\in B(H)$
that commutes with $\pi(H^{\infty}(E))$, Theorem 4.4 of \cite{MS98} produces
an operator $Y$ on the Hilbert space $K$ of the minimal isometric dilation
$(V,\rho)$ of $(T,\sigma)$ that commutes with $(V,\rho)$, satisfies the
equation $\left\Vert Y\right\Vert =\left\Vert X\right\Vert $ and satisfies the
equation $X=P_{H}Y|H$. Since $(T,\sigma)$ is $C_{\cdot0}$, Lemma \ref{C0}
implies that $(V,\rho)$ is an induced representation. Theorem \ref{TSisT}
identifies the structure of that induced representation and shows that $\Phi$
implements an equivalence between $(V,\rho)$ and the (covariant)
representation $\tau_{2}^{\mathcal{F}(E)}$. The map $I\otimes v_{2}$ embeds
$\mathcal{F}(E)\otimes_{\tau_{2}}\mathcal{D}_{\ast}$ into $\mathcal{F}%
(E)\otimes_{\tau}\mathcal{G}$ in such a way that $(I\otimes v_{2})\Phi
^{-1}(Y)\Phi(I\otimes v_{2})^{\ast}$ commutes with $\tau^{\mathcal{F}%
(E)}(H^{\infty}(E))$. So, since $U$ is the Fourier transform from
$\mathcal{F}(E)\otimes_{\tau}\mathcal{G}$ to $\mathcal{F}(E^{\tau}%
)\otimes_{\iota}\mathcal{G}$, Theorem \ref{comminduced} guarantees that
$U(I\otimes v_{2})\Phi^{-1}(Y)\Phi(I\otimes v_{2})^{\ast}U^{\ast}$ is an
operator on $\mathcal{F}(E^{\tau})\otimes_{\iota}\mathcal{G}$ that lies in
$\iota^{\mathcal{F}(E^{\tau})}(H^{\infty}(E^{\tau}))$, i.e., $U(I\otimes
v_{2})\Phi^{-1}(Y)\Phi(I\otimes v_{2})^{\ast}U^{\ast}=\Xi\otimes
I_{\mathcal{G}}$ for a $\Xi\in H^{\infty}(E^{\tau})$. Hence, as a calculation
reveals, $U_{0}^{\ast}(\Xi\otimes I_{\mathcal{G}})U_{0}=X$ and $\left\Vert
\Xi\right\Vert \leq\left\Vert Y\right\Vert =\left\Vert X\right\Vert $.
\end{proof}

\begin{remark}
\label{SarasonCompare}If $M=\mathbb{C}=E$, and if $(T,\sigma)$ is a
$C_{\cdot0}$ representation with $1$-dimensional defect spaces, then Theorem
\ref{commutant} gives Sarason's original commutant lifting theorem \cite{dS67}.
\end{remark}

\section{Invariant Subspaces}

In the theory of models for single operators, invariant subspaces are
determined by factorizations of the characteristic operator functions. The
same is true in our setting. To keep the presentation as simple as possible,
we shall restrict our attention to $C_{\cdot0}$ representations. We shall need
to consider factorizations, i.e., compositions, $\Theta=\Theta_{1}\Theta_{2}$,
where $\Theta$ is the necessarily inner characteristic function associated
with a $C_{\cdot0}$-representation and where each $\Theta_{i}$, $i=1,2$, is an
inner characteristic function that is \emph{not} necessarily purely
contractive. Two such compositions $\Theta=\Theta_{1}\Theta_{2}=\Theta
_{1}^{\prime}\Theta_{2}^{\prime}$ are said to be \emph{equivalent} if
$\Theta_{1}^{\prime}=\Theta_{1}(I\otimes V_{0})$ and $\Theta_{2}^{\prime
}=(I\otimes V_{0}^{\ast})\Theta_{2}$ for a suitable unitary operator $V_{0}$.

\begin{theorem}
\label{invariant}Let $(T,\sigma)$ a $C_{\cdot0}$-representation of $E$ on $H$,
with $T\times\sigma$ denoting the associated representation of $H^{\infty}%
(E)$, and let $\Theta:=\hat{\Theta}_{T}$ be the inner characteristic function
of this representation. Then there is a bijection between the subspaces of $H$
that are invariant under $(T\times\sigma)(H^{\infty}(E))$ and equivalence
classes of factorizations $\Theta=\Theta_{1}\Theta_{2}$ of $\Theta$ as a
composition of two inner characteristic functions.
\end{theorem}

\begin{proof}
By Theorem \ref{TSisT}, we may assume that $(T,\sigma)$ is $(T_{\Theta}%
,\sigma_{\Theta})$ for the inner characteristic function $(\Theta
,\mathcal{D},\mathcal{D}_{\ast},\tau_{1},\tau_{2})$. Hence, the space $H$ is
$H(\Theta)=(\mathcal{F}(E)\otimes\mathcal{D}_{\ast})\ominus\Theta
(\mathcal{F}(E)\otimes\mathcal{D})$.

Fix a subspace $\mathcal{M}\subseteq H(\Theta)$ that is invariant under
$(T\times\sigma)(H^{\infty}(E))$; that is, for every $\xi\in E$ and $a\in M$,
$T_{\Theta}(\xi)\mathcal{M}\subseteq\mathcal{M}$ and $\sigma_{\Theta
}(a)\mathcal{M}\subseteq\mathcal{M}$. Write $\mathcal{N}=\mathcal{M}%
\oplus\Theta(\mathcal{F}(E)\otimes\mathcal{D})\subseteq K(\Theta)$. Recall
that $T_{\Theta}(\xi)$ (for $\xi\in E$) and $\sigma_{\Theta}(a)$ (for $a\in M
$) are the compressions of $T_{\xi}\otimes I_{\mathcal{D}_{\ast}}$ and
$\varphi_{\infty}(a)\otimes I_{\mathcal{D}_{\ast}}$, respectively, to
$H(\Theta)$. Also recall that $T_{\xi}\otimes I_{\mathcal{D}_{\ast}}$ and
$\varphi_{\infty}(a)\otimes I_{\mathcal{D}_{\ast}}$ leave $\Theta
(\mathcal{F}(E)\otimes\mathcal{D})$ invariant. It follows that $\mathcal{N}$
is invariant under these operators. Thus, defining $S(\xi)$ and $\pi(a)$ (for
$\xi\in E$ and $a\in M$) to be the restrictions of $T_{\xi}\otimes
I_{\mathcal{D}_{\ast}}$ and $\varphi_{\infty}(a)\otimes I_{\mathcal{D}_{\ast}%
}$, respectively, to $\mathcal{N}$, we get an isometric representation of $E $
on $\mathcal{N}$. Since this is the restriction of a pure representation in
the sense of \cite{MS99}, meaning that condition (ii) of Lemma \ref{C0} is
satisfied, it is also pure. It follows from the equivalence of (ii) and (iv)
in Lemma \ref{C0} $(S,\pi)$ is induced. That is, there is a representation
$\rho$ of $M$ on a Hilbert space $H_{0}$ such that $(S,\pi)$ is unitarily
equivalent to the induced representation on $\mathcal{F}(E)\otimes_{\rho}%
H_{0}$. Hence, there is a unitary operator $\Theta_{1}$ from $\mathcal{F}%
(E)\otimes_{\rho}H_{0}$ onto $\mathcal{N}$ intertwining the induced
representation and $(S,\pi)$. It is then easy to check that $(\Theta_{1}%
,H_{0},\mathcal{D}_{\ast},\rho,\tau_{2})$ is an inner characteristic function.
(Recall that it is not assumed to be purely contractive).

We now write $\Theta_{2}=\Theta_{1}^{\ast}\Theta:\mathcal{F}(E)\otimes
_{\tau_{1}}\mathcal{D}\rightarrow\mathcal{F}(E)\otimes_{\rho}H_{0}$. Clearly,
$\Theta_{2}$ is an isometry (note that the range of $\Theta$ is contained in
the range of $\Theta_{1}$) and since $\Theta_{2}$ evidently intertwines
$\rho^{\mathcal{F}(E)}$ and $\tau_{1}^{\mathcal{F}(E)}$, we see that
$(\Theta_{2},\mathcal{D},H_{0},\tau_{1},\rho)$ is an inner characteristic
function (where, again, we do not assume that it is purely contractive).We
have $\Theta=\Theta_{1}\Theta_{2}$.

So far, starting with an invariant subspace $\mathcal{M}$ of $H(\Theta)$, we
obtained a factorization of $\Theta$. Note also that
\begin{equation}
\mathcal{F}(E)\otimes_{\tau_{2}}\mathcal{D}_{\ast}\ominus\Theta_{1}%
(\mathcal{F}(E)\otimes_{\rho}H_{0})=H(\Theta)\ominus\mathcal{M}. \label{M}%
\end{equation}

Now assume that $(\Theta_{1},H_{0},\mathcal{D}_{\ast},\rho,\tau_{2})$ and
$(\Theta_{2},\mathcal{D},H_{0},\tau_{1},\rho)$ are two characteristic
functions (not necessarily purely contractive) such that $\Theta=\Theta
_{1}\Theta_{2}$. Clearly $\Theta(\mathcal{F}(E)\otimes_{\tau_{1}}%
\mathcal{D})\subseteq\Theta_{1}(\mathcal{F}(E)\otimes_{\rho}H_{0})$. Set
\[
\mathcal{M}=\Theta_{1}(\mathcal{F}(E)\otimes_{\rho}H_{0})\ominus
\Theta(\mathcal{F}(E)\otimes_{\tau_{1}}\mathcal{D}).
\]
Then $\mathcal{M}\subseteq H(\Theta)$. Since $\mathcal{M}$ is clearly
invariant for $\sigma_{\Theta}(M)$, we need to show that it is invariant for
$T_{\Theta}(\xi)$ for $\xi\in E$. Fix an $h\in\mathcal{M}$ and $\xi\in E$.
Since $h$ is in the range of $\Theta_{1}$ and $\Theta_{1}$ intertwines
$T_{\xi}\otimes I_{H_{0}}$ and $T_{\xi}\otimes I_{\mathcal{D}_{\ast}}$,
$(T_{\xi}\otimes I_{\mathcal{D}_{\ast}})h$ is also in the range of $\Theta
_{1}$. Thus $T_{\Theta}(\xi)h=P_{H(\Theta)}(T_{\xi}\otimes I_{\mathcal{D}%
_{\ast}})h$ lies in $\mathcal{M}$. Hence $\mathcal{M} $ is an invariant
subspace of $H(\Theta)$. Note also that if we start with an equivalent
factorization $\Theta=\Theta_{1}^{\prime}\Theta_{2}^{\prime}$ we get the same
subspace $\mathcal{M}$.

It is clear from the decomposition (\ref{M}) that if we start with an
invariant subspace $\mathcal{M}$ and find the factorization $\Theta=\Theta
_{1}\Theta_{2}$ as above, then the invariant subspace associated to this
factorization is the space $\mathcal{M}$ we started with.

Now start with a factorization $\Theta=\Theta_{1}\Theta_{2}$ and associate
with it the subspace $\mathcal{M}=\Theta_{1}(\mathcal{F}(E)\otimes_{\rho}%
H_{0})\ominus\Theta(\mathcal{F}(E)\otimes_{\tau_{1}}\mathcal{D})$ as above. To
this subspace we apply the argument at the beginning of the proof to get a
factorization $\Theta=\Theta_{1}^{\prime}\Theta_{2}^{\prime}$. To do this, we
write $\mathcal{N}=\mathcal{M}\oplus\Theta(\mathcal{F}(E)\otimes_{\tau_{1}%
}\mathcal{D})$ ($=\Theta_{1}(\mathcal{F}(E)\otimes_{\rho}H_{0})$) and find a
representation $\rho^{\prime}$ on $H_{0}^{\prime}$ and a unitary operator
$\Theta_{1}^{\prime}:\mathcal{F}(E)\otimes_{\rho^{\prime}}H_{0}^{\prime
}\rightarrow\mathcal{N}$ that implements a unitary equivalence of the induced
representation on $\mathcal{F}(E)\otimes_{\rho^{\prime}}H_{0}^{\prime}$ and
the restriction to $\mathcal{N}$ of the induced representation on
$\mathcal{F}(E)\otimes_{\tau_{2}}\mathcal{D}_{\ast}$. Setting $V=\Theta
_{1}^{\ast}\Theta_{1}^{\prime}$ we get a unitary operator from $\mathcal{F}%
(E)\otimes_{\rho^{\prime}}H_{0}^{\prime}$ onto $\mathcal{F}(E)\otimes_{\rho
}H_{0}$ that intertwines the induced representations. It is easy to see that
such a unitary operator is of the form $I_{\mathcal{F}(E)}\otimes V_{0}$ for
some unitary operator $V_{0}$ from $H_{0}^{\prime}$ onto $H_{0}$ (roughly,
$V_{0}$ is the restriction of $V$ to $H_{0}^{\prime}$ viewed as the wandering
subspace of $\mathcal{F}(E)\otimes H_{0}^{\prime}$). We thus have $\Theta
_{1}(I_{\mathcal{F}(E)}\otimes V_{0})=\Theta_{1}^{\prime}$.
\end{proof}

\section{An Example: Analytic crossed products}

In this section we illustrate some of the results of the previous sections as
applied to the special case of correspondences induced from endomorphisms. We
shall fix an endomorphism $\alpha$ of a $W^{\ast}$-algebra $M$ and we shall
let $E$ be the $W^{\ast}$-correspondence $_{\alpha}M$. That is, as a (right)
$W^{\ast}$-module over $M$, $E$ is $M$ with the inner product defined by the
formula $\langle\xi_{1},\xi_{2}\rangle=\xi_{1}^{\ast}\xi_{2}$, $\xi_{1}$,
$\xi_{2}\in E$, but the left action is given by $\alpha$, i.e., $a\cdot
\xi\ (=\varphi(a)\xi):=\alpha(a)\xi$, for $\xi\in E$ and $a\in M$.

The associated Hardy algebra, $H^{\infty}(E)$, has a particularly attractive
description, which we shall develop. Note that for each $k\geq1$, the
correspondence $E^{\otimes k}$ can be identified with $_{\alpha^{k}}M$. The
map implementing the isomorphism takes $\xi_{1}\otimes\cdots\otimes\xi_{k}$ to
$\alpha^{k-1}(\xi_{1})\alpha^{k-2}(\xi_{2})\cdots\xi_{k}$. Thus $\mathcal{F}%
(E)$ can be identified with the direct sum $\sum_{k=0}^{\infty}\oplus
\,_{\alpha^{k}}M$ (where $\alpha^{0}$ is the identity map, and the zeroth
summand, $_{\alpha^{0}}M$, is simply $M$, viewed as the identity
correspondence from $M$ to $M$).

The action of $M$ on $\mathcal{F}(E)$ given in this form, $\varphi_{\infty} $,
now written $\alpha_{\infty}$, is familiar from the theory of crossed
products: for $a\in M$, $\alpha_{\infty}(a)(\xi_{k})=(\alpha^{k}(a)\xi_{k})$
for $(\xi_{k})\in\mathcal{F}(E)$. On the other hand for $\xi\in E$, the
creation operator is given by the formula $T_{\xi}(\xi_{k})=(\theta_{k})$
where $\theta_{k}=\alpha^{k-1}(\xi)\xi_{k-1}$. Note that since%
\[
T_{a\xi b}=\alpha_{\infty}(a)T_{\xi}\alpha_{\infty}(b)\text{,}%
\]
$a$, $b\in M$ and $\xi\in E$, the operators $T_{\xi}$ are completely
determined by $T_{1}$, where $1$ is the identity element of $M$ viewed as a
vector in $E$. Evidently, $T_{1}$ is a power partial isometry, and assuming
that $\alpha$ is unital, which we shall, $T_{1}$ is an isometry. We shall
write $w$ for $T_{1}$. Then $H^{\infty}(E)$ is simply the $\sigma$-weakly
closed subalgebra of the $W^{\ast}$-algebra $\mathcal{L}(\mathcal{F}(E))$
generated by $\alpha_{\infty}(M)$ and $w$. For historical reasons we shall
call this Hardy algebra \emph{the analytic crossed product} determined by $M$
and $\alpha$ and denote it by $M\rtimes_{\alpha}\mathbb{Z}_{+}$.

Non-self-adjoint algebras of this form (and closely related algebras) have a
long history going back to work of Kadison and Singer \cite{KS60} and Arveson
\cite{wA67a, wA67b}. In these papers and in most of the subsequent literature,
$\alpha$ is assumed to be an \emph{automorphism }of $M$. However, in
\cite{jP84}, Peters studied a related structure associated to an endomorphism
of a commutative $C^{\ast}$-algebra and proposed the name \emph{semi-crossed
products} for these. They turn out to be examples of tensor algebras and are
discussed from this point of view in \cite{MS98}. The term,
\emph{non-self-adjoint crossed product} was introduced in \cite{MMS78}, but
was changed to \emph{analytic crossed product} some years later in \cite{MS86}
to reflect better their function theoretic aspects. Since we are trying to
promote the view that all Hardy algebras as \emph{bona fide} spaces of
analytic functions, we shall adopt the term \textquotedblleft analytic crossed
product\textquotedblright\ to describe algebras of the form $M\rtimes_{\alpha
}\mathbb{Z}_{+}$.

Fix a (not-necessarily faithful) representation $\sigma$ of $M$ on the Hilbert
space $H$. Since $E^{\otimes n}$ may be identified with $_{\alpha^{n}}M$ for
all $n\geq0$, the spaces $E^{\otimes n}\otimes_{\sigma}H$ may each be
identified with $H$ via the Hilbert space isomorphism $W_{k}$ defined by the
formulae%
\begin{equation}
W_{k}(\xi_{1}\otimes\cdots\otimes\xi_{k}\otimes h)=\left\{
\begin{array}
[c]{cc}%
\sigma(\alpha^{k-1}(\xi_{1})\alpha^{k-2}(\xi_{2})\cdots\xi_{k})h\text{,} &
k>0\\
\sigma(\xi_{0})h & k=0
\end{array}
\right.  \text{,} \label{inducediso}%
\end{equation}
$\xi_{i}\in E$, $h\in H$. Then the direct sum $W:=\sum_{k\geq0}\oplus W_{k}$
is a Hilbert space isomorphism from $\mathcal{F}(E)\otimes_{\sigma}H$ onto
$\ell^{2}(\mathbb{Z}_{+},H)$, where $\ell^{2}(\mathbb{Z}_{+},H):=\{\xi
:\mathbb{Z}_{+}\rightarrow H\mid\sum_{k\geq0}\Vert\xi(k)\Vert^{2}<\infty\}$.
(It will be convenient below to indicate the dependence of $W$ and the $W_{k}$
on $\sigma$ by writing $W^{\sigma}$ and $W_{k}^{\sigma}$, but we omit this
until necessary.) Define a covariant representation of $E$ on $\ell
^{2}(\mathbb{Z}_{+},H)$, denoted $(S_{H},\psi_{H})$, by the equations
\[
(S_{H}(\xi)x)(k)=\sigma(\alpha^{k-1}(\xi))x(k-1),\;\xi\in E=\,_{\alpha
}M\text{, }x\in\ell^{2}(\mathbb{Z}_{+},H)
\]
and%
\[
(\psi_{H}(a)x)(k)=\sigma(\alpha^{k}(a))x(k),\;a\in M\text{, }x\in\ell
^{2}(\mathbb{Z}_{+},H)\text{.}%
\]
Thus, $S_{H}(1)$ is the unilateral shift (of appropriate multiplicity). Then a
moment's reflection using the definition of the representation induced by
$\sigma$, Definition \ref{inducedRep}, and equations (\ref{ind1}) and
(\ref{ind2}), reveals that $W$ implements a unitary equivalence between the
representation of $(M,~_{\alpha}M)$ induced by $\sigma$ and $(S_{H},\psi_{H}%
)$. That is
\[
W\sigma^{\mathcal{F}(_{\alpha}M)}(w)W^{\ast}=S_{H}(1)
\]
and%
\[
W\sigma^{\mathcal{F}(_{\alpha}M)}(\alpha_{\infty}(a))W^{\ast}=\psi
_{H}(a)\text{,}%
\]
$a\in M$.

Consider next an operator $R\in B(\ell^{2}(\mathbb{Z}_{+},H))$ that commutes
with the representation $S_{H}\times\psi_{H}(M\rtimes_{\alpha}\mathbb{Z}_{+}%
)$. Then since $R$ commutes with the shift $S_{H}(1)$, it is well known and
easy to verify that $R$ must be a block analytic Toeplitz operator. That is,
the matrix of $R$ with the direct sum decomposition of $\ell^{2}%
(\mathbb{Z}_{+},H)$ has this form:%
\begin{equation}
R=\left(
\begin{array}
[c]{ccccc}%
R_{0} & R_{1} & R_{2} &  & \cdots\\
0 & R_{0} & R_{1} & R_{2} & \\
0 & 0 & R_{0} & R_{1} & \ddots\\
& \ddots & \ddots & \ddots & \ddots\\
\vdots &  &  &  &
\end{array}
\right)  \text{,} \label{RToep}%
\end{equation}
where each $R_{k}\in B(H)$. On the other hand, since $R$ commutes with
$\psi_{H}(M)$, a straightforward calculation reveals that each $R_{k}$
satisfies the equation
\begin{equation}
\sigma(a)R_{k}=R_{k}\sigma(\alpha^{k}(a))\text{,} \label{Rcomm}%
\end{equation}
for all $a\in M$, i.e., $R_{k}$ intertwines $\sigma$ and $\sigma\circ
\varphi^{k}$. And conversely, every bounded operator $R$ on $\ell
^{2}(\mathbb{Z}_{+},H)$ whose matrix with respect to the direct sum
decomposition of $\ell^{2}(\mathbb{Z}_{+},H)$ is a block Toeplitz matrix, as
in equation (\ref{RToep}), whose entries satisfy equation (\ref{Rcomm}), must
commute with the image of $S_{H}\times\psi_{H}$.

Suppose now that $\sigma$ is faithful, so we may form the $\sigma$-dual of
$E=\,_{\alpha}M$ and note that $(E^{\sigma})^{\otimes k}$ is the $\sigma$-dual
correspondence of $E^{\otimes k}=\,_{\alpha^{k}}M$. Hence
\[
(E^{\sigma})^{\otimes k}=\{\eta:H\rightarrow\,_{\alpha^{k}}M\otimes H\mid
\eta\sigma(a)=(\alpha^{k}(a)\otimes I)\eta,a\in M\}.
\]
It follows from the definition of the maps $W_{k}$ in equation
(\ref{inducediso}) that
\[
W_{k}\cdot(E^{\sigma})^{\otimes k}:=\{W_{k}\eta\mid\eta\in(E^{\sigma
})^{\otimes k}\}=\{z\in B(H)\mid z\sigma(a)=\sigma(\alpha^{k}(a))z,\ a\in M\}
\]
Thus we have substantially proved the following proposition. We leave the
remaining details to the reader.

\begin{proposition}
\label{CommutantReal1}Suppose $E=\,_{\alpha}M$, for an endomorphism $\alpha$
of $M$, and that $\sigma$ is a faithful representation of $M$ on the Hilbert
space $H$. If $W=\sum_{k\geq0}\oplus W_{k}$ is the Hilbert space isomorphism
from $\mathcal{F}(E)\otimes_{\sigma}H$ to $\ell^{2}(\mathbb{Z}_{+},H)$, where
the $W_{k}$ are defined in equation (\ref{inducediso}) and if $U:\mathcal{F}%
(E)\otimes_{\sigma}H\rightarrow\mathcal{F}(E^{\sigma})\otimes_{\iota}H$ is the
Fourier transform determined by $\sigma$, then for all $\eta\in(E^{\sigma
})^{\otimes k}$, $\sigma(a)W_{k}\eta=\sigma(\alpha^{k}(a))W_{k}\eta$, for all
$a\in M$, and%
\[
WU^{\ast}(T_{\eta}\otimes I_{H})UW^{\ast}=\left(
\begin{array}
[c]{cccccc}%
0 & \cdots & W_{k}\eta & 0 & \cdots & \\
0 & 0 & \cdots & W_{k}\eta & 0 & \cdots\\
& 0 & 0 & \ddots & W_{k}\eta & \ddots\\
&  & 0 & \ddots &  & \ddots\\
&  &  & \ddots &  & \\
&  &  &  &  &
\end{array}
\right)  \text{.}%
\]
Further, $WU^{\ast}(H^{\infty}(E^{\sigma})\otimes I_{H})UW^{\ast}=\{R\in
B(\ell^{2}(\mathbb{Z}_{+},H))\mid R$ satisfies equations (\ref{RToep}) and
(\ref{Rcomm})$\}$, which is the commutant of $S_{H}\times\psi_{H}(H^{\infty
}(E))$.
\end{proposition}

Suppose now that $\pi$ is a completely contractive $\sigma$-weakly continuous
representation of any Hardy algebra, $H^{\infty}(E)$, on a Hilbert space $H$,
then the associated covariant representation $(T,\sigma)$ of $E$ is given by
the formulae $\sigma=\pi\circ\varphi_{\infty}$ and $T(\xi)=\pi(T_{\xi})$,
$\xi\in E$. Consequently, in the present setting, where the Hardy algebra is
$M\rtimes_{\alpha}\mathbb{Z}_{+}$, if $\pi$ is a completely contractive
$\sigma$-weakly continuous representation of $M\rtimes_{\alpha}\mathbb{Z}_{+}$
on the Hilbert space $H$, the covariant representation $(T,\sigma)$ of
$_{\alpha}M$ on $H$ is determined entirely by $\sigma$ and the contraction
operator $t:=T(1)=\pi(w)$. If we let $W_{k}:E^{\otimes k}\otimes_{\sigma
}H\rightarrow H$ be the Hilbert space isomorphism from equation
(\ref{inducediso}) and compute, we find that
\begin{multline*}
\tilde{T}_{k}W_{k}^{\ast}\sigma(\alpha^{k-1}(\xi_{1})\alpha^{k-2}(\xi
_{2})\cdots\xi_{k}))h=\tilde{T}_{k}(\xi_{1}\otimes\xi_{2}\cdots\otimes\xi
_{k}\otimes h)\\
=T(\xi_{1})T(\xi_{2})\cdots T(\xi_{k})h=T(1)\sigma(\xi_{1})T(1)\sigma(\xi
_{2})\cdots T(1)\sigma(\xi_{k})h\\
=t^{k}\sigma(\alpha^{k-1}(\xi_{1})\alpha^{k-2}(\xi_{2})\cdots\xi_{k})h.
\end{multline*}
(In the last equality we used the fact that $t=T(1)$ and the covariance
property of the representation). Thus the generalized powers of $\tilde{T}$
are related to the ordinary powers of $t$ through the equation $\tilde{T}%
_{k}W_{k}^{\ast}=t^{k}$ for $k\geq1$. In particular, we see that $\Vert
\tilde{T}_{k}^{\ast}h\Vert=\Vert t^{k\ast}h\Vert$ for all $h\in H$. It follows
that $(T,\sigma)$ is a $C_{\cdot0}$-representation or a c.n.c. representation
if and only if $t$ is a $C_{\cdot0}$-operator or a completely non-coisometric operator.

Also, the defect operators of $(T,\sigma)$ are related to the defect operators
of $t$ via the formulae $(I_{H}-\tilde{T}\tilde{T}^{\ast})^{1/2}%
=(I_{H}-tt^{\ast})^{1/2}$ and $W_{1}(I_{E\otimes_{\sigma}H}-\tilde{T}^{\ast
}\tilde{T})^{1/2}W_{1}^{\ast}=(I_{H}-t^{\ast}t)^{1/2}$. Hence, if we form
$\tau_{1}:=\sigma\circ\alpha|\mathcal{D}$ where, as usual, $\mathcal{D}%
=\overline{(I_{E\otimes_{\sigma}H}-\tilde{T}^{\ast}\tilde{T})^{1/2}%
(E\otimes_{\sigma}H)}$, and if we form $W^{\tau_{1}}:\mathcal{F}%
(E)\otimes_{\tau_{1}}\mathcal{D}\rightarrow\ell^{2}(\mathbb{Z}_{+}%
,\mathcal{D})$ and follow it with $I\otimes W_{1}$ mapping $\ell
^{2}(\mathbb{Z}_{+},\mathcal{D})$ onto $\ell^{2}(\mathbb{Z}_{+},\mathcal{D}%
_{t})$, where $\mathcal{D}_{t}=\overline{(I_{H}-t^{\ast}t)^{1/2}H}$ is the
defect space of $t$, then $\mathcal{W}_{v}:=I_{H}\oplus(I\otimes W_{1}%
)W^{\tau_{1}}$ is a Hilbert space isomorphism mapping the Hilbert space of the
minimal isometric dilation $(V,\rho)$ of $(T,\sigma)$ onto the Hilbert space
of the minimal isometric dilation of $t$, vis., $H\oplus\ell^{2}%
(\mathbb{Z}_{+},\mathcal{D}_{t})$. Further, we have $\mathcal{W}%
_{v}V(1)\mathcal{W}_{v}^{\ast}=v$, where
\[
v=\left(
\begin{array}
[c]{clll}%
t & 0 & 0 & \cdots\\
d & 0 & 0 & \cdots\\
0 & I_{\mathcal{D}_{t}} & 0 & \\
0 & 0 & I_{\mathcal{D}_{t}} & \\
&  &  & \ddots\\
&  &  &
\end{array}
\right)  \text{,}%
\]
and $d:=(I_{H}-t^{\ast}t)^{1/2}$.

Now consider the characteristic operator of $(T,\sigma)$, $(\Theta
_{T},\mathcal{D},\mathcal{D}_{\ast},\tau_{1},\tau_{2})$ and identify
$(T,\sigma)$ with its canonical model using Theorem \ref{TSisT}. Recall from
Remark \ref{CharOpBackground} that our notation remains consistent; this new
$\tau_{1}$ is still the restriction of $\sigma\circ\alpha$ to $\mathcal{D}$;
$\tau_{2}$ is the restriction of $\sigma$ to $\mathcal{D}_{\ast}$. Even though
the defect space $\mathcal{D}_{\ast}$ for $(T,\sigma)$ is the same as the
defect space $\mathcal{D}_{\ast t}:=\overline{(I-tt^{\ast})^{1/2}H}$, we shall
continue to distinguish notationally between them. Thus $t=T(1)$ is the
operator which, in the notation of Theorem \ref{TSisT}, would be denoted
$T_{\hat{\Theta}_{T}}(1)$ and similarly the minimal isometric dilation
$(V,\rho)$ of $(T,\sigma)$ would be denoted $(V_{\hat{\Theta}_{T}},\rho
_{\hat{\Theta}_{T}})$, etc. However, this notation is ponderous and so we
shall drop the subscript $\hat{\Theta}_{T}$. We shall write $\mathcal{W}%
_{\ast}$ for $\mathcal{W}^{\tau_{2}}$, so that $\mathcal{W}_{\ast}$ is a
Hilbert space isomorphism from $\mathcal{F}(E)\otimes_{\tau_{2}}%
\mathcal{D}_{\ast}$ onto $\ell^{2}(\mathbb{Z}_{+},\mathcal{D}_{\ast t})$ such
that
\[
\mathcal{W}_{\ast}(w\otimes I_{\mathcal{D}_{\ast}})=S_{\mathcal{D}_{\ast t}%
}\mathcal{W}_{\ast}%
\]
where $S_{\mathcal{D}_{\ast t}}$ is the unilateral shift on $\ell
^{2}(\mathbb{Z}_{+},\mathcal{D}_{\ast t})$. We also write $\mathcal{W}_{1}$
for $(I\otimes W_{1}^{\tau_{1}})W^{\tau_{1}}$, which is a Hilbert space
isomorphism from $\mathcal{F}(E)\otimes_{\tau_{1}}\mathcal{D}$ onto $\ell
^{2}(\mathbb{Z}_{+},\mathcal{D}_{t})$ that satisfies the equation%
\[
\mathcal{W}_{1}(w\otimes I_{\mathcal{D}})=S_{\mathcal{D}_{t}}\mathcal{W}%
_{1}\text{,}%
\]
where for $S_{\mathcal{D}_{t}}$ is the unilateral shift on $\ell
^{2}(\mathbb{Z}_{+},\mathcal{D}_{t})$. The characteristic operator $\Theta
_{T}$ maps $\mathcal{F}(E)\otimes_{\tau_{1}}\mathcal{D}$ to $\mathcal{F}%
(E)\otimes_{\tau_{2}}\mathcal{D}_{\ast}$ and intertwines the induced
representations, $\tau_{1}^{\mathcal{F}(E)}$ and $\tau_{2}^{\mathcal{F}(E)} $.
Thus, if we set $\Theta:=\mathcal{W}_{\ast}\Theta_{T}\mathcal{W}_{1}^{-1}$, we
obtain a contraction from $\ell^{2}(\mathbb{Z}_{+},\mathcal{D}_{t})$ to
$\ell^{2}(\mathbb{Z}_{+},\mathcal{D}_{\ast t})$ that intertwines
$S_{\mathcal{D}_{\ast t}}$ and $S_{\mathcal{D}_{t}}$. We shall write
$\Delta_{T}$ for $(I-\Theta_{T}^{\ast}\Theta_{T})^{1/2}$ and $\Delta$ for
$(I-\Theta^{\ast}\Theta)^{1/2}$, so $\mathcal{W}_{1}\Delta_{T}\mathcal{W}%
_{1}^{-1}=\Delta$. Also, we shall write $\mathcal{W}_{\Delta}$ for the
restriction of $\mathcal{W}_{1}$ to $\overline{\Delta_{T}(\mathcal{F}%
(E)\otimes_{\tau_{1}}\mathcal{D})}$, obtaining a Hilbert space isomorphism
from this space onto $\overline{\Delta\ell^{2}(\mathbb{Z}_{+},\mathcal{D}%
_{t})}$. Consequently, $\mathcal{W}$ which we shall define to be
$\mathcal{W}_{\ast}\oplus\mathcal{W}_{\Delta}$ is a Hilbert space isomorphism
from $K(\Theta_{T})$, which recall from Theorem \ref{TSisT} is $(\mathcal{F}%
(E)\otimes_{\tau_{2}}\mathcal{D}_{\ast})\oplus\overline{\Delta_{T}%
(\mathcal{F}(E)\otimes_{\tau_{1}}\mathcal{D})}$, onto $\ell^{2}(\mathbb{Z}%
_{+},\mathcal{D}_{\ast t})\oplus\overline{\Delta\ell^{2}(\mathbb{Z}%
_{+},\mathcal{D}_{t})}$.

Recall next the definition of $S_{\hat{\Theta}_{T}}(\cdot):=S$, from Lemma
\ref{charep}, and write $S$ for the isometry $S(1)$. (Actually, $S$ is unitary
as we shall see in a moment.) Then if $\tilde{S}$ is defined on $\overline
{\Delta\ell^{2}(\mathbb{Z}_{+},\mathcal{D}_{t})}$ by the formula $\tilde
{S}(\Delta\xi)=\Delta S_{\mathcal{D}_{t}}\xi$, then, as an easy calculation
shows, $\mathcal{W}_{\Delta}$ implements a unitary equivalence between $S$ and
$\tilde{S}$. Consequently, $\mathcal{W}$ implements a unitary equivalence
between $S_{\mathcal{D}_{\ast t}}\oplus\tilde{S}$ acting on $\ell
^{2}(\mathbb{Z}_{+},\mathcal{D}_{\ast t})\oplus\overline{\Delta\ell
^{2}(\mathbb{Z}_{+},\mathcal{D}_{t})}$. Thus, it looks like $\mathcal{W}$
implements a unitary equivalence between the minimal isometric dilation
$v=V(1)$ for $t$ and the isometry that occurs in the Sz.-Nagy-Foia\c{s} model
for $t$ in \cite{szNF70}.\footnote{Strictly speaking to identify fully the
constructs of the Sz.-Nagy-Foia\c{s} theory, we need to transfer the
discussion from $\ell^{2}$-spaces on $\mathbb{Z}$ to $L^{2}$-spaces on
$\mathbb{T}$ via the Fourier transform. We omit this detail. However, the
whole theory has been developed on $\mathbb{Z}$ by Douglas in \cite{rD68}.}
But $S_{\mathcal{D}_{\ast t}}\oplus\tilde{S}$ is not quite the
Sz.-Nagy-Foia\c{s} model isometry. The point is that the model that Sz.-Nagy
and Foia\c{s} produce acts on $\ell^{2}(\mathbb{Z}_{+},\mathcal{D}_{\ast
t})\oplus\overline{\tilde{\Delta}\ell^{2}(\mathbb{Z},\mathcal{D}_{t})}$, where
$\ell^{2}(\mathbb{Z},\mathcal{D}_{t})$ consists of all square summable
$\mathcal{D}_{t}$-valued functions on the \emph{integers} $\mathbb{Z}$,
$\tilde{\Delta}$ is an operator that we describe in a second and the part of
the model that acts on $\overline{\tilde{\Delta}\ell^{2}(\mathbb{Z}%
,\mathcal{D}_{t})}$ is the (restriction of the) bilateral shift. The
difference lies in the definition of $\tilde{\Delta}$. Note that since
$\Theta$ intertwines $S_{\mathcal{D}_{\ast t}}$ and $S_{\mathcal{D}_{t}}$,
$\Theta$ has a \emph{unique} extension to an operator $\tilde{\Theta}$ from
$\ell^{2}(\mathbb{Z},\mathcal{D}_{t})$ to $\ell^{2}(\mathbb{Z},\mathcal{D}%
_{\ast t})$ that intertwines the two \emph{bilateral }shifts. We simply let
$\tilde{\Delta}=(I-\tilde{\Theta}^{\ast}\tilde{\Theta})^{1/2}$. Then the piece
Sz.-Nagy and Foia\c{s} build for their model is $\overline{\tilde{\Delta}%
\ell^{2}(\mathbb{Z},\mathcal{D}_{t})}$. However, in terms of $\tilde{\Theta}$,
$\Delta=(I-P\tilde{\Theta}^{\ast}P\tilde{\Theta})^{1/2}|_{\ell^{2}%
(\mathbb{Z}_{+},\mathcal{D}_{t})}$, so on the face of it, one would expect
$\overline{\Delta\ell^{2}(\mathbb{Z}_{+},\mathcal{D}_{t})}$ and $\overline
{\tilde{\Delta}\ell^{2}(\mathbb{Z},\mathcal{D}_{t})}$ to be different.
Neverthless, if we assume that our representation $(T,\sigma)$ is c.n.c., as
we shall, then the map that takes $\ell^{2}(\mathbb{Z}_{+},\mathcal{D}_{t})$
to $\overline{\tilde{\Delta}\ell^{2}(\mathbb{Z},\mathcal{D}_{t})}$ by sending
a vector of the form $\Delta\xi$ to $\tilde{\Delta}\tilde{\xi}$, where
$\tilde{\xi}$ is the extension of $\xi$ to all of $\mathbb{Z}$, which is zero
on the negative integers, is in fact a Hilbert space isomorphism that
intertwines $\tilde{S}$ on $\ell^{2}(\mathbb{Z}_{+},\mathcal{D}_{t})$ and the
restriction of the bilateral shift to $\overline{\tilde{\Delta}\ell
^{2}(\mathbb{Z},\mathcal{D}_{t})}$. This is the content, really, of part (ii)
of Lemma \ref{cond}, which gives meaning to the term \textquotedblleft
predictable\textquotedblright. Thus, if we incorporate this additional Hilbert
space isomorphism ($\Delta\xi\mapsto\tilde{\Delta}\tilde{\xi})$ into the
definition of $\mathcal{W}$, then we have proved most of the following
theorem. The remaining details are easy to supply and so will be omitted.

\begin{theorem}
\label{Malpha}Let $\pi$ be a completely contractive, $\sigma$-weakly
continuous representation of the analytic crossed product $M\rtimes_{\alpha
}\mathbb{Z}_{+}$ on a Hilbert space $H$ such that $t=\pi(w)=T(1)$ is a c.n.c.
contraction, where $(T,\sigma)$ is the associated covariant representation,
and let $(\Theta_{T},\mathcal{D},\mathcal{D}_{\ast},\tau_{1},\tau_{2})$ be the
characteristic operator attached to this representation. Then the Hilbert
space isomorphism $\mathcal{W}$ just described, viewed as a map from the space
$K(\hat{\Theta}_{T})$ of the minimal isometric dilation of $(T,\sigma)$ to the
shift space $\ell^{2}(\mathbb{Z}_{+},\mathcal{D}_{\ast t})\oplus
\overline{\tilde{\Delta}\ell^{2}(\mathbb{Z},\mathcal{D}_{t})}$ maps all parts
of the model space for $(T,\sigma)$ to the corresponding parts of
Sz.-Nagy-Foia\c{s} model space for $t$, i.e., the operator $\Theta
=\mathcal{W}_{\ast}\Theta_{T}\mathcal{W}^{-1}$ described above is equivalent
to the characteristic operator function of the operator $t$ described in
\cite{szNF70}.
\end{theorem}

\begin{conclremarks}
\label{Coclremarks}\hfill

\begin{enumerate}
\item[(i)] In view of Theorem \ref{Malpha}, it appears that for analytic
crossed products, at least, one may extend the model developed in Theorem
\ref{TSisT} to get a \emph{unitary }dilation for a c.n.c. representation
$(T,\sigma)$ of the algebra. That is, thinking of the isometric dilation
$(V,\rho)$ for $(T,\sigma)$ as acting on $\ell^{2}(\mathbb{Z}_{+}%
,\mathcal{D}_{\ast t})\oplus\overline{\tilde{\Delta}\ell^{2}(\mathbb{Z}%
,\mathcal{D}_{t})}$, $v:=V(1)$ is an isometry that satisfies the equation
$v\rho\circ\alpha(a)=\rho(a)v$ for all $a\in M$. The minimal unitary
\emph{extension} of $v$ is the (restriction of the) bilateral shift acting
$\ell^{2}(\mathbb{Z},\mathcal{D}_{\ast t})\oplus\overline{\tilde{\Delta}%
\ell^{2}(\mathbb{Z},\mathcal{D}_{t})}$. However, while $v$ extends to a
unitary $w$, say, on $\ell^{2}(\mathbb{Z},\mathcal{D}_{\ast t})\oplus
\overline{\tilde{\Delta}\ell^{2}(\mathbb{Z},\mathcal{D}_{t})}$, it may not be
possible to extend $\rho$ to a representation $\tilde{\rho}$ on this space so
that the equation $w\tilde{\rho}\circ\alpha(a)=\tilde{\rho}(a)w$ also holds
for all $a\in M$. If such a $\tilde{\rho}$ were to exist, then it would have a
natural extension to the\ $C^{\ast}$-inductive limit of the system built from
$M$ and the powers of $\alpha$ as described in \cite{pS93}. Simple examples
show that this need not be the case. We intend to take this matter up in a
future study.

\item[(ii)] The example studied in this section may seem very special.
However, thanks to our investigation in \cite{MS00}, we may assert that under
technical conditions that we ignore here, every $W^{\ast}$-correspondence over
a von Neumann algebra is Morita equivalent to one that comes from an
endomorphism of another, possibly different, von Neumann algebra. Thus, up to
Morita equivalence, all Hardy algebras are analytic crossed products. We
intend take this matter up also in a future study.

\item[(iii)] As we noted in Theorem \ref{Malpha}, the characteristic function
$\hat{\Theta}_{T}$ of the representation $(T,\sigma)$ is equivalent to the
characteristic operator function $\Theta$ of $t=T(1)$ (after one takes the
Fourier transform that identifies $\ell^{2}$ with $L^{2}(\mathbb{T)}$ and
identifies $\Theta$ as a function, rather than as an operator.). Classically,
$\Theta$ is an analytic function from the open unit disc $\mathbb{D}$ in
$\mathbb{C}$ to $B(\mathcal{D}_{t},\mathcal{D}_{\ast t})$. On the other hand,
because $\hat{\Theta}_{T}$ is an element of $H^{\infty}(E^{\tau})$, where
$(\mathcal{G},\tau)$ is the supplement of $\tau_{1}$ and $\tau_{2}$ that we
fixed in the discussion just before equation (\ref{DecompFEG}), $\hat{\Theta
}_{T}$ has a Taylor or Fourier expansion%
\[
\hat{\Theta}_{T}\sim T_{\eta_{0}}+T_{\eta_{1}}+\cdots\text{,}%
\]
where the $\eta_{i}\in(E^{\tau})^{\otimes i}$. As we show in \cite{MS04} using
the gauge group, the arithmetic means of this series converge weak-$\ast$ to
$\hat{\Theta}_{T}$. As we noted above, $W_{i}\cdot(E^{\sigma})^{\otimes
i}:=\{W_{i}\eta\mid\eta\in(E^{\sigma})^{\otimes i}\}=\{z\in B(\mathcal{G})\mid
z\sigma(a)=\sigma(\alpha^{i}(a))z,\ a\in M\}$. To compute the $W_{i}\eta
_{i}\in B(\mathcal{G})$, we may appeal to the analysis leading to Theorem
\ref{expression} or to the result of the calculation there to conclude that
$W_{0}\eta_{0}=-t|\mathcal{D}$, $W_{1}\eta_{1}=\Delta_{\ast}\Delta
|\mathcal{D}$, $W_{2}\eta_{2}=\Delta_{\ast}t^{\ast}\Delta|\mathcal{D}$,
$\cdots$. So, if we evaluate $\Theta_{T}$ on the open unit ball of
$E=~_{\alpha}M$ using the formula from Theorem \ref{expression}, then a
straightforward calculation based on the analysis we have made and the
definition of the characteristic operator function for $t$ from \cite{szNF70}
shows that if $\xi_{0}$ denotes the identity operator in $M$, but viewed as a
vector in $E$, then for all complex numbers $z$, $|z|<1$,
\[
\Theta(\overline{z})=\hat{\Theta}_{T}(z\xi_{0})\text{.}%
\]
(The reason for $\overline{z}$ and not $z$ is an artifact of the role that
elements in the dual play in the representations of the algebras and need not
concern us here.) Thus, $\hat{\Theta}_{T}$ is effectively determined on the
one dimensional slice $\{z\xi_{0}\mid|z|<1\}$. Of course, this is fairly
evident from Theorem \ref{expression} and the fact that $\xi_{0}$ is a cyclic
vector for $E$ as a right module over $M$.
\end{enumerate}
\end{conclremarks}

\end{document}